%% file: sm.tex
\documentclass{svmult}

\usepackage{makeidx}         
\usepackage{graphicx}        
\usepackage{multicol}        
\usepackage[bottom]{footmisc}

\makeindex             

\usepackage[all]{xy}
\usepackage{amssymb}
\usepackage{amsmath}

\newcommand{\id}{{1\hspace{-1.2mm}1}}
 \newcommand{\lon}{\longrightarrow}
 
 \newcommand{\rar}{\rightarrow}

 \newcommand{\End}{{\mathcal E n d}}
\newcommand{\p}{{\partial}}
\newcommand{\Id}{{\mathrm I\mathrm d}}
\newcommand{\no}{{\noindent}}

 \newcommand{\Z}{{\mathbb Z}}
 \newcommand{\bS}{{\mathbb S}}

 \newcommand{\R}{{\mathbb R}}
 
 \newcommand{\K}{{\mathbb K}}
 \newcommand{\ot}{\otimes}

\newcommand{\sV}{{\mathsf V}}

\newcommand{\out}{{\mathsf  O\mathsf u\mathsf t}}
\newcommand{\inn}{{\mathsf  I\mathsf n}}
\newcommand{\sgn}{{\mathit s  \mathit g\mathit  n}}

\newcommand{\DefQ}{{\mathcal D \mathit e\mathit f\mathcal Q}}
\newcommand{\PV}{{\mathcal P \mathit o\mathit l\mathit y \mathcal V}}

 \newcommand{\Beq}{\begin{equation}}
 \newcommand{\Eeq}{\end{equation}}
 \newcommand{\Beqr}{\begin{eqnarray}}
 \newcommand{\Eeqr}{\end{eqnarray}}
 \newcommand{\Beqrn}{\begin{eqnarray*}}
 \newcommand{\Eeqrn}{\end{eqnarray*}}
 \newcommand{\Ba}{\begin{array}}
 \newcommand{\Ea}{\end{array}}
 \newcommand{\Bi}{\begin{itemize}}
 \newcommand{\Ei}{\end{itemize}}
 \newcommand{\Bc}{\begin{center}}
 \newcommand{\Ec}{\end{center}}

 \newcommand{\fG}{{\mathfrak G}}

 \newcommand{\f}{{\mathcal O}}
 \newcommand{\cA}{{\mathcal A}}
 \newcommand{\cB}{{\mathcal B}}
 \newcommand{\cC}{{\mathcal C}}
 \newcommand{\caD}{{\mathcal D}}
 \newcommand{\cE}{{\mathcal E}}
 \newcommand{\cF}{{\mathcal F}}

 \newcommand{\cP}{{\mathcal P}}
\newcommand{\cQ}{{\mathcal Q}}
 \newcommand{\cS}{{\mathcal S}}
 \newcommand{\cT}{{\mathcal T}}

 \newcommand{\al}{\alpha}
 
 \newcommand{\ga}{\gamma}
 
 \newcommand{\Ga}{\Gamma}
 
 \newcommand{\var}{\varepsilon}

 \newcommand{\Hom}{{\mathrm H\mathrm o\mathrm m}}
 
\newcommand{\diff}{{\mathit d \mathit i \mathit f \mathit f}}
 \newcommand{\sip}{\smallskip}
 \newcommand{\bip}{\bigskip}
 
\begin{document}

\title*{Permutahedra, HKR isomorphism and polydifferential Gerstenhaber-Schack complex}
\titlerunning{Permutahedra, HKR isomorphism and  Gerstenhaber-Schack complex}
\author{S.A.\ Merkulov}
\institute{Department of Mathematics, Stockholm University, 10691 Stockholm, Sweden
\texttt{sm@math.su.se}}
%
%
\maketitle

\vspace{-3mm}
\begin{center}
{\em
To Murray Gerstenhaber and Jim Stasheff
}
\end{center}
\vspace{-5mm}

\hyphenation{pro-perads}

\section{Introduction}


This paper aims to give a short but self-contained introduction into the theory of (wheeled) props,
properads, dioperads and operads, and illustrate some of its key ideas in terms of a prop(erad)ic interpretation of
simplicial and permutahedra cell complexes  with subsequent  applications to the Hochschild-Kostant-Rosenberg
type isomorphisms.

Let $V$ be a graded vector space over a field $\K$ and $\f_V:=\odot^\bullet V^*$ the free graded commutative algebra generated by
the dual vector space $V^*:=\Hom_\K(V,\K)$. One can interpret $\f_V$ as  the algebra of polynomial functions on the  space $V$.
The classical
Hochschild-Kostant-Rosenberg theorem asserts that the Hochschild cohomology of $\f_V$ (with coefficients in $\f_V$)
is isomorphic to the space,  $\wedge^\bullet \cT_V$, of polynomial polyvector fields on $V$ which in turn is isomorphic
as a vector space to $\wedge^\bullet V\ot \odot^\bullet V^*$,
\Beq\label{introduction-KTH}
HC^\bullet(\f_V)\simeq \wedge^\bullet \cT_V\simeq  \wedge^\bullet V\ot \odot^\bullet V^*.
\Eeq
The  Hochschild complex,
$
C^\bullet(\f_V)=\oplus_{k\geq 0} \Hom(\f_V^{\ot k}, \f_V)[1-k],
$
of $\f_V$ has a natural subcomplex,
$C^\bullet_{\diff}(\f_V)\subset C^\bullet(\f_V)$, spanned by polydifferential operators. It was proven
in \cite{Ko} (see also \cite{CFL}) that again
\vspace{-1mm}
\Beq\label{introduction-KTH-poly}
HC^\bullet_{\diff}(\f_V)\simeq  \wedge^\bullet \cT_V \simeq \wedge^\bullet V\ot \odot^\bullet V^*.
\vspace{-1mm}
\Eeq
The first result (\ref{introduction-KTH}) actually fails  for a ring of smooth functions, $\f_M$, on a generic graded manifold
$M$ while the second one (\ref{introduction-KTH-poly}) stays always  true \cite{CFL}. Thus one must, in general, be careful in distinguishing ordinary and polydifferential
Hochschild cohomology for smooth functions.

The vector space $\f_V$  has  a natural (co)commutative bialgebra structure so that one can also associate
to $\f_V$  a Gerstenhaber-Schack complex \cite{GS},
$
\cC^{\bullet,\bullet}(\f_V):=\bigoplus_{m,n\geq 1} \Hom(\f_V^{\ot m}, \f_V^{\ot n})[2-m-n]
$.
Is cohomology was computed in \cite{GS2,LM},
\Beq\label{introduction-KTH-GS}
H\cC^{\bullet,\bullet}(\f_V)\simeq  \wedge^{\bullet \geq 1} V\ot \wedge^{\bullet\geq 1} V^*.
\Eeq
In this paper we introduce  (more precisely, {\em deduce}\, from the permutahedra cell complex)  a
relatively non-obvious {\em polydifferential}\, subcomplex,
 $\cC^{\bullet,\bullet}_{\diff}(\f_V)\subset \cC^{\bullet,\bullet}(\f_V)$, such that
$
\cC^{\bullet,\bullet}_{\diff}(\f_V)\cap C^\bullet(\f_V)= C^\bullet_{diff}(\f_V)
$,
and prove that this inclusion is a quasi-isomorphism,
\Beq\label{introduction-KTH-GS-poly}
H\cC^{\bullet,\bullet}_{\diff}(\f_V)\simeq  \wedge^{\bullet\geq 1} V\ot \wedge^{\bullet\geq 1} V^*.
\Eeq
In fact, we show in this paper very simple pictorial proofs of all four results,
(\ref{introduction-KTH})-(\ref{introduction-KTH-GS-poly}), mentioned above:
first we interpret Saneblidze-Umble's  \cite{SU} permutahedra cell complex as a differential graded (dg, for short) properad $\cP$,
then we use tensor powers of $\cP$ to create a couple of other dg props, $\caD$  and $\cQ$, whose
cohomology we know immediately by their  very constructions, and then, studying representations of $\caD$ and $\cQ$
 in an arbitrary vector space
$V$ we {\em obtain}\, (rather than {\em define}) the well-known polydifferential subcomplex
of the Hochschild complex for $\f_V$ and, respectively, a new  polydifferential subcomplex of the Gerstenhaber-Schack
complex whose cohomologies are given, in view of contractibility of the permutahedra,  by formulae (\ref{introduction-KTH-poly}) and (\ref{introduction-KTH-GS-poly}). Finally, using again the language of props
we deduce from (\ref{introduction-KTH-poly}) and (\ref{introduction-KTH-GS-poly})  formulae
(\ref{introduction-KTH}) and, respectively, (\ref{introduction-KTH-GS}). As a corollary to (\ref{introduction-KTH-GS-poly}) and (\ref{introduction-KTH-GS}) we show a slight sharpening of the famous Etingof-Kazhdan theorem:
{\em for any Lie bialgebra structure on a vector space $V$ there exists its bialgebra quantization within the class
of polydifferential operators from $\cC^{\bullet,\bullet}_{poly}(\f_V)$.}

The paper is organized as follows. In \S 2 we give a short but self-contained introduction into the
theory of (wheeled) props, properads, dioperads and operads.
In \S 3 we prove formulae
(1)-(4) using properadic interpretation of the permutahedra cell complex. In \S 4 we study a dg prop, $\DefQ$,
whose representations in a dg space $V$ are in one-to-one correspondence with unital $A_\infty$-structures
on $\f_V$, and use it to give a new pictorial proof of another classical result that isomorphisms (\ref{introduction-KTH}) and
(\ref{introduction-KTH-poly}) extend to isomorphisms of Lie algebras, with $\wedge^\bullet \cT_M$ assumed to be equipped with
Schouten  brackets.

We work over a field $\K$ of characteristic zero.  If $V=\oplus_{i\in \Z} V^i$ is a graded vector space, then
$V[k]$ is a graded vector space with $V[k]^i:=V^{i+k}$. We denote $\ot^\bullet V:= \bigoplus_{n\geq 0}\ot^n V$,
$\ot^{\bullet\geq 1} V:= \bigoplus_{n\geq 1}\ot^n V$, and similarly for symmetric and skewsymmetric tensor powers,
$\odot^\bullet V$ and $\wedge^\bullet V$. The symbol $[n]$ stands for an ordered set $\{1,2,\ldots,n\}$.


\section{An introduction to the theory of (wheeled) props}
\label{ch1}


\noindent{\bf 2.1\, An associative algebra as a morphism of graphs.}
Recall that  an {\em associative algebra}\, structure on a vector space $E$
is a linear map $E\otimes E \rightarrow E$ satisfying the associativity condition,
$(a_1a_2)a_3 = a_1(a_2a_3)$, for any $a_1,a_2,a_2\in E$. Let us represent a typical element,
$a_1\ot a_2\ot \ldots \ot a_n\in  \ot^n E$, of  the tensor algebra, $\ot^\bullet E$,
of $E$
as a {\em decorated directed graph},
$$
G\langle a_1, \ldots a_n\rangle:= \ \ \
\begin{xy}
<0mm,0mm>*{\cdot};
<0mm,1mm>*{\cdot};
<0mm,-1mm>*{\cdot};
<0mm,-2mm>*{};<0mm,-5mm>*{}**@{-},
 <0mm,-5mm>*{\bullet};
<0mm,-5mm>*{};<0mm,-8mm>*{}**@{-},
%
<0mm,2mm>*{};<0mm,5mm>*{}**@{-},
 <0mm,5mm>*{\bullet};
<0mm,5mm>*{};<0mm,8mm>*{}**@{-},
 <0mm,9mm>*{\bullet};
<0mm,9mm>*{};<0mm,12mm>*{}**@{-},
<4mm,9mm>*{a_1};
<4mm,5mm>*{a_2};
<5.5mm,-5mm>*{a_{n}};
\end{xy}, \vspace{-1mm}
$$
where the adjective {\em decorated}\, means that each vertex of the shown graph $G$ is equipped  with an  element
of $E$ and the adjective {\em directed}\, means that the graph $G$ is equipped with the flow running by default
(unless otherwise is explicitly shown)
from the bottom to the top. 
Let $G\langle E\rangle$ be the vector space spanned
by all such decorations, $G\langle a_1, \ldots a_n\rangle$, of the shown chain-like graph $G$ modulo 
the relations of the type,  \vspace{-0.5mm}
$$
\begin{xy}
<0mm,0mm>*{\cdot};
<0mm,1mm>*{\cdot};
<0mm,-1mm>*{\cdot};
<0mm,-2mm>*{};<0mm,-5mm>*{}**@{-},
 <0mm,-5mm>*{\bullet};
<0mm,-5mm>*{};<0mm,-8mm>*{}**@{-},
%
<0mm,2mm>*{};<0mm,5mm>*{}**@{-},
 <0mm,5mm>*{\bullet};
<0mm,5mm>*{};<0mm,8mm>*{}**@{-},
 <0mm,9mm>*{\bullet};
<0mm,9mm>*{};<0mm,12mm>*{}**@{-},
<10.6mm,5mm>*{\lambda_1 a_2'+ \lambda_2 a_2''};
<4mm,9mm>*{a_1};
<5.5mm,-5mm>*{a_{n}};
\end{xy}
=
\lambda_1\
\begin{xy}
<0mm,0mm>*{\cdot};
<0mm,1mm>*{\cdot};
<0mm,-1mm>*{\cdot};
<0mm,-2mm>*{};<0mm,-5mm>*{}**@{-},
 <0mm,-5mm>*{\bullet};
<0mm,-5mm>*{};<0mm,-8mm>*{}**@{-},
%
<0mm,2mm>*{};<0mm,5mm>*{}**@{-},
 <0mm,5mm>*{\bullet};
<0mm,5mm>*{};<0mm,8mm>*{}**@{-},
 <0mm,9mm>*{\bullet};
<0mm,9mm>*{};<0mm,12mm>*{}**@{-},
<4mm,9mm>*{a_1};
<4mm,5mm>*{a_2'};
<5.5mm,-5mm>*{a_{n}};
\end{xy}
\ + \
\lambda_2 \
\begin{xy}
<0mm,0mm>*{\cdot};
<0mm,1mm>*{\cdot};
<0mm,-1mm>*{\cdot};
<0mm,-2mm>*{};<0mm,-5mm>*{}**@{-},
 <0mm,-5mm>*{\bullet};
<0mm,-5mm>*{};<0mm,-8mm>*{}**@{-},
%
<0mm,2mm>*{};<0mm,5mm>*{}**@{-},
 <0mm,5mm>*{\bullet};
<0mm,5mm>*{};<0mm,8mm>*{}**@{-},
 <0mm,9mm>*{\bullet};
<0mm,9mm>*{};<0mm,12mm>*{}**@{-},
<4mm,9mm>*{a_1};
<4mm,5mm>*{a_2''};
<5.5mm,-5mm>*{a_{n}};
\end{xy}
\ \ \  \forall \lambda_1,\lambda_2\in \K,
 \vspace{-1mm}
$$
which identify $G\langle E\rangle$ with $\ot^n E$. Note that if $G$ has only one internal vertex
(we call such graphs (1,1)-{\em corollas}\,), then
$G\langle E \rangle=E$.
The multiplication operation in $E$ gets encoded
in this picture as a contraction of an internal edge, e.g.  \vspace{-1mm}
$$
\begin{xy}
<0mm,0mm>*{\cdot};
<0mm,1mm>*{\cdot};
<0mm,-1mm>*{\cdot};
<0mm,-2mm>*{};<0mm,-5mm>*{}**@{-},
 <0mm,-5mm>*{\bullet};
<0mm,-5mm>*{};<0mm,-8mm>*{}**@{-},
%
<0mm,2mm>*{};<0mm,5mm>*{}**@{-},
 <0mm,5mm>*{\bullet};
<0mm,5mm>*{};<0mm,8mm>*{}**@{-},
 <0mm,9mm>*{\bullet};
<0mm,9mm>*{};<0mm,12mm>*{}**@{-},
<4mm,9mm>*{a_1};
<4mm,5mm>*{a_2};
<5.5mm,-5mm>*{a_{n}};
\end{xy} \ \
\lon \ \ \
\begin{xy}
<0mm,0mm>*{\cdot};
<0mm,1mm>*{\cdot};
<0mm,-1mm>*{\cdot};
<0mm,-2mm>*{};<0mm,-5mm>*{}**@{-},
 <0mm,-5mm>*{\bullet};
<0mm,-5mm>*{};<0mm,-8mm>*{}**@{-},
%
<0mm,2mm>*{};<0mm,5mm>*{}**@{-},
 <0mm,5mm>*{\bullet};
<0mm,5mm>*{};<0mm,8mm>*{}**@{-},
<5.5mm,5mm>*{a_1a_2};
<5.5mm,-5mm>*{a_{n}};
\end{xy}
\vspace{-0.5mm}
$$
which upon repetition gives a {\em contraction}\,  map
$
\mu_G: G\langle E\rangle \rar E.
$
Moreover, the associativity conditions for the multiplication assures us that the map $\mu_G$ is canonical,
i.e.\ it does {\em not}\, depend on a particular sequence of contractions of the graph $G$ into a corolla and is uniquely
determined by the graph $G$ itself.

\sip

Actually there is no need to be specific about contracting  precisely  {\em two}\, vertices --- any connected subset of vertices
will do!
Denoting the set of all
possible  directed connected chain-like graphs with one input leg and one output leg by $\fG_1^1$,
one can equivalently define an associative algebra structure on a vector space $E$  as
a collection of linear maps, $\{\mu_G: G\langle E\rangle\rar E\}_{G\in \fG_1^1}$,  which satisfy
 the condition,
$$
\mu_G=\mu_{G/H}\circ \mu_H',
$$
 for any subgraph $H\subset G$,  $H\in \fG_1^1$. Here $\mu_H': G\langle E \rangle \rar (G/H)\langle E\rangle$ is the map
which equals $\mu_H$ on the decorated vertices lying in $H$ and which is identity on all other vertices,
while
$\mu_{G/H}: (G/H)\langle E\rangle\rar E$ is the contraction map associated with the graph $G/H$ obtained from $G$ by contracting all vertices lying
in the subgraph $H$ into a single corolla.

\bip

\no{\bf 2.2\, Families of directed labelled  graphs.}
Thus the notion of an associative algebra  can be encoded
into the family of  graphs  $\fG_1^1$ with morphisms of graphs given by contractions along
(admissible) subgraphs. This interpretation
of an associative algebra structure has a strong potential for generalization leading us directly to the notions of
{\em wheeled
props, props, properads, dioperads}\, and {\em operads}\, depending on the way we choose to enlarge the above rather small
and primitive family of graphs $\fG_1^1$. There are several natural enlargements  of $\fG_1^1$:

\no{\sf (i)} $\fG^\circlearrowright$  is, by definition,  the family of {\em arbitrary}\, (not necessarily connected)
 directed graphs built
step-by-step
 from the so called
$(m,n)$-{\em corollas},
\Beq\label{corolla}
\begin{xy}
 <0mm,0mm>*{\bullet};
 <0mm,0mm>*{};<-8mm,3mm>*{}**@{-},
 <0mm,0mm>*{};<-4.5mm,3mm>*{}**@{-},
 <0mm,0mm>*{};<0mm,2.6mm>*{\ldots}**@{},
 <0mm,0mm>*{};<4.5mm,3mm>*{}**@{-},
 <0mm,0mm>*{};<8mm,3mm>*{}**@{-},
<0mm,0mm>*{};<-8mm,-3mm>*{}**@{-},
 <0mm,0mm>*{};<-4.5mm,-3mm>*{}**@{-},
 <0mm,0mm>*{};<0mm,-2.6mm>*{\ldots}**@{},
 <0mm,0mm>*{};<4.5mm,-3mm>*{}**@{-},
 <0mm,0mm>*{};<8mm,-3mm>*{}**@{-};
<0mm,5mm>*{\overbrace{\ \ \ \ \ \ \ \ \ \ \ \ \ \  }};
<0mm,-5mm>*{\underbrace{\ \ \ \ \ \ \ \ \ \ \ \ \ \ }};
<0mm,7mm>*{^{m\ \ output\ legs}};
<0mm,-7mm>*{_{n\ \ input\ legs}};
 \end{xy}, \ \ \ m,n\geq 0,
\Eeq
by taking their disjoint unions and/or gluing some output legs of one corolla with the same number of input legs of
another corolla.
 This is the largest possible
enlargement of $\fG_1^1$ in the class of {\em directed}\, graphs. We have
$\fG^\circlearrowright=\coprod_{m,n\geq 0} \fG^\circlearrowright(m,n)$, where $\fG^\circlearrowright(m,n)
\subset \fG^\circlearrowright$ is the subset of graphs having $m$
output legs and $n$ input legs, e.g.\vspace{-2mm}
\Beq\label{Graph-examples}
\begin{xy}
 <0.4mm,0.0mm>*{};<2.4mm,2.1mm>*{}**@{-},
 <-0.38mm,-0.2mm>*{};<-2.8mm,2.5mm>*{}**@{-},
<0mm,-0.8mm>*{\bullet};
<0mm,-1.0mm>*{};<0mm,-3.6mm>*{}**@{-},
 <2.96mm,2.4mm>*{\bullet};
 <2.4mm,2.8mm>*{};<0mm,5mm>*{}**@{-},
  <3.4mm,3.1mm>*{};<5.1mm,5mm>*{}**@{-},
%
<-2.8mm,2.5mm>*{};<0mm,5mm>*{\bullet}**@{},
<-2.8mm,2.5mm>*{};<0mm,5mm>*{}**@{-},
<0mm,5mm>*{};<0mm,8.6mm>*{}**@{-},
\end{xy}  \in \fG^\circlearrowright(2,1),
\ \ \ \
\begin{xy}
 <0mm,2.47mm>*{};<0mm,-0.5mm>*{}**@{-},
 <0.5mm,3.5mm>*{};<2.2mm,5.2mm>*{}**@{-},
 <-0.48mm,3.48mm>*{};<-2.2mm,5.2mm>*{}**@{-},
 <0mm,3mm>*{\bullet};<0mm,3mm>*{}**@{},
  <0mm,-0.8mm>*{\bullet};<0mm,-0.8mm>*{}**@{},
<0mm,-0.8mm>*{};<-2.2mm,-3.5mm>*{}**@{-},
<0mm,-0.8mm>*{};<2.2mm,-3.5mm>*{}**@{-},
 <-2.5mm,5.7mm>*{\bullet};<0mm,0mm>*{}**@{},
<-2.5mm,5.7mm>*{};<-2.5mm,9.4mm>*{}**@{-},
<-2.5mm,5.7mm>*{};<-5mm,3mm>*{}**@{-},
<-5mm,3mm>*{};<-5mm,-0.8mm>*{}**@{-},
 <-2.5mm,-4.2mm>*{\bullet};<0mm,3mm>*{}**@{},
 <-2.8mm,-3.6mm>*{};<-5mm,-0.8mm>*{}**@{-},
 <-2.5mm,-4.6mm>*{};<-2.5mm,-7.3mm>*{}**@{-},
 (-2.2,9.4)*{}
   \ar@{->}@(ur,dr) (-2.2,-7.4)*{}
\end{xy}
 \in \fG^\circlearrowright(1,1),
\ \ \ \
\begin{xy}
 <0.4mm,0.0mm>*{};<2.4mm,2.1mm>*{}**@{-},
 <-0.38mm,-0.2mm>*{};<-2.8mm,2.5mm>*{}**@{-},
<0mm,-0.8mm>*{\bullet};
<0mm,-1.0mm>*{};<0mm,-3.6mm>*{}**@{-},
 <2.96mm,2.4mm>*{\bullet};
 <2.4mm,2.8mm>*{};<0mm,5mm>*{}**@{-},
  <3.4mm,3.1mm>*{};<5.1mm,5mm>*{}**@{-},
%
<-2.8mm,2.5mm>*{};<0mm,5mm>*{\bullet}**@{},
<-2.8mm,2.5mm>*{};<0mm,5mm>*{}**@{-},
<0mm,5mm>*{};<0mm,8.6mm>*{}**@{-},
\end{xy}
\ \
\begin{xy}
 <0mm,2.47mm>*{};<0mm,-0.5mm>*{}**@{-},
 <0.5mm,3.5mm>*{};<2.2mm,5.2mm>*{}**@{-},
 <-0.48mm,3.48mm>*{};<-2.2mm,5.2mm>*{}**@{-},
 <0mm,3mm>*{\bullet};<0mm,3mm>*{}**@{},
  <0mm,-0.8mm>*{\bullet};<0mm,-0.8mm>*{}**@{},
<0mm,-0.8mm>*{};<-2.2mm,-3.5mm>*{}**@{-},
<0mm,-0.8mm>*{};<2.2mm,-3.5mm>*{}**@{-},
 <-2.5mm,5.7mm>*{\bullet};<0mm,0mm>*{}**@{},
<-2.5mm,5.7mm>*{};<-2.5mm,9.4mm>*{}**@{-},
<-2.5mm,5.7mm>*{};<-5mm,3mm>*{}**@{-},
<-5mm,3mm>*{};<-5mm,-0.8mm>*{}**@{-},
 <-2.5mm,-4.2mm>*{\bullet};<0mm,3mm>*{}**@{},
 <-2.8mm,-3.6mm>*{};<-5mm,-0.8mm>*{}**@{-},
 <-2.5mm,-4.6mm>*{};<-2.5mm,-7.3mm>*{}**@{-},
 (-2.2,9.4)*{}
   \ar@{->}@(ur,dr) (-2.2,-7.4)*{}
\end{xy}
\in \fG^\circlearrowright(3,2).
\vspace{-2mm}
\Eeq


\no{\sf (ii)}
$\fG_c^\circlearrowright=\coprod_{m,n\geq 0} \fG_c^\circlearrowright(m,n)$ is a subset of $\fG^\circlearrowright$
consisting of {\em connected}\, graphs. For example, the first two graphs in (\ref{Graph-examples}) belong to $\fG_c^\circlearrowright$
while the third one (which is the disjoint union of the first two graphs) does not.


\no{\sf (iii)}
$\fG^\uparrow=\coprod_{m,n\geq 0} \fG^\uparrow(m,n)$ is a subset of $\fG^\circlearrowright$ consisting of directed graphs
with no {\em closed}\, directed paths of internal edges which begin and end at the same vertex,
 e.g.\ the first graph in (\ref{Graph-examples}) belongs to
$\fG^\uparrow$, while the other two do not.


\no{\sf (iv)}
$\fG_c^\uparrow:= \fG^\uparrow \cap \fG_c^\circlearrowright$.

\no{\sf (v)}
$\fG^\upharpoonleft_{c,0}$ is a subset of $\fG_c^\uparrow$ consisting of graphs of {\em genus zero}\, (as $1$-dimensional
$CW$ complexes).

\no{\sf (vi)}
$\fG^1$ is a subset of $\fG^\uparrow_{c,0}$ built from corollas (\ref{corolla}) of type $(1,n)$ only, $n\geq 1$.
We have $\fG^1=\coprod_{n\geq 1} \fG^1(1,n)$ and we further abbreviate $\fG^1(n):=\fG^1(1,n)$; thus $\fG^1(n)$  is the subset of
$\fG^1$ consisting of graphs with precisely $n$ input legs. All graphs in $\fG^1$ have precisely one output leg.


\no
Let $\fG^\checkmark$ be any of the above mentioned families of graphs. {\em
We assume from now on that input and out output legs (if any) of graphs from $\fG^\checkmark(m,n)\subset \fG^\checkmark$
are bijectively {\em labelled}\, by elements of the sets $[n]$ and $[m]$ respectively}. Hence the group $\bS_m\times \bS_n$
naturally acts on the set $\fG^\checkmark(m,n)$ by permuting the labels\footnote{in the case of $\fG^1$ the action of
the factor $\bS_m$ is, of course, trivial.}.

\sip

\no{\bf 2.3\, Decorations of directed labelled  graphs.}
Next we have to think on what to use for {\em decorations}\, of
the vertices of a graph $G\in \fG^\checkmark(m,n)$. The presence of the family of the permutation groups
$\{\bS_m\times \bS_n\}_{m,n\geq 0}$ suggests the following notion:
 an $\bS$-{\em bimodule}, $E$, is, by definition, a collection
of graded vector spaces, $\{E(m,n)\}_{m,n\geq 0}$, equipped with a left
action of the group $\bS_m$  and with a right action of $\bS_n$ which commute
with each other. For example, for any graded vector space $V$ the collection,
${\End}\langle V\rangle
=\{{\End}\langle V\rangle(m,n):= \Hom(V^{\ot n}, V^{\ot m})\}_{m,n\geq 0}$,
is naturally an $\bS$-bimodule.

Let $E$ be an $\bS$-bimodule and  $G\in \fG^\circlearrowright(m,n)$ an arbitrary graph. The graph $G$ is built by definition
from a number of various $(p,q)$-corollas
constituting a set which we denote by ${\sV}(G)$ and call
 the set of {\em vertices}\, of
$G$; the set of output (resp.\ input) legs of a vertex $v\in \sV(G)$ is denoted by $\out_v$ (resp.\ by $\inn_v$).
Let $\langle [p]\rar \out_v\rangle$ be the $p!$-dimensional vector space generated over $\K$ by the set
of all bijections from  $[p]$ to  $\out_v$, i.e.\ by the set of all possible labeling of $\out_v$ by integers;
it is naturally a right $\bS_p$-module; we define analogously a left $\bS_q$-module $\langle \inn_v\rar [n]\rangle$ and
then define a vector space,
$$
E(\out_v, \inn_v):= \langle [p]\rar \out_v\rangle \otimes_{\bS_p} E(p,q)
 \otimes_{\bS_q} \langle \inn_v\rar [q]\rangle.
$$
An element of  $E(\out_v, \inn_v)$ is  called a {\em decoration of the vertex}\, $v\in \sV(G)$.
To define next a space of {\em decorations of a graph}\, $G$ we should think of
 taking  the tensor product of the constructed vector spaces $E(\out_v, \inn_v)$
over all vertices $v\in \sV(G)$ but face a problem
that the set $\sV(G)$ is unordered so that the ordinary definition of the tensor product of vector spaces does not immediately
apply.
The solution is to consider first {\em all}\, possible linear orderings, $\ga: [k]\rar \sV(G)$,   $k:=|V(G)|$,
 of the set $\sV(G)$ and then take coinvariants,
$$
\otimes_{v\in \sV(G)} E(\out_v, \inn_v):= (\oplus_{\ga}
 E\left(\out_{\ga(1)}, \inn_{\ga(1)})
\ot\ldots \ot
 E(\out_{\ga(k)}, \inn_{\ga(k)})\right)_{\bS_k},
$$
with respect to the natural action of the group $\bS_{k}$ permuting the orderings. Now we are ready to define the
vector space of {\em decorations}\, of the graph $G$ as a quotient of the unordered tensor product,
$$
G\langle E\rangle:= (
\otimes_{v\in \sV(G)} E(\out_v, \inn_v))_{Aut G}
$$
with respect to the automorphism group of the graph $G$ which is, by definition, the subgroup of the
 symmetry group of the 1-dimensional
$CW$-complex underlying the graph $G$ which fixes the legs.
An element of $G\langle E\rangle$ is   called
a {\em graph $G$ with internal vertices decorated by elements of $E$}, Thus a decorated graph is essentially a pair,
$(G, [a_1\ot \ldots \ot a_{k}])$, consisting of a graph $G$ with $k=|\sV(G)|$ and an equivalence class
of tensor products of elements $a_\bullet \in E$.
Note that if $E=\{E(m,n)\}$ is a {\em dg}\, $S$-bimodule (i.e.\ each $E(m,n)$ is a complex
equipped  with an $\bS_m\times \bS_n$-equivariant differential $\delta$) then $G\langle E\rangle$ is naturally
a dg vector space with the differential
$$
\delta_G \left(G, \left[a_1\ot \ldots \ot a_{k|}\right]\right):=
(G, [\sum_{i=1}^{k=|\sV(G)|}(-1)^{a_1+\ldots + a_{i-1}}a_1\ot \ldots
\ot \delta a_i\ot \ldots \ot a_{k}]).
$$
Note also that if $G$ is an $(m,n)$-graph with one internal vertex (i.e.\ an $(m,n)$-corolla), then
$G\langle E\rangle$ is canonically isomorphic to $E(m,n)$.

\sip

\no{\bf 2.4 Props, properads, dioperads and  operads} Let $\fG^\checkmark$ be one of the families of graphs introduced in
\S~{2.2}. A subgraph $H\subset G$ of a graph $G\in \fG^\checkmark$ is called {\em admissible}\,
if both $H$ and $G/H$ also belong to $\fG^\checkmark$, where $G/H$ is the graph obtained from $G$ by shrinking all vertices and
all internal edges of $H$ into a new single vertex.

\sip

\no{\bf 2.4.1\, Definition.}
A $\fG^\checkmark$-{\em algebra}\, is an $\bS$-bimodule $E=\{E(m,n)\}$ together with
a collection of linear maps,
$
\left\{\mu_G: G\langle E\rangle\rar E\right\}_{G\in \fG^\checkmark}$,
satisfying the ``associativity" condition,
\Beq\label{graph-associativity}
\mu_G=\mu_{G/H}\circ \mu_H',
\Eeq
 for any admissible subgraph $H\subset G$, where $\mu_H': G\langle E \rangle \rar (G/H)\langle E\rangle$ is the map
which equals $\mu_H$ on the decorated vertices lying in $H$ and which is identity on all other vertices of $G$.
If the $\bS$-bimodule $E$ underlying a $\fG^\checkmark$-algebra has a differential $\delta$ satisfying,
for any $G\in \fG^\checkmark$, the condition $\delta\circ \mu_G=\mu_G\circ \delta_G$, then the $\fG^\checkmark$-algebra
is called {\em differential}.

\sip

\no{\bf 2.4.2\, Remarks.}
(a)
For the family of graphs $\fG^1$ the condition (\ref{graph-associativity}) is void for elements in $E(m,n)$ with $m\neq 1$.
Thus we may assume without loss of generality that a $\fG^1$-algebra $E$ satisfies an extra condition
that $E(m,n)= 0$ unless $m=1$. For the same reason we may assume that a $\fG_1^1$-algebra $E$ satisfies
$E(m,n)=0$ unless $m=n=1$.

\sip

\no (b) As we have an obvious identity $\mu_G= \mu_{G/G}\circ \mu_G'$, the ``associativity" condition
(\ref{graph-associativity}) can  be equivalently reformulated as follows: for any two admissible
subgraphs $H_1,H_2\subset G$ one has
\Beq\label{graph-associativity-2}
\mu_{G/H_1}\circ \mu_{H_1}'=\mu_{G/H_2}\circ \mu_{H_2}',
\Eeq
i.e.\ the contraction of a decorated graph $G$ into a decorated corolla along a family admissible subgraphs
does not depend on particular
choices of these subgraphs (if there are any). This is indeed a natural extension of the notion of
 associativity from 1 dimension to 3 dimensions,
and hence we can omit double commas in the term.

\sip

\no{\bf 2.4.3\, Definitions} (see, e.g.,  \cite{MSS, Va, BM,  Me1} and references cited there)

(i) An $\fG^\circlearrowright$-algebra $E$ is called a {\em wheeled prop}.

 (ii) An $\fG^\circlearrowright_c$-algebra is called a {\em wheeled properad}.

(ii) An $\fG^\uparrow$-algebra  is called a {\em prop}.

(iv) An $\fG^\uparrow_c$-algebra  is called a {\em properad}.

(v) An $\fG^\uparrow_{0,c}$-algebra  is called a {\em dioperad}.

(vi) An $\fG^1$-algebra  is called an {\em operad}.

(vii) An $\fG^1_1$-algebra  is called an {\em associative algebra}.

\sip

\no{\bf 2.4.4\, Remarks.}
(a) We have an obvious chain of inclusions of the  categories of $\fG^\checkmark$-algebras,
$$
\mbox{
(vii)$\subset$(vi)$\subset$(v)$\subset$(iv)$\subset$(iii)$\subset$(ii)$\subset$(i).
}
$$


\no (b)
Note that {\em every}\, subgraph of a graph in $\fG^\circlearrowright$ is admissible. In this sense wheeled props
are the most general and natural algebraic structures associated with the class of {\em directed}\, graphs.
The set of independent operations in a $\fG^\circlearrowright$-algebra is generated by one vertex graphs
with at least one {\em loop} (that is,
an internal edge beginning and ending at the vertex) and by two vertex graphs without closed directed paths
(i.e.\ the ones belonging to $\fG^\uparrow$).
\sip

\no
(c) By contrast to wheeled props, the set of operations in an ordinary prop, i.e.\ in a $\fG^\uparrow$-algebra,
is generated by the set of two vertex graphs only, and, as it is not hard to check, if the associativity condition hold for three vertex
graphs, then it holds for arbitrary graphs in $\fG^\uparrow$. 
This is {\em not}\, true for $\fG^\circlearrowright$-algebras
which is a first indication that the homotopy theory for wheeled props should be substantially different from the one
 for ordinary props.

\sip

\no (d)
If we forget orientations (i.e.\ the flow) on edges and  work instead with a family of undirected graphs, $\fG$,
built, by definition, from corollas with $m\geq 1$ undirected legs via their gluings, then we get
 a notion of  $\fG$-{\em algebra}\,
which is closely related to the notion of {\em modular operad} \cite{GeKa}.

\sip

\no{\bf 2.5.\, First basic example: endomorphism $\fG^\checkmark$-algebras.}
 For any finite-dimensional vector space $V$ the $\bS$-bimodule
$\End_V=\{ \Hom(V^{\ot n}, V^{\ot m})\}$ is naturally a
 $\fG^\checkmark$-algebra called the {\em endomorphism
$\fG^\checkmark$-algebra
of $V$}.\footnote{For $\fG^1$-algebras, that is for operads, it is
enough to restrict oneself to the case $m=1$ only,
 i.e.\  set, by default, $\End_V:=\{ \Hom(V^{\ot n}, V\}_{n\geq 1}$ (cf.\ \S 2.4.2(a)).}
For any two vertex graph $G\in \fG^\uparrow_c$ the associated composition $\mu_G: G\langle \End_V \rangle \rar \End_V$
is the ordinary composition of two linear maps;
for a one vertex graph $G\in \fG^\circlearrowright$ with say $k$ loops the associated map $\mu_G$ is the ordinary $k$-fold trace
of a linear map; for a two vertex disconnected graph $G\in \fG^\uparrow$ the associated map $\mu_G$ is the ordinary
tensor product of linear maps. It is easy to see that all the axioms are satisfied.

\sip

Note that for all $\fG^\checkmark$-algebras except $\fG^\circlearrowright$ and $\fG^\circlearrowright_c$
the basic algebraic operations $\mu_G$ do {\em not}\, involve traces so that the above assumption on
 {\em finite-dimensionality}
\,  of $V$ can be dropped for endomorphism props, properads, dioperads and operads.
If $V$ is a (finite-dimensional) {\em dg}\, vector space, then $\End_V$ is naturally a {\em dg}\,
 $\fG^\checkmark$-algebra.

\sip

\no{\bf 2.6\, Second basic example: a free $\fG^\checkmark$-algebra.}
For an  $\bS$-bimodule,
 $E=\{E(m,n)\}$, one can construct another $\bS$-bimodule,
$\cF^\checkmark\hspace{-0.6mm} \langle E\rangle=\{\cF^\checkmark \hspace{-0.6mm}\langle E\rangle (m,n)\}$ with
$$
\cF^\checkmark\hspace{-0.6mm} \langle E\rangle (m,n):= \bigoplus_{G\in \fG^\checkmark(m,n)} G\langle E\rangle.
$$
This $\bS$-bimodule $\cF^\checkmark\hspace{-0.6mm} \langle E\rangle$ has a natural $\fG^\checkmark$-algebra
structure  with the contraction maps
$\mu_G$ being tautological. The $\fG^\checkmark$-algebra  $\cF^\checkmark\hspace{-0.6mm} \langle E\rangle$ is called the
{\em free $\fG^\checkmark$-algebra} (i.e., respectively, {\em the free wheeled prop, the free prop, the free dioperad} etc)
{\em generated
by the $\bS$-bimodule $E$}.

\sip

\no{\bf 2.7\, Morphisms of $\fG^\checkmark$-algebras.}  A morphisms of $\fG^\checkmark$-algebras, $\rho: \cP_1\rar \cP_2$,
is a morphism of the underlying $\bS$-bimodules such that, for any graph $G\in \fG^\checkmark$,
 one has $\rho\circ \mu_G= \mu_G\circ (\rho^{\ot G})$,
where $\rho^{\ot G}$ means a map, $G\langle \cP_1\rangle \rar G\langle \cP_2\rangle$, which changes decorations of each vertex
 in $G$ in accordance with $\rho$. It is often assumed by default that a morphism $\rho$ is homogeneous which 
(almost always)
 implies that $\rho$ has degree $0$. Unless otherwise is explicitly stated we do {\em not}\,  assume
in this paper that morphisms of $\fG^\checkmark$-algebras are homogeneous so that they can have nontrivial parts
in degrees other than zero.
A morphism of $\fG^\checkmark$-algebras, $\cP\rar \End\langle V\rangle$, is called a {\em representation}\, of the
$\fG^\checkmark$-algebra $\cP$ in a graded
vector space $V$.
If $\cP_1$ is a free $\fG^\checkmark$-algebra, $\cF^\checkmark\hspace{-0.6mm}
 \langle E\rangle$, generated by some $\bS$-bimodule $E$,
then the set of morphisms $\fG^\checkmark$-algebras, $\{\rho: \cP_1\rar \cP_2\}$, is in one-to-one correspondence with the
set of morphisms of $\bS$-bimodules, $\{\rho|_E: E\rar \cP_2\}$, i.e.\ a $\fG^\checkmark$-morphism is uniquely
determined by its values on the generators. In particular, the set of morphisms,  $\cF^\checkmark\hspace{-0.6mm}
 \langle E\rangle\rar \cP_2$, has a graded vector space structure for any $\cP_2$.

\sip

 A {\em free resolution}\, of a dg $\fG^\checkmark$-algebra
$\cP$ is, by definition, a dg free $\fG^\checkmark$-algebra, $(\cF^\checkmark \hspace{-0.6mm}\langle E \rangle, \delta)$,
generated by some $\bS$-bimodule $E$ together with a degree zero morphism of dg $\fG^\checkmark$-algebras,
$\pi: (\cF^\checkmark \hspace{-0.6mm}\langle E \rangle, \delta) \rar \cP$, which induces a cohomology isomorphism.
If the differential $\delta$ in $\cF^\checkmark \hspace{-0.6mm}\langle \cE \rangle$ is
decomposable with respect to compositions $\mu_G$, then $\pi: (\cF^\checkmark \hspace{-0.6mm}\langle E \rangle, \delta) \rar
\cP$ is called a {\em minimal model}\, of $\cP$. In this case the
free algebra $\cF^\checkmark \hspace{-0.6mm}\langle E \rangle$ is  often denoted by
$\cP_\infty$.

\sip

\no{\bf 2.8 Props and properads.}
We shall work in this paper only with $\fG^\uparrow$- and $\fG_c^\uparrow$-algebras,
i.e.\ with props and properads. For later use we mention several useful constructions with these graph-algebras.

\sip

(i) There is a functor, $\Psi$, which associates canonically to an arbitrary dg properad, $\cP$, an associated dg prop
$\Psi(\cP)$ \cite{Va}. As we are working over a field of characteristic $0$, this functor is, by K\"unneth theorem,
exact, i.e.\ $\Psi(H(\cP))= H(\Psi(\cP))$. For example, if $\cP$ is a dg free properad 
$(\cF^\uparrow_c\langle E\rangle, \delta)$,
then $\Psi(\cP)$ is precisely $\cF^\uparrow\langle E\rangle$ with the same differential (as given on the
generators).

(ii) The above mentioned functor, $\cF^\uparrow:(E, \delta)\rar
(\cF^\uparrow\langle E\rangle, \delta)$,
in the category of dg $\bS$-bimodules is also exact,
$H(\cF^\uparrow\langle E\rangle)= \cF^\uparrow\langle H(E)\rangle$.  Moreover,
if we set in this situation $\mbox{Rep}_V(\cF^\uparrow\langle E\rangle)$
  for the vector
space of all possible representations, $\{\rho: \cF^\uparrow\langle E\rangle\rar \End_V\}\simeq \Hom(E, \End_V)$,
and define a differential $\delta$ in
$\mbox{Rep}_V(\cF^\uparrow\langle E\rangle)$ by the formula $\delta\rho:= \rho\circ \delta$, then the resulting functor,
$(E,\delta)\rar (\mbox{Rep}_V(\cF^\uparrow\langle E\rangle), \delta)$, in the category of complexes is exact, i.e.\
\Beq\label{main}
H(\mbox{Rep}_V(\cF^\uparrow\langle E\rangle))= \mbox{Rep}_V(H(\cF^\uparrow\langle E\rangle))=
\mbox{ Rep}_V(\cF^\uparrow\langle H(E)\rangle).\hspace{-1mm}
\Eeq
Indeed, as we are working over a field of characteristic zero, we can always choose an equivariant chain homotopy
between complexes $(E,\delta)$ and $(H(E), 0)$. This chain homotopy induces a chain homotopy between complexes
$(\mbox{Rep}_V(\cF^\uparrow\langle E\rangle), \delta)$ and $(\mbox{Rep}_V(\cF^\uparrow\langle H(E)\rangle), 0)$
proving formula (\ref{main}).

\sip
(iii)
There is also a natural {\em parity change}\, functor, $\Pi$, which associates with a dg prop(erad) $\cP$
a dg prop(erad) $\Pi\cP$ with the following property: every representation of $\Pi\cP$ in a graded vector
space $V$ is equivalent to
a representation of $\cP$ in $V[1]$. This functor is also exact. If, for example, $\cP$ is a dg free prop
$\cF^\uparrow\langle E\rangle$ generated by an $\bS$-bimodule $E=\{E(m,n)\}$, then, as it is not hard to check,
$\Pi\cP= \cF^\uparrow\langle \check{E}\rangle$, where $\check{E}:=\{\sgn_m\ot E(m,n)\ot \sgn_n[m-n]\}$ and $\sgn_m$
stands for the one-dimensional sign representation of $\bS_m$.


\section{Simplicial and permutahedra cell complexes as dg properads}
\label{ch3}

\no{\bf 3.1 Simplices as a
 dg properad.} A geometric $(n\hspace{-0.5mm}-\hspace{-0.5mm} 1)$-{\em simplex}, $\Delta_{n-1}$,
 is, by definition, a  subset in $\R^{n}=\{x^1, \ldots, x^{n}\}$, $n\hspace{-0.5mm}\geq\hspace{-0.5mm} 1$,
  satisfying the equation $\sum_{i=1}^nx^i =1$, $x^i\geq 0$ for all $i$.
 To define its cell complex one has to choose an orientation on $\Delta_{n-1}$ which is the same as to choose an orientation on the hyperplane
$\sum_{i=1}^n x^i=1$. We induce it from the standard orientation on $\R^{n+1}$ by requiring that the
manifold with boundary defined by the equation $\sum_{i=1}x^i\leq 1$ is naturally oriented.
Let $(C_\bullet(\Delta_{n-1})=\oplus_{k=0}^{n-1}C_{-k}(\Delta_{n-1}), \delta)$ be the standard (non-positively graded) cell complex of $\Delta_{n-1}$.
By definition, $C_{-k}(\Delta_{n-1})$ is a $\binom nk$-dimensional vector space spanned by $k$-dimensional cells,
$$
\triangle^I_{n-1}:=\left\{(x^1, \ldots, x^n)\in \Delta_{n-1}|x^i=0, i\in I\right\},
$$
parameterized by all possible subsets $I$ of $[n]$ of cardinality $n-k-1$ and equipped with the natural orientations (which we describe explicitly below).

Note that the action , $(x^1,\ldots,x^n)\rar (x^{\sigma(1)}, \ldots,
x^{\sigma(n)})$, of the permutation group $\bS_n$ on $\R^n$ leaves $\Delta_{n-1}$ invariant as a subset but {\em not}\,
as an oriented manifold with boundary. As an $\bS_n$-module, one can obviously identify $C_{1-n}(\Delta_{n-1})$ with $\sgn_n[n-1]$,
and hence one can represent pictorially the oriented cell $\triangle^\emptyset_{n-1}$ as a labelled $(0,n)$-corolla,
$\xy
 <0mm,0mm>*{\blacktriangle};
<0.3mm,-0.3mm>*{};<1.5mm,4mm>*{}**@{-},
<0.3mm,-0.3mm>*{};<3.5mm,4mm>*{}**@{-},
<5.5mm,3.7mm>*{...};
<0.3mm,-0.3mm>*{};<9mm,4mm>*{}**@{-},
 <1.5mm,5.5mm>*{{_1}};
 <4mm,5.5mm>*{{_2}};
 <9.5mm,5.5mm>*{{_n}};
\endxy
$, with the symmetry condition
\Beq\label{simplex-top-cell}
\xy
 <0mm,0mm>*{\blacktriangle};
<0.3mm,-0.3mm>*{};<1.5mm,4mm>*{}**@{-},
<0.3mm,-0.3mm>*{};<3.5mm,4mm>*{}**@{-},
<5.5mm,3.7mm>*{...};
<0.3mm,-0.3mm>*{};<9mm,4mm>*{}**@{-},
 <1.5mm,5.5mm>*{{_1}};
 <4mm,5.5mm>*{{_2}};
 <9.5mm,5.5mm>*{{_n}};
\endxy
= (-1)^\sigma \ \
\xy
 <0mm,0mm>*{\blacktriangle};
<0.3mm,-0.3mm>*{};<1.5mm,4mm>*{}**@{-},
<0.3mm,-0.3mm>*{};<3.5mm,4mm>*{}**@{-},
<5.7mm,3.7mm>*{...};
<0.3mm,-0.3mm>*{};<9mm,4mm>*{}**@{-},
 <0mm,5.9mm>*{{_{\sigma(1)}}};
 <5.5mm,5.9mm>*{{_{\sigma(2)}}};
 <13mm,5.9mm>*{{_{\sigma(n)}}};
\endxy
, \ \ \
\forall \sigma\in \bS_n.
\Eeq
The boundary of $\triangle^\emptyset_{n-1}$ is a union of $n$ cells, $\triangle_{n-1}^i$, $i=1,\ldots, n$, of dimension $n-2$.
The permutation group $\bS_n$ permutes, in general, these cells and changes their natural orientations while keeping
their linear span $C_{2-n}(\Delta_{n-1})$ invariant. It is obvious that the subgroup
$G_i:=\{\sigma \in \bS_n\mid \sigma(i)=i\}\simeq \bS_{n-1}$  of $\bS_n$ is a symmetry group of the cell $\triangle_{n-1}^i$ as an
{\em unoriented}\, manifold with boundary. If we take the orientation into account, then the vector subspace of
$C_{2-n}(\Delta_{n-1})$  spanned by $\triangle_{n-1}^i$ can be identified as an $\bS_{n-1}$-module with $\sgn_{n-1}$,
and hence the  $n$-dimensional space $C_{2-n}(\Delta_{n-1})$ itself can be identified as an $\bS_n$-module with
$\K[\bS_n]\ot_{\bS_{n-1}} \sgn_{n-1}[n-2]$. Its basis elements, $\triangle_{n-1}^i$, can be pictorially represented as
 $(0,n)$-corollas
$\xy
 <0mm,0mm>*{\blacktriangle};
<-0.3mm,0.3mm>*{};<-2.5mm,4mm>*{}**@{-},
<0.3mm,-0.3mm>*{};<1.5mm,4mm>*{}**@{-},
<0.3mm,-0.3mm>*{};<3.mm,4mm>*{}**@{-},
<3.4mm,5.5mm>*{...};
<8mm,5.5mm>*{...};
<0.3mm,-0.3mm>*{};<8mm,4mm>*{}**@{-},
 <1.2mm,5.5mm>*{{_1}};
 <-3.3mm,5.5mm>*{{_i}};
<5.8mm,5.8mm>*{{_{\hat{i}}}};
 <10.3mm,5.5mm>*{{_n}};
\endxy
$
with the legs in the right bunch being skew symmetric with respect to the change of labellings by an element
$\sigma \in G_i$ (cf.\ (\ref{simplex-top-cell})). The boundary operator $\delta: C_{1-n}(\Delta_{n-1}) \rar
C_{2-n}(\Delta_{n-1})$ is equivariant with respect to the $\bS_n$-action and is given on the generators by the formula,
$$
\delta\
\xy
 <0mm,0mm>*{\blacktriangle};
<0.3mm,-0.3mm>*{};<1.5mm,4mm>*{}**@{-},
<0.3mm,-0.3mm>*{};<3.5mm,4mm>*{}**@{-},
<4.7mm,3.3mm>*{...};
<0.3mm,-0.3mm>*{};<8mm,4mm>*{}**@{-},
 <1.5mm,5.5mm>*{{_1}};
 <4mm,5.5mm>*{{_2}};
 <8.5mm,5.5mm>*{{_n}};
\endxy
=\sum_{i=1}^n (-1)^{i+1}
\xy
 <0mm,0mm>*{\blacktriangle};
<-0.3mm,-0.3mm>*{};<-2.5mm,4mm>*{}**@{-},
<0.3mm,-0.3mm>*{};<1.5mm,4mm>*{}**@{-},
<0.3mm,-0.3mm>*{};<3.mm,4mm>*{}**@{-},
<3.4mm,5.5mm>*{...};
<8mm,5.5mm>*{...};
<0.3mm,-0.3mm>*{};<8mm,4mm>*{}**@{-},
 <1.2mm,5.5mm>*{{_1}};
 <-3.3mm,5.5mm>*{{_i}};
<5.8mm,5.5mm>*{{_{\hat{i}}}};
 <10.3mm,5.5mm>*{{_n}};
\endxy, \ \ \ \hat{i}\ \mbox{omitted},
$$
More generally, the symmetry group of, say, a cell $\triangle_{1-n}^I\in
C_{1-k-n}(\Delta_{n-1})$ with $I=\{i_1<i_2<\ldots < i_{n-k}\}$  is $\bS_{k}\times \bS_{n-k}\subset \bS_n$
with $\bS_{k}\times \Id$ leaving the orientation of $\triangle_{n-1}^I$ invariant and $\Id\times \bS_{n-k}$
changing the orientation  via the sign
representation. Thus we can identify the oriented cell $\triangle_{n-1}^I$, as an
element of the $\bS_n$-module $C_{1-k-n}(\Delta_{n-1})$, with a $(0,n)$-corolla,
$
\xy
 <0mm,0mm>*{\blacktriangle};
<-0.3mm,-0.3mm>*{};<-1.5mm,4mm>*{}**@{-},
<-0.3mm,-0.3mm>*{};<-5mm,4mm>*{}**@{-},
<-0.3mm,-0.3mm>*{};<-8mm,4mm>*{}**@{-},
<-6.4mm,5.3mm>*{...};
<0.3mm,-0.3mm>*{};<1.5mm,4mm>*{}**@{-},
<0.3mm,-0.3mm>*{};<3.4mm,4mm>*{}**@{-},
<4.6mm,5.3mm>*{...};
<0.3mm,-0.3mm>*{};<8mm,4mm>*{}**@{-},
 <2mm,5.5mm>*{{_{j_1}}};
 <9.5mm,5.5mm>*{{_{j_{n-k}}}};
<-9mm,5.5mm>*{{_{i_1}}};
<-2mm,5.5mm>*{{_{i_k}}};
\endxy$, which has `symmetric' output legs in the left bunch and `skewsymmetric' output legs in the right one. Here
$\{j_1< \ldots <j_{n-k}\}:=[n]\setminus I$. The $\bS_n$-module, $C_{1-k-n}(\Delta_{n-1})$ is then canonically isomorphic to $
E_{k}(n):=\K[\bS_n]_{\bS_k\times \bS_{n-k}} (\id_k\ot \sgn_{n-k})$,
where $\id_k$ stands for the trivial 1-dimensional representation of $\bS_k$.
The boundary operator $\delta: C_{1-k-n}(\Delta_{n-1}) \rar
C_{2-k-n}(\Delta_{n-1})$ is equivariant with respect to the $\bS_n$-action and is given on the
generators by the formula,
\Beq\label{simplex-differential}
\delta\
\xy
 <0mm,0mm>*{\blacktriangle};
<-0.3mm,-0.3mm>*{};<-2.5mm,4mm>*{}**@{-},
<-0.3mm,-0.3mm>*{};<-5mm,4mm>*{}**@{-},
<-0.3mm,-0.3mm>*{};<-8mm,4mm>*{}**@{-},
<-6.4mm,5.3mm>*{...};
<0.3mm,-0.3mm>*{};<1.5mm,4mm>*{}**@{-},
<0.3mm,-0.3mm>*{};<5mm,4mm>*{}**@{-},
<6.2mm,5.3mm>*{...};
<0.3mm,-0.3mm>*{};<8mm,4mm>*{}**@{-},
 <2mm,5.5mm>*{{_{k+1}}};
 <9.5mm,5.5mm>*{{_{{n}}}};
<-9mm,5.5mm>*{{_{1}}};
<-3mm,5.5mm>*{{_{k}}};
\endxy
=
\sum_{i=k+1}^n (-1)^{i+1}
\xy
 <0mm,0mm>*{\blacktriangle};
<-0.3mm,-0.3mm>*{};<-3mm,4mm>*{}**@{-},
<-0.3mm,-0.3mm>*{};<-6.5mm,4mm>*{}**@{-},
<-0.3mm,-0.3mm>*{};<-9mm,4mm>*{}**@{-},
<-0.3mm,-0.3mm>*{};<-12mm,4mm>*{}**@{-},
<-10.4mm,5.3mm>*{...};
<0.3mm,-0.3mm>*{};<3mm,4mm>*{}**@{-},
<0.3mm,-0.3mm>*{};<9mm,4mm>*{}**@{-},
<6.2mm,5.3mm>*{...};
<11.2mm,5.3mm>*{...};
<0.3mm,-0.3mm>*{};<12mm,4mm>*{}**@{-},
 <2mm,5.5mm>*{{_{k+1}}};
 <8.7mm,5.5mm>*{{_{{\hat{i}}}}};
<14mm,5.5mm>*{{_{{n}}}};
<-13mm,5.5mm>*{{_{1}}};
<-7mm,5.5mm>*{{_{k}}};
<-4mm,5.5mm>*{{_{i}}};
\endxy
\Eeq
Thus we proved the following

\sip

\no{\bf 3.1.1 Proposition.} {\em  (i) The standard simplicial cell complex is canonically isomorphic to a dg free properad,
$\cS:= \cF_c^\uparrow\langle E\rangle$, generated by an $\bS$-bimodule, $E=\{E(m,n)\}$,
$$
E(m,n)=\left\{
\Ba{lr}
\bigoplus_{k=0}^{m-1}E_k(m)[m-k-1] = \mbox{span}
\langle\hspace{-2mm}
\xy
 <0mm,0mm>*{\blacktriangle};
<-0.3mm,-0.3mm>*{};<-1.5mm,4mm>*{}**@{-},
<-0.3mm,-0.3mm>*{};<-5mm,4mm>*{}**@{-},
<-0.3mm,-0.3mm>*{};<-8mm,4mm>*{}**@{-},
<-6.4mm,5.3mm>*{...};
<0.3mm,-0.3mm>*{};<1.5mm,4mm>*{}**@{-},
<0.3mm,-0.3mm>*{};<3.4mm,4mm>*{}**@{-},
<4.6mm,5.3mm>*{...};
<0.3mm,-0.3mm>*{};<8mm,4mm>*{}**@{-},
 <2mm,5.5mm>*{{_{j_1}}};
 <9.5mm,5.5mm>*{{_{j_{m-k}}}};
<-9mm,5.5mm>*{{_{i_1}}};
<-2mm,5.5mm>*{{_{i_k}}};
\endxy
\hspace{-4mm}
\rangle_{0\leq k\leq m-1}
 & \mbox{for}\ n=0,  \vspace{3mm}\\
0 & \mbox{for}\  n\geq 1,
\Ea
\right.
$$
and equipped with the differential given on the generators by (\ref{simplex-differential}).

(ii) The cohomology of $(\cS, \delta)$ is concentrated in degree zero and equals the free properad generated by the
following  degree
zero graphs with `symmetric' legs,
\Beq\label{simplex-cohomology}
\xy
 <0mm,0mm>*{\bullet};
<0mm,0mm>*{};<-4mm,4mm>*{}**@{-},
<0mm,0mm>*{};<-2mm,4mm>*{}**@{-},
<0mm,0mm>*{};<4mm,4mm>*{}**@{-},
<1mm,4mm>*{...};
 <-4.5mm,5.6mm>*{_1};
<-2mm,5.6mm>*{_2};
<4.5mm,5.6mm>*{_m};
\endxy
:=
\sum_{\sigma\in \bS_m}
\xy
 <0mm,0mm>*{\blacktriangle};
<-0.3mm,-0.3mm>*{};<-1.5mm,4mm>*{}**@{-},
<-0.3mm,-0.3mm>*{};<-5mm,4mm>*{}**@{-},
<-0.3mm,-0.3mm>*{};<-8mm,4mm>*{}**@{-},
<-7mm,5.3mm>*{...};
<0.3mm,-0.3mm>*{};<8mm,4mm>*{}**@{-},
 <12mm,5.5mm>*{{_{\sigma(m)}}};
<-12mm,5.5mm>*{{_{\sigma(1)}}};
<0mm,6mm>*{{_{\sigma(m-1})}};
\endxy,
\ \ m\geq 1.
\Eeq

}
\sip

Claim~3.1.1(ii) follows from the contractibilty  of simplices, and, for each $m$,  graph (\ref{simplex-cohomology})
represents the sum
of all vertices of $\Delta_{m-1}$.
\sip

\no{\bf 3.1.2 From simplicia to Koszul complex.} The  vector space, $\mbox{Rep}_V(\cS)$, of representations,
$\rho: \cS\rar \End_V$, of the simplicial properad in a vector space $V$ can be identified with
$\sum_{k=0}^{m-1} \odot^k V\ot \wedge^{m-k} V[m-k-1]$. We can naturally make the latter into complex by setting
$d\rho:= \rho\circ \p$ (cf.\ \S 2.8ii). It is easy to
see that we get in this way, for each $m\geq 1$,
$$
\wedge^m V\stackrel{d}{\lon} V\ot \wedge^{m-1}V \stackrel{d}{\lon}  \odot^2 V\ot \wedge^{m-2}V
 \stackrel{d}{\lon}\ldots  \stackrel{d}{\lon}  \odot^{m-1} V\ot V,
 $$
 the classical Koszul complex. Hence Proposition~3.1.1(ii) and isomorphism (\ref{main})
 imply the well-known result that its
 cohomology is concentrated in degree zero and equals $\odot^n V$. Thus the Koszul complex is nothing but a representation
of the simplicial cell complex in a particular vector space $V$.

 \sip

\no{\bf 3.2 Permutahedra as a dg properad.}
An $(n\hspace{-0.5mm}-\hspace{-0.5mm} 1)$-dimensional {\em permutahedron}, $P_{n-1}$,
 is, by definition, the convex hull of $n!$
  points $(\sigma(1), \sigma(2), \ldots, \sigma(n))$, $\forall \sigma\in \bS_n$,
in $\R^{n}=\{x^1, \ldots, x^{n}\}$.
 To define its cell complex one has to choose an orientation on $P_{n-1}$ which is the same as to choose an orientation on
the hyperplane $\sum_{i=1}^n x^i=n(n+1)/2$ to which $P_{n-1}$ belongs. We induce it from the standard orientation on $\R^{n}$
by requiring that the manifold with boundary defined by the equation $\sum_{i=1}x^i\leq n(n+1)/2$ is naturally oriented.
Let $(C_\bullet(P_{n-1})=\oplus_{k=0}^{n-1}C_{-k}(P_{n-1}), \delta)$ stand for the associated (non-positively graded)
complex of oriented cells of $P_{n-1}$. \
Its $(n\hspace{-1mm}-\hspace{-1mm} k\hspace{-1mm} -\hspace{-1mm} 1)$-dimensional cells, $P_{n-1}^{I_1,\ldots, I_p}$,
 are indexed by all possible partitions, $[n]=I_1\sqcup I_2 \sqcup \ldots \sqcup I_k$, of the ordered set $[n]$
into $k$ disjoint ordered nonempty subsets (see \cite{SU}).
The natural action , $(x^1,\ldots,x^n)\rar (x^{\sigma(1)}, \ldots,
x^{\sigma(n)})$, of the permutation group $\bS_n$ on $\R^n$ leaves $P_{n-1}$ invariant, and hence makes the
cell complex $C_\bullet(P_{n-1})$ into an $\bS_n$-module. We obviously have, for example, $C_{1-n}(P_{n-1})=\sgn_n$,
so that we can identify the top cell $P_{n-1}^{[n]}$ {\em as an element of the ${\bS_n}$-module}\, with the $(n,0)$-corolla
$
\xy
*=<3mm,2mm>\txt{}*\frm{-},
*=<3mm,1mm>\txt{}*\frm{-};
<-3.5mm,6.3mm>*{_1}**@{},
<-1mm,6.3mm>*{_2}**@{},
<4mm,6.3mm>*{_n}**@{},
<0mm,1mm>*{};<3mm,5mm>*{}**@{-},
  <0mm,1mm>*{};<-1mm,5mm>*{}**@{-},
  <0mm,1mm>*{};<-3mm,5mm>*{}**@{-},
<1mm,4.6mm>*{.\hspace{-0.4mm}.\hspace{-0.4mm}.}**@{},
\endxy
$
with skewsymmetric output legs, $
\xy
*=<3mm,2mm>\txt{}*\frm{-};
*=<3mm,1mm>\txt{}*\frm{-};
<-3.5mm,6.3mm>*{_1}**@{},
<-1mm,6.3mm>*{_2}**@{},
<4mm,6.3mm>*{_n}**@{},
<0mm,1mm>*{};<3mm,5mm>*{}**@{-},
  <0mm,1mm>*{};<-1mm,5mm>*{}**@{-},
  <0mm,1mm>*{};<-3mm,5mm>*{}**@{-},
<1mm,4.6mm>*{.\hspace{-0.4mm}.\hspace{-0.4mm}.}**@{},
\endxy
=
(-1)^\sigma
\xy
*=<3mm,2mm>\txt{}*\frm{-};
*=<3mm,1mm>\txt{}*\frm{-};
<-6mm,6.3mm>*{_{\sigma(1)}}**@{},
<0mm,6.3mm>*{_{\sigma(2)}}**@{},
<6.8mm,6.3mm>*{_{\sigma(n)}}**@{},
<0mm,1mm>*{};<3mm,5mm>*{}**@{-},
  <0mm,1mm>*{};<-1mm,5mm>*{}**@{-},
  <0mm,1mm>*{};<-3mm,5mm>*{}**@{-},
<1mm,4.6mm>*{.\hspace{-0.4mm}.\hspace{-0.4mm}.}**@{},
\endxy
$
, $\forall \sigma\in \bS_n$. More generally, a  simple analysis (similar to the simplicial case in \S 3.1) of how
the action of $\bS_n$ on $\R^n$ permutes the cells and changes their orientation implies that $C_{1-k-n}(P_{n-1})$ is canonically
isomorphic as an  $\bS_n$-module to
 $$
W_k(n):= \bigoplus_{[n]= I_1\sqcup I_2 \sqcup \ldots \sqcup I_k} \K[\bS_n]\ot_{\bS_{I_1}\times \ldots\times S_{I_k}}
\left(\sgn_{I_1}\ot\ldots \ot \sgn_{I_k}\right)[k+n-1]
$$
and that the cells $P_{n-1}^{I_1,\ldots, I_p}$ can be identified {\em as elements of the ${\bS_n}$-module}\, with the
$(n,0)$-corollas,
$
\xy
 <-2mm,0mm>*{\mbox{$\xy *=<16mm,3mm>\txt{}*\frm{-}\endxy$}};<0mm,0mm>*{}**@{},
 <-2mm,0mm>*{\mbox{$\xy *=<16mm,2mm>\txt{}*\frm{-}\endxy$}};<0mm,0mm>*{}**@{},
  <-10mm,1.5mm>*{};<-12mm,7mm>*{}**@{-},
  <-10mm,1.5mm>*{};<-11mm,7mm>*{}**@{-},
  <-10mm,1.5mm>*{};<-9.5mm,6mm>*{}**@{-},
  <-10mm,1.5mm>*{};<-8mm,7mm>*{}**@{-},
 <-10mm,1.5mm>*{};<-9.5mm,6.6mm>*{.\hspace{-0.4mm}.\hspace{-0.4mm}.}**@{},
 <0mm,0mm>*{};<-6.5mm,3.6mm>*{.\hspace{-0.1mm}.\hspace{-0.1mm}.}**@{},
  <-3mm,1.5mm>*{};<-5mm,7mm>*{}**@{-},
  <-3mm,1.5mm>*{};<-4mm,7mm>*{}**@{-},
  <-3mm,1.5mm>*{};<-2.5mm,6mm>*{}**@{-},
  <-3mm,1.5mm>*{};<-1mm,7mm>*{}**@{-},
 <-3mm,1.5mm>*{};<-2.5mm,6.6mm>*{.\hspace{-0.4mm}.\hspace{-0.4mm}.}**@{},
 <0mm,0mm>*{};<2mm,3.6mm>*{.\hspace{-0.1mm}.\hspace{-0.1mm}.}**@{},
<6mm,1.5mm>*{};<4mm,7mm>*{}**@{-},
  <6mm,1.5mm>*{};<5mm,7mm>*{}**@{-},
  <6mm,1.5mm>*{};<6.5mm,6mm>*{}**@{-},
  <6mm,1.5mm>*{};<8mm,7mm>*{}**@{-},
 <6mm,1.5mm>*{};<6.5mm,6.6mm>*{.\hspace{-0.4mm}.\hspace{-0.4mm}.}**@{},
<0mm,0mm>*{};<-9.5mm,8.2mm>*{^{I_{ 1}}}**@{},
<0mm,0mm>*{};<-3mm,8.2mm>*{^{I_{ i}}}**@{},
<0mm,0mm>*{};<6mm,8.2mm>*{^{I_{ k}}}**@{},
\endxy
$, where legs in each $I_i$-bunch are skewsymmetric and  the labels from $ I_i$ are assumed to be distributed over
them in the increasing order from the left to the right. The boundary operator $\delta: C_{1-k-n}(P_{n-1})\rar
C_{2-k-n}(P_{n-1})$ is given on generators by (cf.\ \cite{SU})

\Beq\label{perm-original-differential}
\delta
\xy
 <-2mm,0mm>*{\mbox{$\xy *=<16mm,3mm>\txt{}*\frm{-}\endxy$}};<0mm,0mm>*{}**@{},
<-2mm,0mm>*{\mbox{$\xy *=<16mm,2mm>\txt{}*\frm{-}\endxy$}};<0mm,0mm>*{}**@{},
  <-10mm,1.5mm>*{};<-12mm,7mm>*{}**@{-},
  <-10mm,1.5mm>*{};<-11mm,7mm>*{}**@{-},
  <-10mm,1.5mm>*{};<-9.5mm,6mm>*{}**@{-},
  <-10mm,1.5mm>*{};<-8mm,7mm>*{}**@{-},
 <-10mm,1.5mm>*{};<-9.5mm,6.6mm>*{.\hspace{-0.4mm}.\hspace{-0.4mm}.}**@{},
 <0mm,0mm>*{};<-6.5mm,3.6mm>*{.\hspace{-0.1mm}.\hspace{-0.1mm}.}**@{},
  <-3mm,1.5mm>*{};<-5mm,7mm>*{}**@{-},
  <-3mm,1.5mm>*{};<-4mm,7mm>*{}**@{-},
  <-3mm,1.5mm>*{};<-2.5mm,6mm>*{}**@{-},
  <-3mm,1.5mm>*{};<-1mm,7mm>*{}**@{-},
 <-3mm,1.5mm>*{};<-2.5mm,6.6mm>*{.\hspace{-0.4mm}.\hspace{-0.4mm}.}**@{},
 <0mm,0mm>*{};<2mm,3.6mm>*{.\hspace{-0.1mm}.\hspace{-0.1mm}.}**@{},
<6mm,1.5mm>*{};<4mm,7mm>*{}**@{-},
  <6mm,1.5mm>*{};<5mm,7mm>*{}**@{-},
  <6mm,1.5mm>*{};<6.5mm,6mm>*{}**@{-},
  <6mm,1.5mm>*{};<8mm,7mm>*{}**@{-},
 <6mm,1.5mm>*{};<6.5mm,6.6mm>*{.\hspace{-0.4mm}.\hspace{-0.4mm}.}**@{},
<0mm,0mm>*{};<-9.5mm,8.2mm>*{^{I_{ 1}}}**@{},
<0mm,0mm>*{};<-3mm,8.2mm>*{^{I_{ i}}}**@{},
<0mm,0mm>*{};<6mm,8.2mm>*{^{I_{ k}}}**@{},
\endxy
=
\sum_{i=1}^k\sum_{I_i=i_i'\sqcup I_i'' \atop |I_i'|, |I_i''|\geq 1}
(-1)^{\var +\sigma_{I_i'\sqcup I_i''}}
\xy
 <0mm,0mm>*{\mbox{$\xy *=<20mm,3mm>\txt{}*\frm{-}\endxy$}};<0mm,0mm>*{}**@{},
<0mm,0mm>*{\mbox{$\xy *=<20mm,2mm>\txt{}*\frm{-}\endxy$}};<0mm,0mm>*{}**@{},
  <-10mm,1.5mm>*{};<-12mm,7mm>*{}**@{-},
  <-10mm,1.5mm>*{};<-11mm,7mm>*{}**@{-},
  <-10mm,1.5mm>*{};<-9.5mm,6mm>*{}**@{-},
  <-10mm,1.5mm>*{};<-8mm,7mm>*{}**@{-},
 <-10mm,1.5mm>*{};<-9.5mm,6.6mm>*{.\hspace{-0.4mm}.\hspace{-0.4mm}.}**@{},
 <0mm,0mm>*{};<-6.5mm,3.6mm>*{.\hspace{-0.1mm}.\hspace{-0.1mm}.}**@{},
  <-3mm,1.5mm>*{};<-5mm,7mm>*{}**@{-},
  <-3mm,1.5mm>*{};<-4mm,7mm>*{}**@{-},
  <-3mm,1.5mm>*{};<-2.5mm,6mm>*{}**@{-},
  <-3mm,1.5mm>*{};<-1mm,7mm>*{}**@{-},
 <-3mm,1.5mm>*{};<-2.5mm,6.6mm>*{.\hspace{-0.4mm}.\hspace{-0.4mm}.}**@{},
  <2mm,1.5mm>*{};<0mm,7mm>*{}**@{-},
  <2mm,1.5mm>*{};<1mm,7mm>*{}**@{-},
  <2mm,1.5mm>*{};<2.5mm,6mm>*{}**@{-},
  <2mm,1.5mm>*{};<4mm,7mm>*{}**@{-},
 <2mm,1.5mm>*{};<2.5mm,6.6mm>*{.\hspace{-0.4mm}.\hspace{-0.4mm}.}**@{},
 <0mm,0mm>*{};<6mm,3.6mm>*{.\hspace{-0.1mm}.\hspace{-0.1mm}.}**@{},
<10mm,1.5mm>*{};<8mm,7mm>*{}**@{-},
  <10mm,1.5mm>*{};<9mm,7mm>*{}**@{-},
  <10mm,1.5mm>*{};<10.5mm,6mm>*{}**@{-},
  <10mm,1.5mm>*{};<12mm,7mm>*{}**@{-},
 <10mm,1.5mm>*{};<10.5mm,6.6mm>*{.\hspace{-0.4mm}.\hspace{-0.4mm}.}**@{},
<0mm,0mm>*{};<-9.5mm,8.2mm>*{^{I_{ 1}}}**@{},
<0mm,0mm>*{};<-3mm,8.2mm>*{^{I'_{ i}}}**@{},
<0mm,0mm>*{};<2mm,8.2mm>*{^{I''_{ i}}}**@{},
<0mm,0mm>*{};<10mm,8.2mm>*{^{I_{ k}}}**@{},
\endxy
\Eeq
 where $\var:=  i+1+I_1+\ldots + I_{i-1}+ I_i'$ and
$(-1)^{\sigma_{I_i'\sqcup I_i''}}$ is the sign of the permutation $[I_i]\rar I_i'\sqcup I_i''$. Thus we proved the following

\sip
\no{\bf 3.2.1 Proposition.} {\em  (i) The Saneblidze-Umble permutahedra cell complex is canonically isomorphic to a
dg free properad, $\cP_\bullet:= \cF_c^\uparrow\langle W\rangle$, generated by an $\bS$-bimodule, $W=\{W(m,n)\}$,
$$
W(m,n):=\left\{
\Ba{lc}
\bigoplus_{k=1}^{m}W_k(m) = \mbox{span}
\langle
\xy
 <-2mm,0mm>*{\mbox{$\xy *=<16mm,3mm>\txt{}*\frm{-}\endxy$}};<0mm,0mm>*{}**@{},
<-2mm,0mm>*{\mbox{$\xy *=<16mm,2mm>\txt{}*\frm{-}\endxy$}};<0mm,0mm>*{}**@{},
  <-10mm,1.5mm>*{};<-12mm,7mm>*{}**@{-},
  <-10mm,1.5mm>*{};<-11mm,7mm>*{}**@{-},
  <-10mm,1.5mm>*{};<-9.5mm,6mm>*{}**@{-},
  <-10mm,1.5mm>*{};<-8mm,7mm>*{}**@{-},
 <-10mm,1.5mm>*{};<-9.5mm,6.6mm>*{.\hspace{-0.4mm}.\hspace{-0.4mm}.}**@{},
 <0mm,0mm>*{};<-6.5mm,3.6mm>*{.\hspace{-0.1mm}.\hspace{-0.1mm}.}**@{},
  <-3mm,1.5mm>*{};<-5mm,7mm>*{}**@{-},
  <-3mm,1.5mm>*{};<-4mm,7mm>*{}**@{-},
  <-3mm,1.5mm>*{};<-2.5mm,6mm>*{}**@{-},
  <-3mm,1.5mm>*{};<-1mm,7mm>*{}**@{-},
 <-3mm,1.5mm>*{};<-2.5mm,6.6mm>*{.\hspace{-0.4mm}.\hspace{-0.4mm}.}**@{},
 <0mm,0mm>*{};<2mm,3.6mm>*{.\hspace{-0.1mm}.\hspace{-0.1mm}.}**@{},
<6mm,1.5mm>*{};<4mm,7mm>*{}**@{-},
  <6mm,1.5mm>*{};<5mm,7mm>*{}**@{-},
  <6mm,1.5mm>*{};<6.5mm,6mm>*{}**@{-},
  <6mm,1.5mm>*{};<8mm,7mm>*{}**@{-},
 <6mm,1.5mm>*{};<6.5mm,6.6mm>*{.\hspace{-0.4mm}.\hspace{-0.4mm}.}**@{},
<0mm,0mm>*{};<-9.5mm,8.2mm>*{^{I_{ 1}}}**@{},
<0mm,0mm>*{};<-3mm,8.2mm>*{^{I_{ i}}}**@{},
<0mm,0mm>*{};<6mm,8.2mm>*{^{I_{ k}}}**@{},
\endxy
\rangle
 & \ \ \mbox{for}\ n=0,  \vspace{3mm}\\
0 & \ \ \mbox{for}\  n\geq 1,
\Ea
\right.
$$
and equipped with the differential given on the generators by (\ref{perm-original-differential}).

(ii) The cohomology of $(\cP_\bullet, \p)$ is concentrated in degree zero and equals a free
properad generated by the following degree
zero graphs,
$$
\xy
 <0mm,0mm>*{\blacksquare};
<0mm,0mm>*{};<-4mm,4mm>*{}**@{-},
<0mm,0mm>*{};<-2mm,4mm>*{}**@{-},
<0mm,0mm>*{};<4mm,4mm>*{}**@{-},
<1mm,4mm>*{...};
 <-4.5mm,5.6mm>*{_1};
<-2mm,5.6mm>*{_2};
<4.5mm,5.6mm>*{_m};
\endxy
:=
\sum_{\sigma\in \bS_m}
\xy
 <-2mm,0mm>*{\mbox{$\xy *=<16mm,3mm>\txt{}*\frm{-}\endxy$}};<0mm,0mm>*{}**@{},
 <-2mm,0mm>*{\mbox{$\xy *=<16mm,2mm>\txt{}*\frm{-}\endxy$}};<0mm,0mm>*{}**@{},
  <-10mm,1.5mm>*{};<-11mm,7mm>*{}**@{-},
 <0mm,0mm>*{};<-6.5mm,3.6mm>*{.\hspace{-0.1mm}.\hspace{-0.1mm}.}**@{},
  <-3mm,1.5mm>*{};<-3mm,7mm>*{}**@{-},
 <0mm,0mm>*{};<2mm,3.6mm>*{.\hspace{-0.1mm}.\hspace{-0.1mm}.}**@{},
  <6mm,1.5mm>*{};<7mm,7mm>*{}**@{-},
<0mm,0mm>*{};<-14.5mm,8.2mm>*{^{I_1=\sigma(1)}}**@{},
<0mm,0mm>*{};<-3mm,8.2mm>*{^{I_i=\sigma(i)}}**@{},
<0mm,0mm>*{};<10mm,8.2mm>*{^{I_m=\sigma(m)}}**@{},
\endxy
,
\ \ m\geq 1.
$$

}
\sip

Claim~3.2.1(ii) follows from the contractibility  of permutahedra, and, for each $n$,  the above
graph 
represents the sum
of all vertices of $P_{n-1}$.

\sip

\no{\bf 3.2.2\, From permutahedra to a cobar construction.}  Baranovsky made in \cite{Ba} a remarkable observation
that the permutahedra
cell complex can be used to compute the cohomology of the cobar construction, $\Omega(\wedge^{\bullet} V)$, where
$\wedge^{\bullet} V\simeq \odot^{\bullet}(V[1])$ is interpreted as a graded commutative coalgebra
generated by a vector space $V[1]$.
 In our approach this result
follows immediately from the following two observations: (i) the graded space $\Omega(\wedge^{\bullet} V)$
can be identified with
 $\K\oplus \mbox{Rep}_V(\cP_\bullet)$, where $\mbox{Rep}_V(\cP_\bullet)$ is the space
 of all possible representations,  $\rho: \cP_\bullet\rar \End_V$; the differential in
$\Omega(\wedge^{\bullet} V)$
induced from $\cP_\bullet$ by the formula $d\rho:= \rho\circ \delta$ is precisely the differential of the cobar construction.
Then Proposition~3.2.2(ii) and isomorphism (\ref{main}) imply that the cohomology, $H(\Omega(\wedge^\bullet V))$,
of the cobar construction  equals
$\K\oplus \odot^{\bullet\geq 1}V=\odot^\bullet V$.

\sip

\no{\bf 3.2.3\, Permutahedra cochain complex.} Exactly in the same way as in \S 3.2  one can construct a dg properad,
$(\cP^\bullet, \delta)$, out of the permutahedra {\em cochain}\, complex, $C^\bullet(P_{n-1}):= \Hom_\K(C_{\bullet}(P_{n-1}), \K)$,
with the differential $\delta$ dual to the one given in (\ref{perm-original-differential}). Very remarkably, Chapoton has shown
in \cite{Ch} that one can make the $\bS$-module $\{C^\bullet(P_{n-1})\}$ into a dg {\em operad}, and that this operad is a
quadratic one.
 We do not use this interesting fact in our paper and continue interpreting instead permutahedra as a dg
{\em properad}\, from  which we shall build below more complicated dg props with nice geometric and/or algebraic meaning.
Let us apply first the parity chain functor, $\Pi$, to the dg properad $\cP^\bullet$ (see \S 2.8iii). The result, $\cP := \Pi \cP^\bullet$,
is a dg free properad,
$\cF_c^\uparrow \langle  \check{W} \rangle$, generated by an $\bS$-bimodule, $\check{W}=\{\check{W}(m,n)\}$
with $\check{W}(m,n)=0$ for $n\neq 0$ and with $\check{W}(m,0)$ equal to
\Beq\label{W}
Y(m):=
\bigoplus_{k=1}^{m}
\bigoplus_{[m]= I_1\sqcup I_2 \sqcup \ldots \sqcup I_k \atop |I_\bullet|\geq 1} \K[\bS_n]\ot_{\bS_{I_1}\times \ldots\times S_{I_k}}
\left(\id_{I_1}\ot\ldots \ot \id_{I_k}\right)[k].
\vspace{-5mm}
\Eeq
If we represent the generators of $\cP$ by corollas
$
\xy
 <-2mm,0mm>*{\mbox{$\xy *=<16mm,3mm>\txt{}*\frm{-}\endxy$}};<0mm,0mm>*{}**@{},
  <-10mm,1.5mm>*{};<-12mm,7mm>*{}**@{-},
  <-10mm,1.5mm>*{};<-11mm,7mm>*{}**@{-},
  <-10mm,1.5mm>*{};<-9.5mm,6mm>*{}**@{-},
  <-10mm,1.5mm>*{};<-8mm,7mm>*{}**@{-},
 <-10mm,1.5mm>*{};<-9.5mm,6.6mm>*{.\hspace{-0.4mm}.\hspace{-0.4mm}.}**@{},
 <0mm,0mm>*{};<-6.5mm,3.6mm>*{.\hspace{-0.1mm}.\hspace{-0.1mm}.}**@{},
  <-3mm,1.5mm>*{};<-5mm,7mm>*{}**@{-},
  <-3mm,1.5mm>*{};<-4mm,7mm>*{}**@{-},
  <-3mm,1.5mm>*{};<-2.5mm,6mm>*{}**@{-},
  <-3mm,1.5mm>*{};<-1mm,7mm>*{}**@{-},
 <-3mm,1.5mm>*{};<-2.5mm,6.6mm>*{.\hspace{-0.4mm}.\hspace{-0.4mm}.}**@{},
 <0mm,0mm>*{};<2mm,3.6mm>*{.\hspace{-0.1mm}.\hspace{-0.1mm}.}**@{},
<6mm,1.5mm>*{};<4mm,7mm>*{}**@{-},
  <6mm,1.5mm>*{};<5mm,7mm>*{}**@{-},
  <6mm,1.5mm>*{};<6.5mm,6mm>*{}**@{-},
  <6mm,1.5mm>*{};<8mm,7mm>*{}**@{-},
 <6mm,1.5mm>*{};<6.5mm,6.6mm>*{.\hspace{-0.4mm}.\hspace{-0.4mm}.}**@{},
<0mm,0mm>*{};<-9.5mm,8.2mm>*{^{I_{ 1}}}**@{},
<0mm,0mm>*{};<-3mm,8.2mm>*{^{I_{ i}}}**@{},
<0mm,0mm>*{};<6mm,8.2mm>*{^{I_{ k}}}**@{},
\endxy
$ with {\em symmetric}\, legs in each $I_i$-bunch, then the induced differential in $\cP$ is given by,
\Beq\label{P-diffl}
 \delta\ \xy
 <0mm,0mm>*{\mbox{$\xy *=<20mm,3mm>\txt{}*\frm{-}\endxy$}};<0mm,0mm>*{}**@{},
  <-10mm,1.5mm>*{};<-12mm,7mm>*{}**@{-},
  <-10mm,1.5mm>*{};<-11mm,7mm>*{}**@{-},
  <-10mm,1.5mm>*{};<-9.5mm,6mm>*{}**@{-},
  <-10mm,1.5mm>*{};<-8mm,7mm>*{}**@{-},
 <-10mm,1.5mm>*{};<-9.5mm,6.6mm>*{.\hspace{-0.4mm}.\hspace{-0.4mm}.}**@{},
 <0mm,0mm>*{};<-6.4mm,3.6mm>*{.\hspace{-0.1mm}.\hspace{-0.1mm}.}**@{},
  <-3mm,1.5mm>*{};<-5mm,7mm>*{}**@{-},
  <-3mm,1.5mm>*{};<-4mm,7mm>*{}**@{-},
  <-3mm,1.5mm>*{};<-2.5mm,6mm>*{}**@{-},
  <-3mm,1.5mm>*{};<-1mm,7mm>*{}**@{-},
 <-3mm,1.5mm>*{};<-2.5mm,6.6mm>*{.\hspace{-0.4mm}.\hspace{-0.4mm}.}**@{},
  <2mm,1.5mm>*{};<0mm,7mm>*{}**@{-},
  <2mm,1.5mm>*{};<1mm,7mm>*{}**@{-},
  <2mm,1.5mm>*{};<2.5mm,6mm>*{}**@{-},
  <2mm,1.5mm>*{};<4mm,7mm>*{}**@{-},
 <2mm,1.5mm>*{};<2.5mm,6.6mm>*{.\hspace{-0.4mm}.\hspace{-0.4mm}.}**@{},
 <0mm,0mm>*{};<6mm,3.6mm>*{.\hspace{-0.1mm}.\hspace{-0.1mm}.}**@{},
<10mm,1.5mm>*{};<8mm,7mm>*{}**@{-},
  <10mm,1.5mm>*{};<9mm,7mm>*{}**@{-},
  <10mm,1.5mm>*{};<10.5mm,6mm>*{}**@{-},
  <10mm,1.5mm>*{};<12mm,7mm>*{}**@{-},
 <10mm,1.5mm>*{};<10.5mm,6.6mm>*{.\hspace{-0.4mm}.\hspace{-0.4mm}.}**@{},
<0mm,0mm>*{};<-9.5mm,8.2mm>*{^{I_{ 1}}}**@{},
<0mm,0mm>*{};<-3mm,8.2mm>*{^{I_{ i}}}**@{},
<0mm,0mm>*{};<2mm,8.2mm>*{^{I_{ i+1}}}**@{},
<0mm,0mm>*{};<10mm,8.2mm>*{^{I_{ k}}}**@{},
\endxy =\sum_{i=1}^k(-1)^{i+1}
\xy
 <0mm,0mm>*{\mbox{$\xy *=<20mm,3mm>\txt{}*\frm{-}\endxy$}};<0mm,0mm>*{}**@{},
  <-10mm,1.5mm>*{};<-12mm,7mm>*{}**@{-},
  <-10mm,1.5mm>*{};<-11mm,7mm>*{}**@{-},
  <-10mm,1.5mm>*{};<-9.5mm,6mm>*{}**@{-},
  <-10mm,1.5mm>*{};<-8mm,7mm>*{}**@{-},
 <-10mm,1.5mm>*{};<-9.5mm,6.6mm>*{.\hspace{-0.4mm}.\hspace{-0.4mm}.}**@{},
 <0mm,0mm>*{};<-5.5mm,3.6mm>*{.\hspace{-0.1mm}.\hspace{-0.1mm}.}**@{},
%
  <0mm,1.5mm>*{};<-4mm,7mm>*{}**@{-},
  <0mm,1.5mm>*{};<-2.0mm,6mm>*{}**@{-},
  <0mm,1.5mm>*{};<-1mm,7mm>*{}**@{-},
 <0mm,1.5mm>*{};<-2.3mm,6.6mm>*{.\hspace{-0.4mm}.\hspace{-0.4mm}.}**@{},
%
  <0mm,1.5mm>*{};<1mm,7mm>*{}**@{-},
  <0mm,1.5mm>*{};<2.0mm,6mm>*{}**@{-},
  <0mm,1.5mm>*{};<4mm,7mm>*{}**@{-},
 <0mm,1.5mm>*{};<2.3mm,6.6mm>*{.\hspace{-0.4mm}.\hspace{-0.4mm}.}**@{},
 <0mm,0mm>*{};<6mm,3.6mm>*{.\hspace{-0.1mm}.\hspace{-0.1mm}.}**@{},
<10mm,1.5mm>*{};<8mm,7mm>*{}**@{-},
  <10mm,1.5mm>*{};<9mm,7mm>*{}**@{-},
  <10mm,1.5mm>*{};<10.5mm,6mm>*{}**@{-},
  <10mm,1.5mm>*{};<12mm,7mm>*{}**@{-},
 <10mm,1.5mm>*{};<10.5mm,6.6mm>*{.\hspace{-0.4mm}.\hspace{-0.4mm}.}**@{},
%
<0mm,0mm>*{};<-9.5mm,8.2mm>*{^{I_{ 1}}}**@{},
<0mm,0mm>*{};<-2.5mm,8.2mm>*{{^{I_{ i}\sqcup}}}**@{},
<0mm,0mm>*{};<2.7mm,8.2mm>*{^{I_{ i+1}}}**@{},
<0mm,0mm>*{};<10mm,8.2mm>*{^{I_{ k}}}**@{},
\endxy
\Eeq

\no{\bf 3.2.4 Theorem.} {\em  The cohomology of the dg properad $\cP$ is a free properad generated by the degree $-m$ corollas
with skewsymmetric legs,}\vspace{-0.5mm}
$$\label{P-cohomology}
\xy
 <0mm,0mm>*{\bullet};
<0mm,0mm>*{};<-4mm,4mm>*{}**@{-},
<0mm,0mm>*{};<-2mm,4mm>*{}**@{-},
<0mm,0mm>*{};<4mm,4mm>*{}**@{-},
<1mm,4mm>*{...};
 <-4.5mm,5.6mm>*{_1};
<-2mm,5.6mm>*{_2};
<4.5mm,5.6mm>*{_m};
\endxy
:=
\sum_{\sigma\in \bS_m}(-1)^\sigma
\xy
 <-2mm,0mm>*{\mbox{$\xy *=<16mm,3mm>\txt{}*\frm{-}\endxy$}};<0mm,0mm>*{}**@{},
  <-10mm,1.5mm>*{};<-11mm,7mm>*{}**@{-},
 <0mm,0mm>*{};<-6.5mm,3.6mm>*{.\hspace{-0.1mm}.\hspace{-0.1mm}.}**@{},
  <-3mm,1.5mm>*{};<-3mm,7mm>*{}**@{-},
 <0mm,0mm>*{};<2mm,3.6mm>*{.\hspace{-0.1mm}.\hspace{-0.1mm}.}**@{},
  <6mm,1.5mm>*{};<7mm,7mm>*{}**@{-},
<0mm,0mm>*{};<-14.5mm,8.2mm>*{^{I_1=\sigma(1)}}**@{},
<0mm,0mm>*{};<-3mm,8.2mm>*{^{I_i=\sigma(i)}}**@{},
<0mm,0mm>*{};<10mm,8.2mm>*{^{I_m=\sigma(m)}}**@{},
\endxy
,
\ \ m\geq 1. \vspace{-1mm}
$$
\no{\bf Proof.} By  exactness of the functor $\Pi$, the statement follows  from Proposition~3.2.1(ii).
\hfill $\Box$

\sip

\no{\bf 3.2.5 Corollary.} {\em The cohomology of the bar construction, $B(\odot^\bullet V)$, of the graded
commutative algebra
generated by a vector space $V$, is equal to $\wedge^\bullet V$.}

\sip

\no {\bf Proof.} By definition (see, e.g., \cite{Ba}), $B(\odot^\bullet V)$ is a free tensor coalgebra,
$\ot^\bullet(\odot^{\bullet\geq 1}V)[-1]))$ with the differential $d$ induced from the ordinary multiplication in
$\odot^{\bullet\geq 1}V$. On the other hand, it is easy to see that as a vector space $B(\odot^\bullet V)$
can be identified with $\K\oplus\mbox{Rep}_V(\cP)$,
where $\mbox{Rep}_V(\cP)$ is the space of all possible representations,
$\rho: \cP\rar \End_V$, of the parity shifted properad of permutahedra cochains.
Moreover, the bar differential $d$ is given precisely by $d\rho:= \rho\circ \delta$. Thus
 isomorphism (\ref{main}) and Theorem~3.2.4 imply the required result.
\hfill $\Box$

\sip

\no{\bf 3.3  From permutahedra to polydifferential Hochschild complex.} Let us consider a dg $\bS$-module,
$D=\{D(m,n)\}_{m\geq 1, n\geq 0}$, with $D(m,n):= Y(m)\ot \id_n[-1]$ and with the differential,
$\delta: D(m,n)\rar D(m,n)$ being equal to (\ref{P-diffl}) on the tensor factor $Y(m)$ and identity on the
factor $\id_n[-1]$. Let $(\caD:=\cF^\uparrow\langle D\rangle, \delta)$ be the associated dg free prop.
Its generators can be identified
with  $(m,n)$-corollas\vspace{-1mm}
\Beq\label{corollas-D}
\xy
 <0mm,0mm>*{\mbox{$\xy *=<14mm,3mm>\txt{}*\frm{-}\endxy$}};<0mm,0mm>*{}**@{},
  <-7mm,1.5mm>*{};<-9mm,7mm>*{}**@{-},
  <-7mm,1.5mm>*{};<-8mm,7mm>*{}**@{-},
  <-7mm,1.5mm>*{};<-6.5mm,6mm>*{}**@{-},
  <-7mm,1.5mm>*{};<-5mm,7mm>*{}**@{-},
 <-7mm,1.5mm>*{};<-6.5mm,6.6mm>*{.\hspace{-0.4mm}.\hspace{-0.4mm}.}**@{},
 <0mm,0mm>*{};<-3.5mm,3.6mm>*{.\hspace{-0.1mm}.\hspace{-0.1mm}.}**@{},
  <0mm,1.5mm>*{};<-2mm,7mm>*{}**@{-},
  <0mm,1.5mm>*{};<-1mm,7mm>*{}**@{-},
  <0mm,1.5mm>*{};<0.5mm,6mm>*{}**@{-},
  <0mm,1.5mm>*{};<2mm,7mm>*{}**@{-},
 <0mm,1.5mm>*{};<0.5mm,6.6mm>*{.\hspace{-0.4mm}.\hspace{-0.4mm}.}**@{},
 <0mm,0mm>*{};<3.5mm,3.6mm>*{.\hspace{-0.1mm}.\hspace{-0.1mm}.}**@{},
<7mm,1.5mm>*{};<5mm,7mm>*{}**@{-},
  <7mm,1.5mm>*{};<6mm,7mm>*{}**@{-},
  <7mm,1.5mm>*{};<7.5mm,6mm>*{}**@{-},
  <7mm,1.5mm>*{};<9mm,7mm>*{}**@{-},
 <7mm,1.5mm>*{};<7.5mm,6.6mm>*{.\hspace{-0.4mm}.\hspace{-0.4mm}.}**@{},
<0mm,0mm>*{};<-7.5mm,8.2mm>*{^{I_{ 1}}}**@{},
<0mm,0mm>*{};<0mm,8.2mm>*{^{I_{ i}}}**@{},
<0mm,0mm>*{};<7.5mm,8.2mm>*{^{I_{ k}}}**@{},
 <0mm,0mm>*{};<2mm,-4.6mm>*{\ldots}**@{},
 <0mm,0mm>*{};<4mm,-7mm>*{\ldots}**@{},
<-7mm,-1.5mm>*{};<-8mm,-6mm>*{}**@{-},
 <-4.5mm,-1.5mm>*{};<-5mm,-6mm>*{}**@{-},
 <-2mm,-1.5mm>*{};<-2.2mm,-6mm>*{}**@{-},
<5mm,-1.5mm>*{};<5.5mm,-6mm>*{}**@{-},
  <7mm,-1.5mm>*{};<8mm,-6mm>*{}**@{-},
<0mm,0mm>*{};<-8mm,-7.4mm>*{_1}**@{},
<0mm,0mm>*{};<-5mm,-7.4mm>*{_2}**@{},
<0mm,0mm>*{};<-2mm,-7.4mm>*{_3}**@{},
<0mm,0mm>*{};<8.7mm,-7.4mm>*{_n}**@{},
\endxy
\Eeq
 of degree $1-k$, one such a corolla for every partition $[m]=I_1\sqcup \ldots\sqcup I_k$. The differential $\delta$ in
$\caD$
is then given by
\Beq\label{D-differential}
\delta\left(\xy
 <0mm,0mm>*{\mbox{$\xy *=<20mm,3mm>\txt{}*\frm{-}\endxy$}};<0mm,0mm>*{}**@{},
  <-10mm,1.5mm>*{};<-12mm,7mm>*{}**@{-},
  <-10mm,1.5mm>*{};<-11mm,7mm>*{}**@{-},
  <-10mm,1.5mm>*{};<-9.5mm,6mm>*{}**@{-},
  <-10mm,1.5mm>*{};<-8mm,7mm>*{}**@{-},
 <-10mm,1.5mm>*{};<-9.5mm,6.6mm>*{.\hspace{-0.4mm}.\hspace{-0.4mm}.}**@{},
 <0mm,0mm>*{};<-6.4mm,3.6mm>*{.\hspace{-0.1mm}.\hspace{-0.1mm}.}**@{},
  <-3mm,1.5mm>*{};<-5mm,7mm>*{}**@{-},
  <-3mm,1.5mm>*{};<-4mm,7mm>*{}**@{-},
  <-3mm,1.5mm>*{};<-2.5mm,6mm>*{}**@{-},
  <-3mm,1.5mm>*{};<-1mm,7mm>*{}**@{-},
 <-3mm,1.5mm>*{};<-2.5mm,6.6mm>*{.\hspace{-0.4mm}.\hspace{-0.4mm}.}**@{},
  <2mm,1.5mm>*{};<0mm,7mm>*{}**@{-},
  <2mm,1.5mm>*{};<1mm,7mm>*{}**@{-},
  <2mm,1.5mm>*{};<2.5mm,6mm>*{}**@{-},
  <2mm,1.5mm>*{};<4mm,7mm>*{}**@{-},
 <2mm,1.5mm>*{};<2.5mm,6.6mm>*{.\hspace{-0.4mm}.\hspace{-0.4mm}.}**@{},
 <0mm,0mm>*{};<6mm,3.6mm>*{.\hspace{-0.1mm}.\hspace{-0.1mm}.}**@{},
<10mm,1.5mm>*{};<8mm,7mm>*{}**@{-},
  <10mm,1.5mm>*{};<9mm,7mm>*{}**@{-},
  <10mm,1.5mm>*{};<10.5mm,6mm>*{}**@{-},
  <10mm,1.5mm>*{};<12mm,7mm>*{}**@{-},
 <10mm,1.5mm>*{};<10.5mm,6.6mm>*{.\hspace{-0.4mm}.\hspace{-0.4mm}.}**@{},
 <-10mm,-1.5mm>*{};<-12mm,-6mm>*{}**@{-},
 <-7mm,-1.5mm>*{};<-8mm,-6mm>*{}**@{-},
 <-4mm,-1.5mm>*{};<-4.5mm,-6mm>*{}**@{-},
 <0mm,0mm>*{};<0mm,-4.6mm>*{.\hspace{0.1mm}.\hspace{0.1mm}.}**@{},
<10mm,-1.5mm>*{};<12mm,-6mm>*{}**@{-},
 <7mm,-1.5mm>*{};<8mm,-6mm>*{}**@{-},
  <4mm,-1.5mm>*{};<4.5mm,-6mm>*{}**@{-},
<0mm,0mm>*{};<-9.5mm,8.2mm>*{^{I_{ 1}}}**@{},
<0mm,0mm>*{};<-3mm,8.2mm>*{^{I_{ i}}}**@{},
<0mm,0mm>*{};<2mm,8.2mm>*{^{I_{ i+1}}}**@{},
<0mm,0mm>*{};<10mm,8.2mm>*{^{I_{ k}}}**@{},
<0mm,0mm>*{};<-12mm,-7.4mm>*{_1}**@{},
<0mm,0mm>*{};<-8mm,-7.4mm>*{_2}**@{},
<0mm,0mm>*{};<-4mm,-7.4mm>*{_3}**@{},
<0mm,0mm>*{};<6mm,-7.4mm>*{\ldots}**@{},
<0mm,0mm>*{};<12.5mm,-7.4mm>*{_n}**@{},
\endxy
\right) = \sum_{i=1}^{k-1}(-1)^{i+1}
\xy
 <0mm,0mm>*{\mbox{$\xy *=<20mm,3mm>\txt{}*\frm{-}\endxy$}};<0mm,0mm>*{}**@{},
  <-10mm,1.5mm>*{};<-12mm,7mm>*{}**@{-},
  <-10mm,1.5mm>*{};<-11mm,7mm>*{}**@{-},
  <-10mm,1.5mm>*{};<-9.5mm,6mm>*{}**@{-},
  <-10mm,1.5mm>*{};<-8mm,7mm>*{}**@{-},
 <-10mm,1.5mm>*{};<-9.5mm,6.6mm>*{.\hspace{-0.4mm}.\hspace{-0.4mm}.}**@{},
 <0mm,0mm>*{};<-5.5mm,3.6mm>*{.\hspace{-0.1mm}.\hspace{-0.1mm}.}**@{},
%
  <0mm,1.5mm>*{};<-4mm,7mm>*{}**@{-},
  <0mm,1.5mm>*{};<-2.0mm,6mm>*{}**@{-},
  <0mm,1.5mm>*{};<-1mm,7mm>*{}**@{-},
 <0mm,1.5mm>*{};<-2.3mm,6.6mm>*{.\hspace{-0.4mm}.\hspace{-0.4mm}.}**@{},
%
  <0mm,1.5mm>*{};<1mm,7mm>*{}**@{-},
  <0mm,1.5mm>*{};<2.0mm,6mm>*{}**@{-},
  <0mm,1.5mm>*{};<4mm,7mm>*{}**@{-},
 <0mm,1.5mm>*{};<2.3mm,6.6mm>*{.\hspace{-0.4mm}.\hspace{-0.4mm}.}**@{},
 <0mm,0mm>*{};<6mm,3.6mm>*{.\hspace{-0.1mm}.\hspace{-0.1mm}.}**@{},
<10mm,1.5mm>*{};<8mm,7mm>*{}**@{-},
  <10mm,1.5mm>*{};<9mm,7mm>*{}**@{-},
  <10mm,1.5mm>*{};<10.5mm,6mm>*{}**@{-},
  <10mm,1.5mm>*{};<12mm,7mm>*{}**@{-},
 <10mm,1.5mm>*{};<10.5mm,6.6mm>*{.\hspace{-0.4mm}.\hspace{-0.4mm}.}**@{},
 <-10mm,-1.5mm>*{};<-12mm,-6mm>*{}**@{-},
 <-7mm,-1.5mm>*{};<-8mm,-6mm>*{}**@{-},
 <-4mm,-1.5mm>*{};<-4.5mm,-6mm>*{}**@{-},
 <0mm,0mm>*{};<0mm,-4.6mm>*{.\hspace{0.1mm}.\hspace{0.1mm}.}**@{},
<10mm,-1.5mm>*{};<12mm,-6mm>*{}**@{-},
 <7mm,-1.5mm>*{};<8mm,-6mm>*{}**@{-},
  <4mm,-1.5mm>*{};<4.5mm,-6mm>*{}**@{-},
<0mm,0mm>*{};<-9.5mm,8.2mm>*{^{I_{ 1}}}**@{},
<0mm,0mm>*{};<-2.5mm,8.2mm>*{{^{I_{ i}\sqcup}}}**@{},
<0mm,0mm>*{};<2.7mm,8.2mm>*{^{I_{ i+1}}}**@{},
<0mm,0mm>*{};<10mm,8.2mm>*{^{I_{ k}}}**@{},
<0mm,0mm>*{};<-12mm,-7.4mm>*{_1}**@{},
<0mm,0mm>*{};<-8mm,-7.4mm>*{_2}**@{},
<0mm,0mm>*{};<-4mm,-7.4mm>*{_3}**@{},
<0mm,0mm>*{};<6mm,-7.4mm>*{\ldots}**@{},
<0mm,0mm>*{};<12.5mm,-7.4mm>*{_n}**@{},
\endxy
\Eeq

\no{\bf 3.3.1\,  Proposition.} {\em The cohomology of the dg prop $\caD$ is a free prop, $\cF^\uparrow \langle X\rangle$,
 generated by an $\bS$-bimodule $X=\{X(m,n)\}_{m\geq 1, n\geq 0}$ with $X(m,n):=\sgn_m\ot \id_n[m-1]$}.

\sip

\no{\bf Proof.} By \S 2.8(ii), the functor $\cF^\uparrow$ is exact so that we have $H(\caD)=
\cF^\uparrow\langle H(D)
\rangle$. By K\"unneth theorem, $H(D(m,n))= H(Y(m))\ot \id_n[-1]$. Finally, by Theorem 3.2.4, $H(Y(m))=H(\cP)(m,0)=\sgn_m[m]$.
\hfill $\Box$

\sip

The space of all representations, $\rho: \caD\rar \End_V$, of the prop $\caD$ in a (finite-dimensional) vector space $V$
can be obviously identified with
$$
\mbox{\it Rep}_V(\caD):=\bigoplus_{k\geq 1}\Hom\left(\odot^\bullet V, (\odot^{\bullet\geq 1} V)^{\ot k}\right)[1-k]=
   \bigoplus_{k\geq 1} \Hom(\bar{\f}_V^{\ot k},\f_V)[1-k],
$$
where $\f_V:= \odot^\bullet V^*$ is the graded commutative ring of polynomial functions on the space $V$, and
$\bar{\f}_V:= \odot^{\bullet\geq 1} V^*$ is its subring consisting of functions vanishing at $0$. Differential
(\ref{D-differential}) in the prop $\caD$ induces a differential, $\delta$, in the space $\mbox{Rep}(\caD)$ by the formula,
$
\delta\rho:= \rho\circ \delta, \ \ \forall \rho\in \mbox{Rep}(\caD)$.

\sip

\no{\bf 3.3.2\, Proposition.} {\em The complex $(\mbox{\em Rep}(\caD)_V, \delta)$ is canonically isomorphic
to the polydifferential subcomplex, $(C^{\bullet\geq 1}_{\diff}(\f_V), d_H)$  of the standard Hochschild complex,
$(C^\bullet(\f_V), d_H)$,
of the algebra $\f_V$.}

\sip

\no{\bf Proof}. We shall construct a degree $0$ isomorphism of vector spaces, $i:  \mbox{Rep}_V(\caD) \rar
C^{\bullet \geq 1}_{\diff}(\f_V)$,
such that $i\circ \delta= d_H\circ i$, where $d_H$ stands for the Hochschild  differential.
Let $\{e_\al\}$ be a basis of $V$, and $\{x^\al\}$ the associated dual basis of $V^*$. Any $\rho\in \mbox{Rep}_V(\caD)$
is uniquely determined by its values,
$$
\rho
\left(
\xy
 <0mm,0mm>*{\mbox{$\xy *=<14mm,3mm>\txt{}*\frm{-}\endxy$}};<0mm,0mm>*{}**@{},
  <-7mm,1.5mm>*{};<-9mm,7mm>*{}**@{-},
  <-7mm,1.5mm>*{};<-8mm,7mm>*{}**@{-},
  <-7mm,1.5mm>*{};<-6.5mm,6mm>*{}**@{-},
  <-7mm,1.5mm>*{};<-5mm,7mm>*{}**@{-},
 <-7mm,1.5mm>*{};<-6.5mm,6.6mm>*{.\hspace{-0.4mm}.\hspace{-0.4mm}.}**@{},
 <0mm,0mm>*{};<-3.5mm,3.6mm>*{.\hspace{-0.1mm}.\hspace{-0.1mm}.}**@{},
  <0mm,1.5mm>*{};<-2mm,7mm>*{}**@{-},
  <0mm,1.5mm>*{};<-1mm,7mm>*{}**@{-},
  <0mm,1.5mm>*{};<0.5mm,6mm>*{}**@{-},
  <0mm,1.5mm>*{};<2mm,7mm>*{}**@{-},
 <0mm,1.5mm>*{};<0.5mm,6.6mm>*{.\hspace{-0.4mm}.\hspace{-0.4mm}.}**@{},
 <0mm,0mm>*{};<3.5mm,3.6mm>*{.\hspace{-0.1mm}.\hspace{-0.1mm}.}**@{},
<7mm,1.5mm>*{};<5mm,7mm>*{}**@{-},
  <7mm,1.5mm>*{};<6mm,7mm>*{}**@{-},
  <7mm,1.5mm>*{};<7.5mm,6mm>*{}**@{-},
  <7mm,1.5mm>*{};<9mm,7mm>*{}**@{-},
 <7mm,1.5mm>*{};<7.5mm,6.6mm>*{.\hspace{-0.4mm}.\hspace{-0.4mm}.}**@{},
<-9mm,8.2mm>*{^{1}}**@{},
<-7.2mm,8.2mm>*{^{2}}**@{},
<9mm,8.2mm>*{^{m}}**@{},
 <0mm,0mm>*{};<2mm,-4.6mm>*{\ldots}**@{},
 <0mm,0mm>*{};<1mm,-7mm>*{\ldots}**@{},
 <0mm,0mm>*{};<1mm,8.2mm>*{\ldots}**@{},
<-7mm,-1.5mm>*{};<-8mm,-6mm>*{}**@{-},
 <-4.5mm,-1.5mm>*{};<-5mm,-6mm>*{}**@{-},
 <-2mm,-1.5mm>*{};<-2.2mm,-6mm>*{}**@{-},
<5mm,-1.5mm>*{};<5.5mm,-6mm>*{}**@{-},
  <7mm,-1.5mm>*{};<8mm,-6mm>*{}**@{-},
<0mm,0mm>*{};<-8mm,-7.4mm>*{_1}**@{},
<0mm,0mm>*{};<-5mm,-7.4mm>*{_2}**@{},
<0mm,0mm>*{};<8.7mm,-7.4mm>*{_n}**@{},
\endxy
\right)=:
\sum_{I_1,\ldots, I_k, J \atop |I_\bullet|>1, |J|\geq 0}
\Gamma^{I_1,\ldots, I_k}_J x^J \ot e_{I_1}\ot \ldots\ot e_{I_k} \in  \Hom(\bar{\f}_V^{\ot k},\f_V),
\vspace{-1mm}
$$
for some $\Gamma^{I_1,\ldots, I_k}_J \in \K$. Here  the summation runs  over  multi-indices, $I=\al_1\al_2\ldots
\al_{|I|}$, $x^{I}:= x^{\al_1}\odot \ldots \odot x^{\al_{|I|}}$, and  $e_I:=e_{\al_1}\odot\ldots \odot e_{\al_{|I|}}$.
Then the required map $i$ is given explicitly by
\vspace{-1mm}
$$
i(\rho)
:=\sum_{I_1,\ldots, I_k, J}\frac{1}{|J|!|I_1|!\cdots |I_k|!}
\Gamma^{I_1,\ldots, I_k}_J x^J \frac{\p^{|I_1|}}{\p x^{I_1}}\ot \ldots\ot \frac{\p^{|I_k|}}{\p x^{I_k}}
$$
where ${\p^{|I|}}/{\p x^{I}}:= \p^{|I|}/\p x^{\al_1}\ldots
\p x^{\al_{|I|}}$. Now it is an easy calculation to check (using the definition of the Hochschild
differential) that $i\circ \delta= d_H\circ i$.
\hfill $\Box$

\no{\bf 3.3.3 Corollary.} $H(C^{\bullet \geq 1}_{\diff}(\f_V))=\wedge^{\bullet\geq 1}V \ot \odot^\bullet V^*$.

\sip

\no{\bf Proof.} By Proposition~3.3.2, the cohomology $H(C^{\bullet \geq 1}_{\diff}(\f_V))$ of the
polydifferential Hochschild is equal to the cohomology of the complex $(\mbox{Rep}_V(\caD), \delta)$. The latter is
equal, by isomorphism  (\ref{main}), to $\mbox{Rep}_V(H(\caD))$ which in turn is equal, by Proposition~3.3.1,
to $\wedge^{\bullet\geq 1}V \ot \odot^\bullet V^*$.
\hfill $\Box$

\sip

The complex $C^{\bullet}_{\diff}(\f_V)$ is a direct sum, $\f_V\oplus  C^{\bullet \geq 1}_{\diff}(\f_V)$, where $\f_V$
is a trivial subcomplex. Thus Corollary~3.3.3 implies isomorphism (\ref{introduction-KTH-poly}).

\sip

\no{\bf 3.4 Hochschild complex.} Is there a dg prop whose representation complex is
the general (rather than polydifferential) Hochschild complex for polynomial functions? Consider a dg free prop, $\hat{\caD}$,
which is
generated by the same $\bS$-bimodule $D$ as the prop $\caD$ above, but equipped with a different differential, $\hat{\delta}$,
given on the generators by (\ref{D-differential}) and two extra terms,
$$
(\hat{\delta}-\delta)\hspace{-1.5mm}
\xy
 <0mm,0mm>*{\mbox{$\xy *=<14mm,3mm>\txt{}*\frm{-}\endxy$}};<0mm,0mm>*{}**@{},
  <-7mm,1.5mm>*{};<-9mm,7mm>*{}**@{-},
  <-7mm,1.5mm>*{};<-8mm,7mm>*{}**@{-},
  <-7mm,1.5mm>*{};<-6.5mm,6mm>*{}**@{-},
  <-7mm,1.5mm>*{};<-5mm,7mm>*{}**@{-},
 <-7mm,1.5mm>*{};<-6.5mm,6.6mm>*{.\hspace{-0.4mm}.\hspace{-0.4mm}.}**@{},
 <0mm,0mm>*{};<-3.5mm,3.6mm>*{.\hspace{-0.1mm}.\hspace{-0.1mm}.}**@{},
  <0mm,1.5mm>*{};<-2mm,7mm>*{}**@{-},
  <0mm,1.5mm>*{};<-1mm,7mm>*{}**@{-},
  <0mm,1.5mm>*{};<0.5mm,6mm>*{}**@{-},
  <0mm,1.5mm>*{};<2mm,7mm>*{}**@{-},
 <0mm,1.5mm>*{};<0.5mm,6.6mm>*{.\hspace{-0.4mm}.\hspace{-0.4mm}.}**@{},
 <0mm,0mm>*{};<3.5mm,3.6mm>*{.\hspace{-0.1mm}.\hspace{-0.1mm}.}**@{},
<7mm,1.5mm>*{};<5mm,7mm>*{}**@{-},
  <7mm,1.5mm>*{};<6mm,7mm>*{}**@{-},
  <7mm,1.5mm>*{};<7.5mm,6mm>*{}**@{-},
  <7mm,1.5mm>*{};<9mm,7mm>*{}**@{-},
 <7mm,1.5mm>*{};<7.5mm,6.6mm>*{.\hspace{-0.4mm}.\hspace{-0.4mm}.}**@{},
<0mm,0mm>*{};<-7.5mm,8.2mm>*{^{I_{ 1}}}**@{},
<0mm,0mm>*{};<0mm,8.2mm>*{^{I_{ i}}}**@{},
<0mm,0mm>*{};<7.5mm,8.2mm>*{^{I_{ k}}}**@{},
 <0mm,0mm>*{};<2mm,-4.6mm>*{\ldots}**@{},
 <0mm,0mm>*{};<4mm,-7mm>*{\ldots}**@{},
<-7mm,-1.5mm>*{};<-8mm,-6mm>*{}**@{-},
 <-4.5mm,-1.5mm>*{};<-5mm,-6mm>*{}**@{-},
 <-2mm,-1.5mm>*{};<-2.2mm,-6mm>*{}**@{-},
<5mm,-1.5mm>*{};<5.5mm,-6mm>*{}**@{-},
  <7mm,-1.5mm>*{};<8mm,-6mm>*{}**@{-},
<0mm,0mm>*{};<-8mm,-7.4mm>*{_1}**@{},
<0mm,0mm>*{};<-5mm,-7.4mm>*{_2}**@{},
<0mm,0mm>*{};<-2mm,-7.4mm>*{_3}**@{},
<0mm,0mm>*{};<8.7mm,-7.4mm>*{_n}**@{},
\endxy
 =
-\hspace{-4mm} \sum_{[n]=J_1\sqcup J_2\atop |J_1|=|I_1|}
\xy
<-11mm,-5mm>*{};<-11mm,5mm>*{}**@{-},
<-13mm,0mm>*{.\hspace{-0.1mm}.\hspace{-0.1mm}.}**@{},
<-15mm,-5mm>*{};<-15mm,5mm>*{}**@{-},
<-16mm,-5mm>*{};<-16mm,5mm>*{}**@{-},
 <-13.5mm,6.5mm>*{\overbrace{\ \ \  }};
<-13.5mm,8.5mm>*{^{I_1}};
 <-13.5mm,-6.5mm>*{\underbrace{\ \ \  }};
<-13.5mm,-9.1mm>*{_{J_1}};
 <0mm,-6.5mm>*{\underbrace{\hspace{17mm}}};
<0mm,-9.1mm>*{_{J_2}};
 <0mm,0mm>*{\mbox{$\xy *=<14mm,3mm>\txt{}*\frm{-}\endxy$}};
  <-7mm,1.5mm>*{};<-9mm,7mm>*{}**@{-},
  <-7mm,1.5mm>*{};<-8mm,7mm>*{}**@{-},
  <-7mm,1.5mm>*{};<-6.5mm,6mm>*{}**@{-},
  <-7mm,1.5mm>*{};<-5mm,7mm>*{}**@{-},
 <-7mm,1.5mm>*{};<-6.5mm,6.6mm>*{.\hspace{-0.4mm}.\hspace{-0.4mm}.}**@{},
 <0mm,0mm>*{};<-3.5mm,3.6mm>*{.\hspace{-0.1mm}.\hspace{-0.1mm}.}**@{},
  <0mm,1.5mm>*{};<-2mm,7mm>*{}**@{-},
  <0mm,1.5mm>*{};<-1mm,7mm>*{}**@{-},
  <0mm,1.5mm>*{};<0.5mm,6mm>*{}**@{-},
  <0mm,1.5mm>*{};<2mm,7mm>*{}**@{-},
 <0mm,1.5mm>*{};<0.5mm,6.6mm>*{.\hspace{-0.4mm}.\hspace{-0.4mm}.}**@{},
<3.5mm,3.6mm>*{.\hspace{-0.1mm}.\hspace{-0.1mm}.}**@{},
<7mm,1.5mm>*{};<5mm,7mm>*{}**@{-},
  <7mm,1.5mm>*{};<6mm,7mm>*{}**@{-},
  <7mm,1.5mm>*{};<7.5mm,6mm>*{}**@{-},
  <7mm,1.5mm>*{};<9mm,7mm>*{}**@{-},
 <7mm,1.5mm>*{};<7.5mm,6.6mm>*{.\hspace{-0.4mm}.\hspace{-0.4mm}.}**@{},
<0mm,0mm>*{};<-7.5mm,8.2mm>*{^{I_{2}}}**@{},
<0mm,0mm>*{};<0mm,8.2mm>*{^{I_{ i}}}**@{},
<0mm,0mm>*{};<7.5mm,8.2mm>*{^{I_{ k}}}**@{},
 <0mm,0mm>*{};<2mm,-4.6mm>*{\ldots}**@{},
  %
<-7mm,-1.5mm>*{};<-8mm,-5mm>*{}**@{-},
 <-4.5mm,-1.5mm>*{};<-5mm,-5mm>*{}**@{-},
 <-2mm,-1.5mm>*{};<-2.2mm,-5mm>*{}**@{-},
<5mm,-1.5mm>*{};<5.5mm,-5mm>*{}**@{-},
  <7mm,-1.5mm>*{};<8mm,-5mm>*{}**@{-},
\endxy
-\hspace{-2mm} \sum_{[n]=J_1\sqcup J_2\atop |J_2|=|I_k|}\hspace{-1.5mm}
(-1)^k
\xy
<11mm,-5mm>*{};<11mm,5mm>*{}**@{-},
<13mm,0mm>*{.\hspace{-0.1mm}.\hspace{-0.1mm}.}**@{},
<15mm,-5mm>*{};<15mm,5mm>*{}**@{-},
<16mm,-5mm>*{};<16mm,5mm>*{}**@{-},
 <13.5mm,6.5mm>*{\overbrace{\ \ \  }};
<13.5mm,8.5mm>*{^{I_k}};
 <13.5mm,-6.5mm>*{\underbrace{\ \ \  }};
<13.5mm,-9.1mm>*{_{J_2}};
 <0mm,-6.5mm>*{\underbrace{\hspace{17mm}}};
<0mm,-9.1mm>*{_{J_1}};
 <0mm,0mm>*{\mbox{$\xy *=<14mm,3mm>\txt{}*\frm{-}\endxy$}};
  <-7mm,1.5mm>*{};<-9mm,7mm>*{}**@{-},
  <-7mm,1.5mm>*{};<-8mm,7mm>*{}**@{-},
  <-7mm,1.5mm>*{};<-6.5mm,6mm>*{}**@{-},
  <-7mm,1.5mm>*{};<-5mm,7mm>*{}**@{-},
 <-7mm,1.5mm>*{};<-6.5mm,6.6mm>*{.\hspace{-0.4mm}.\hspace{-0.4mm}.}**@{},
 <0mm,0mm>*{};<-3.5mm,3.6mm>*{.\hspace{-0.1mm}.\hspace{-0.1mm}.}**@{},
  <0mm,1.5mm>*{};<-2mm,7mm>*{}**@{-},
  <0mm,1.5mm>*{};<-1mm,7mm>*{}**@{-},
  <0mm,1.5mm>*{};<0.5mm,6mm>*{}**@{-},
  <0mm,1.5mm>*{};<2mm,7mm>*{}**@{-},
 <0mm,1.5mm>*{};<0.5mm,6.6mm>*{.\hspace{-0.4mm}.\hspace{-0.4mm}.}**@{},
<3.5mm,3.6mm>*{.\hspace{-0.1mm}.\hspace{-0.1mm}.}**@{},
<7mm,1.5mm>*{};<5mm,7mm>*{}**@{-},
  <7mm,1.5mm>*{};<6mm,7mm>*{}**@{-},
  <7mm,1.5mm>*{};<7.5mm,6mm>*{}**@{-},
  <7mm,1.5mm>*{};<9mm,7mm>*{}**@{-},
 <7mm,1.5mm>*{};<7.5mm,6.6mm>*{.\hspace{-0.4mm}.\hspace{-0.4mm}.}**@{},
<0mm,0mm>*{};<-7.5mm,8.2mm>*{^{I_{1}}}**@{},
<0mm,0mm>*{};<0mm,8.2mm>*{^{I_{ i}}}**@{},
<0mm,0mm>*{};<7mm,8.2mm>*{^{I_{ k-1}}}**@{},
 <0mm,0mm>*{};<2mm,-4.6mm>*{\ldots}**@{},
  %
<-7mm,-1.5mm>*{};<-8mm,-5mm>*{}**@{-},
 <-4.5mm,-1.5mm>*{};<-5mm,-5mm>*{}**@{-},
 <-2mm,-1.5mm>*{};<-2.2mm,-5mm>*{}**@{-},
<5mm,-1.5mm>*{};<5.5mm,-5mm>*{}**@{-},
  <7mm,-1.5mm>*{};<8mm,-5mm>*{}**@{-},
\endxy
$$
\no{\bf 3.4.1\, Proposition.} $H(\hat{\caD})= H(\caD)=\cF^\uparrow\langle X \rangle$.

\sip

\no{\bf Proof.}
Define a {\em weight}\, of a generating corolla (\ref{corollas-D})
 of the  prop $\hat{\caD}$  to be $\sum_{i=1}^k |I_i|$, and  the {\em weight}, $w(G)$, of a decorated graph $G$ from
$\hat{\caD}$ to be the sum of weights of all constituent corollas of $G$. Then $F_p:=\{span\langle G\rangle |
w(G)\leq p\}_{p\geq 0}$ is a bounded below  exhaustive filtration of the complex $(\hat{\caD},\hat{\delta})$.  By the
classical convergence theorem, the associated spectral sequence $\{E_r, d_r\}_{r\geq 0}$ converges to
$H(\hat{\caD})$. Its $0$th term, $(E_0, d_0)$, is precisely the complex $(\caD, \delta)$. Thus $E_1$, is, by Propsition~3.3.1,
the free prop $\cF^\uparrow\langle X \rangle$ so that $d_1$ vanishes and the spectral sequence degenerates at the first term
completing the proof. \hfill $\Box$

\sip

The complex $(\mbox{Rep}_V(\hat{\caD}), \hat{\delta})$ associated, by \S 2.8(ii),
 to the dg prop $(\hat{\caD}, \hat{\delta})$
is easily seen to be precisely the standard Hochschild complex,
$(C^{\bullet\geq 1}(\bar{\f}_V)=\oplus_{k\geq 1}\Hom(\bar{\f}_V^{\ot k},
\f_V), d_H)$
of the non-unital algebra $\bar{\f}_V$ with coefficients in the unital algebra $\f_V$. Hence Proposition~3.4.1
and isomorphism (\ref{main}) immediately imply that $HC^{\bullet\geq 1}(\bar{\f}_V)=
\wedge^{\bullet \geq 1}V\ot \odot^\bullet V^*$
which in turn implies, with the help of the theory of simplicial modules (see, e.g., Proposition 1.6.5 in \cite{Lo}),
the Hochschild-Kostant-Rosenberg isomorphism (\ref{introduction-KTH}). Using language of dg props, we deduced it,
therefore,
 from the permutahedra cell complex.

\sip

\no{\bf 3.5 From permutahedra to polydifferential Gerstenhaber-Schack complex.} Let us consider a dg $\bS$-module,
$Q=\{Q(m,n)\}_{m\geq 1, n\geq 1}$, with $Q(m,n):= Y(m)\ot Y(n)^*[-2]$ and
the differential, $d=\delta\ot \Id + \Id\ot \delta^*$, where  $\delta$  is given by (\ref{P-diffl}).
 Let $(\cQ:=\cF^\uparrow\langle Q\rangle, d)$ be the associated dg free prop.
Its generators can be identified
with  corollas
$
\xy
 <0mm,0mm>*{\mbox{$\xy *=<14mm,3mm>\txt{}*\frm{-}\endxy$}};<0mm,0mm>*{}**@{},
  <-7mm,1.5mm>*{};<-9mm,7mm>*{}**@{-},
  <-7mm,1.5mm>*{};<-8mm,7mm>*{}**@{-},
  <-7mm,1.5mm>*{};<-6.5mm,6mm>*{}**@{-},
  <-7mm,1.5mm>*{};<-5mm,7mm>*{}**@{-},
 <-7mm,1.5mm>*{};<-6.5mm,6.6mm>*{.\hspace{-0.4mm}.\hspace{-0.4mm}.}**@{},
 <0mm,0mm>*{};<-3.5mm,3.6mm>*{.\hspace{-0.1mm}.\hspace{-0.1mm}.}**@{},
  <0mm,1.5mm>*{};<-2mm,7mm>*{}**@{-},
  <0mm,1.5mm>*{};<-1mm,7mm>*{}**@{-},
  <0mm,1.5mm>*{};<0.5mm,6mm>*{}**@{-},
  <0mm,1.5mm>*{};<2mm,7mm>*{}**@{-},
 <0mm,1.5mm>*{};<0.5mm,6.6mm>*{.\hspace{-0.4mm}.\hspace{-0.4mm}.}**@{},
 <0mm,0mm>*{};<3.5mm,3.6mm>*{.\hspace{-0.1mm}.\hspace{-0.1mm}.}**@{},
<7mm,1.5mm>*{};<5mm,7mm>*{}**@{-},
  <7mm,1.5mm>*{};<6mm,7mm>*{}**@{-},
  <7mm,1.5mm>*{};<7.5mm,6mm>*{}**@{-},
  <7mm,1.5mm>*{};<9mm,7mm>*{}**@{-},
 <7mm,1.5mm>*{};<7.5mm,6.6mm>*{.\hspace{-0.4mm}.\hspace{-0.4mm}.}**@{},
<0mm,0mm>*{};<-7.5mm,8.2mm>*{^{I_{ 1}}}**@{},
<0mm,0mm>*{};<0mm,8.2mm>*{^{I_{ i}}}**@{},
<0mm,0mm>*{};<7.5mm,8.2mm>*{^{I_{ m}}}**@{},
 <-7mm,-1.5mm>*{};<-9mm,-7mm>*{}**@{-},
  <-7mm,-1.5mm>*{};<-8mm,-7mm>*{}**@{-},
  <-7mm,-1.5mm>*{};<-6.5mm,-6mm>*{}**@{-},
  <-7mm,-1.5mm>*{};<-5mm,-7mm>*{}**@{-},
 <-7mm,-1.5mm>*{};<-6.5mm,-6.6mm>*{.\hspace{-0.4mm}.\hspace{-0.4mm}.}**@{},
 <0mm,0mm>*{};<-3.5mm,-3.6mm>*{.\hspace{-0.1mm}.\hspace{-0.1mm}.}**@{},
  <0mm,-1.5mm>*{};<-2mm,-7mm>*{}**@{-},
  <0mm,-1.5mm>*{};<-1mm,-7mm>*{}**@{-},
  <0mm,-1.5mm>*{};<0.5mm,-6mm>*{}**@{-},
  <0mm,-1.5mm>*{};<2mm,-7mm>*{}**@{-},
 <0mm,-1.5mm>*{};<0.5mm,-6.6mm>*{.\hspace{-0.4mm}.\hspace{-0.4mm}.}**@{},
 <0mm,0mm>*{};<3.5mm,-3.6mm>*{.\hspace{-0.1mm}.\hspace{-0.1mm}.}**@{},
<7mm,-1.5mm>*{};<5mm,-7mm>*{}**@{-},
  <7mm,-1.5mm>*{};<6mm,-7mm>*{}**@{-},
  <7mm,-1.5mm>*{};<7.5mm,-6mm>*{}**@{-},
  <7mm,-1.5mm>*{};<9mm,-7mm>*{}**@{-},
 <7mm,-1.5mm>*{};<7.5mm,-6.6mm>*{.\hspace{-0.4mm}.\hspace{-0.4mm}.}**@{},
<0mm,0mm>*{};<-7.5mm,-8.7mm>*{_{J_{ 1}}}**@{},
<0mm,0mm>*{};<0mm,-8.7mm>*{_{J_{ j}}}**@{},
<0mm,0mm>*{};<7.5mm,-8.7mm>*{_{J_{ n}}}**@{},
\endxy
$
of degree $2-m-n$ with symmetric legs in each input and output bunch. The differential $d$ is given on the generators by
$$\label{differential-GS}
d
\xy
 <0mm,0mm>*{\mbox{$\xy *=<20mm,3mm>\txt{}*\frm{-}\endxy$}};<0mm,0mm>*{}**@{},
  <-10mm,1.5mm>*{};<-12mm,7mm>*{}**@{-},
  <-10mm,1.5mm>*{};<-11mm,7mm>*{}**@{-},
  <-10mm,1.5mm>*{};<-9.5mm,6mm>*{}**@{-},
  <-10mm,1.5mm>*{};<-8mm,7mm>*{}**@{-},
 <-10mm,1.5mm>*{};<-9.5mm,6.6mm>*{.\hspace{-0.4mm}.\hspace{-0.4mm}.}**@{},
 <0mm,0mm>*{};<-6.5mm,3.6mm>*{.\hspace{-0.1mm}.\hspace{-0.1mm}.}**@{},
  <-3mm,1.5mm>*{};<-5mm,7mm>*{}**@{-},
  <-3mm,1.5mm>*{};<-4mm,7mm>*{}**@{-},
  <-3mm,1.5mm>*{};<-2.5mm,6mm>*{}**@{-},
  <-3mm,1.5mm>*{};<-1mm,7mm>*{}**@{-},
 <-3mm,1.5mm>*{};<-2.5mm,6.6mm>*{.\hspace{-0.4mm}.\hspace{-0.4mm}.}**@{},
  <2mm,1.5mm>*{};<0mm,7mm>*{}**@{-},
  <2mm,1.5mm>*{};<1mm,7mm>*{}**@{-},
  <2mm,1.5mm>*{};<2.5mm,6mm>*{}**@{-},
  <2mm,1.5mm>*{};<4mm,7mm>*{}**@{-},
 <2mm,1.5mm>*{};<2.5mm,6.6mm>*{.\hspace{-0.4mm}.\hspace{-0.4mm}.}**@{},
 <0mm,0mm>*{};<6mm,3.6mm>*{.\hspace{-0.1mm}.\hspace{-0.1mm}.}**@{},
<10mm,1.5mm>*{};<8mm,7mm>*{}**@{-},
  <10mm,1.5mm>*{};<9mm,7mm>*{}**@{-},
  <10mm,1.5mm>*{};<10.5mm,6mm>*{}**@{-},
  <10mm,1.5mm>*{};<12mm,7mm>*{}**@{-},
 <10mm,1.5mm>*{};<10.5mm,6.6mm>*{.\hspace{-0.4mm}.\hspace{-0.4mm}.}**@{},
%
<0mm,0mm>*{};<-9.5mm,8.2mm>*{^{I_{ 1}}}**@{},
<0mm,0mm>*{};<-3mm,8.2mm>*{^{I_{ i}}}**@{},
<0mm,0mm>*{};<2mm,8.2mm>*{^{I_{ i+1}}}**@{},
<0mm,0mm>*{};<10mm,8.2mm>*{^{I_{ m}}}**@{},
<-10mm,-1.5mm>*{};<-12mm,-7mm>*{}**@{-},
  <-10mm,-1.5mm>*{};<-11mm,-7mm>*{}**@{-},
  <-10mm,-1.5mm>*{};<-9.5mm,-6mm>*{}**@{-},
  <-10mm,-1.5mm>*{};<-8mm,-7mm>*{}**@{-},
 <-10mm,-1.5mm>*{};<-9.5mm,-6.6mm>*{.\hspace{-0.4mm}.\hspace{-0.4mm}.}**@{},
 <0mm,0mm>*{};<-6.5mm,-3.6mm>*{.\hspace{-0.1mm}.\hspace{-0.1mm}.}**@{},
  <-3mm,-1.5mm>*{};<-5mm,-7mm>*{}**@{-},
  <-3mm,-1.5mm>*{};<-4mm,-7mm>*{}**@{-},
  <-3mm,-1.5mm>*{};<-2.5mm,-6mm>*{}**@{-},
  <-3mm,-1.5mm>*{};<-1mm,-7mm>*{}**@{-},
 <-3mm,-1.5mm>*{};<-2.5mm,-6.6mm>*{.\hspace{-0.4mm}.\hspace{-0.4mm}.}**@{},
  <2mm,-1.5mm>*{};<0mm,-7mm>*{}**@{-},
  <2mm,-1.5mm>*{};<1mm,-7mm>*{}**@{-},
  <2mm,-1.5mm>*{};<2.5mm,-6mm>*{}**@{-},
  <2mm,-1.5mm>*{};<4mm,-7mm>*{}**@{-},
 <2mm,-1.5mm>*{};<2.5mm,-6.6mm>*{.\hspace{-0.4mm}.\hspace{-0.4mm}.}**@{},
 <0mm,0mm>*{};<6mm,-3.6mm>*{.\hspace{-0.1mm}.\hspace{-0.1mm}.}**@{},
<10mm,-1.5mm>*{};<8mm,-7mm>*{}**@{-},
  <10mm,-1.5mm>*{};<9mm,-7mm>*{}**@{-},
  <10mm,-1.5mm>*{};<10.5mm,-6mm>*{}**@{-},
  <10mm,-1.5mm>*{};<12mm,-7mm>*{}**@{-},
 <10mm,-1.5mm>*{};<10.5mm,-6.6mm>*{.\hspace{-0.4mm}.\hspace{-0.4mm}.}**@{},
%
<0mm,0mm>*{};<-9.5mm,-9.2mm>*{^{J_{ 1}}}**@{},
<0mm,0mm>*{};<-3mm,-9.2mm>*{^{J_{ j}}}**@{},
<0mm,0mm>*{};<2mm,-9.2mm>*{^{J_{ j+1}}}**@{},
<0mm,0mm>*{};<10mm,-9.2mm>*{^{J_{ n}}}**@{},
\endxy
=\sum_{i=1}^{m-1}(-1)^i
\xy
 <0mm,0mm>*{\mbox{$\xy *=<20mm,3mm>\txt{}*\frm{-}\endxy$}};<0mm,0mm>*{}**@{},
  <-10mm,1.5mm>*{};<-12mm,7mm>*{}**@{-},
  <-10mm,1.5mm>*{};<-11mm,7mm>*{}**@{-},
  <-10mm,1.5mm>*{};<-9.5mm,6mm>*{}**@{-},
  <-10mm,1.5mm>*{};<-8mm,7mm>*{}**@{-},
 <-10mm,1.5mm>*{};<-9.5mm,6.6mm>*{.\hspace{-0.4mm}.\hspace{-0.4mm}.}**@{},
 <0mm,0mm>*{};<-6.5mm,3.6mm>*{.\hspace{-0.1mm}.\hspace{-0.1mm}.}**@{},
  <0mm,1.5mm>*{};<-5mm,7mm>*{}**@{-},
  <0mm,1.5mm>*{};<-4mm,7mm>*{}**@{-},
  <0mm,1.5mm>*{};<-2.5mm,6mm>*{}**@{-},
  <0mm,1.5mm>*{};<-1.3mm,7mm>*{}**@{-},
 <0mm,1.5mm>*{};<-2.5mm,6.6mm>*{.\hspace{-0.4mm}.\hspace{-0.4mm}.}**@{},
  <0mm,1.5mm>*{};<1.3mm,7mm>*{}**@{-},
  <0mm,1.5mm>*{};<2.5mm,6mm>*{}**@{-},
  <0mm,1.5mm>*{};<4mm,7mm>*{}**@{-},
 <0mm,1.5mm>*{};<2.5mm,6.6mm>*{.\hspace{-0.4mm}.\hspace{-0.4mm}.}**@{},
 <0mm,0mm>*{};<6mm,3.6mm>*{.\hspace{-0.1mm}.\hspace{-0.1mm}.}**@{},
<10mm,1.5mm>*{};<8mm,7mm>*{}**@{-},
  <10mm,1.5mm>*{};<9mm,7mm>*{}**@{-},
  <10mm,1.5mm>*{};<10.5mm,6mm>*{}**@{-},
  <10mm,1.5mm>*{};<12mm,7mm>*{}**@{-},
 <10mm,1.5mm>*{};<10.5mm,6.6mm>*{.\hspace{-0.4mm}.\hspace{-0.4mm}.}**@{},
%
<0mm,0mm>*{};<-9.5mm,8.2mm>*{^{I_{ 1}}}**@{},
<0mm,0mm>*{};<-3mm,8.2mm>*{^{I_{ i}}}**@{},
<0mm,0mm>*{};<2mm,8.2mm>*{^{\sqcup I_{ i+1}}}**@{},
<0mm,0mm>*{};<10mm,8.2mm>*{^{I_{ m}}}**@{},
<-10mm,-1.5mm>*{};<-12mm,-7mm>*{}**@{-},
  <-10mm,-1.5mm>*{};<-11mm,-7mm>*{}**@{-},
  <-10mm,-1.5mm>*{};<-9.5mm,-6mm>*{}**@{-},
  <-10mm,-1.5mm>*{};<-8mm,-7mm>*{}**@{-},
 <-10mm,-1.5mm>*{};<-9.5mm,-6.6mm>*{.\hspace{-0.4mm}.\hspace{-0.4mm}.}**@{},
 <0mm,0mm>*{};<-6.5mm,-3.6mm>*{.\hspace{-0.1mm}.\hspace{-0.1mm}.}**@{},
  <-3mm,-1.5mm>*{};<-5mm,-7mm>*{}**@{-},
  <-3mm,-1.5mm>*{};<-4mm,-7mm>*{}**@{-},
  <-3mm,-1.5mm>*{};<-2.5mm,-6mm>*{}**@{-},
  <-3mm,-1.5mm>*{};<-1mm,-7mm>*{}**@{-},
 <-3mm,-1.5mm>*{};<-2.5mm,-6.6mm>*{.\hspace{-0.4mm}.\hspace{-0.4mm}.}**@{},
  <2mm,-1.5mm>*{};<0mm,-7mm>*{}**@{-},
  <2mm,-1.5mm>*{};<1mm,-7mm>*{}**@{-},
  <2mm,-1.5mm>*{};<2.5mm,-6mm>*{}**@{-},
  <2mm,-1.5mm>*{};<4mm,-7mm>*{}**@{-},
 <2mm,-1.5mm>*{};<2.5mm,-6.6mm>*{.\hspace{-0.4mm}.\hspace{-0.4mm}.}**@{},
 <0mm,0mm>*{};<6mm,-3.6mm>*{.\hspace{-0.1mm}.\hspace{-0.1mm}.}**@{},
<10mm,-1.5mm>*{};<8mm,-7mm>*{}**@{-},
  <10mm,-1.5mm>*{};<9mm,-7mm>*{}**@{-},
  <10mm,-1.5mm>*{};<10.5mm,-6mm>*{}**@{-},
  <10mm,-1.5mm>*{};<12mm,-7mm>*{}**@{-},
 <10mm,-1.5mm>*{};<10.5mm,-6.6mm>*{.\hspace{-0.4mm}.\hspace{-0.4mm}.}**@{},
%
<0mm,0mm>*{};<-9.5mm,-9.2mm>*{^{J_{ 1}}}**@{},
<0mm,0mm>*{};<-3mm,-9.2mm>*{^{J_{ j}}}**@{},
<0mm,0mm>*{};<2mm,-9.2mm>*{^{J_{ j+1}}}**@{},
<0mm,0mm>*{};<10mm,-9.2mm>*{^{J_{ n}}}**@{},
\endxy
+
\sum_{j=1}^{n-1}(-1)^j
\xy
 <0mm,0mm>*{\mbox{$\xy *=<20mm,3mm>\txt{}*\frm{-}\endxy$}};<0mm,0mm>*{}**@{},
  <-10mm,1.5mm>*{};<-12mm,7mm>*{}**@{-},
  <-10mm,1.5mm>*{};<-11mm,7mm>*{}**@{-},
  <-10mm,1.5mm>*{};<-9.5mm,6mm>*{}**@{-},
  <-10mm,1.5mm>*{};<-8mm,7mm>*{}**@{-},
 <-10mm,1.5mm>*{};<-9.5mm,6.6mm>*{.\hspace{-0.4mm}.\hspace{-0.4mm}.}**@{},
 <0mm,0mm>*{};<-6.5mm,3.6mm>*{.\hspace{-0.1mm}.\hspace{-0.1mm}.}**@{},
  <-3mm,1.5mm>*{};<-5mm,7mm>*{}**@{-},
  <-3mm,1.5mm>*{};<-4mm,7mm>*{}**@{-},
  <-3mm,1.5mm>*{};<-2.5mm,6mm>*{}**@{-},
  <-3mm,1.5mm>*{};<-1mm,7mm>*{}**@{-},
 <-3mm,1.5mm>*{};<-2.5mm,6.6mm>*{.\hspace{-0.4mm}.\hspace{-0.4mm}.}**@{},
  <2mm,1.5mm>*{};<0mm,7mm>*{}**@{-},
  <2mm,1.5mm>*{};<1mm,7mm>*{}**@{-},
  <2mm,1.5mm>*{};<2.5mm,6mm>*{}**@{-},
  <2mm,1.5mm>*{};<4mm,7mm>*{}**@{-},
 <2mm,1.5mm>*{};<2.5mm,6.6mm>*{.\hspace{-0.4mm}.\hspace{-0.4mm}.}**@{},
 <0mm,0mm>*{};<6mm,3.6mm>*{.\hspace{-0.1mm}.\hspace{-0.1mm}.}**@{},
<10mm,1.5mm>*{};<8mm,7mm>*{}**@{-},
  <10mm,1.5mm>*{};<9mm,7mm>*{}**@{-},
  <10mm,1.5mm>*{};<10.5mm,6mm>*{}**@{-},
  <10mm,1.5mm>*{};<12mm,7mm>*{}**@{-},
 <10mm,1.5mm>*{};<10.5mm,6.6mm>*{.\hspace{-0.4mm}.\hspace{-0.4mm}.}**@{},
%
<0mm,0mm>*{};<-9.5mm,8.2mm>*{^{I_{ 1}}}**@{},
<0mm,0mm>*{};<-3mm,8.2mm>*{^{I_{ i}}}**@{},
<0mm,0mm>*{};<2mm,8.2mm>*{^{I_{ i+1}}}**@{},
<0mm,0mm>*{};<10mm,8.2mm>*{^{I_{ m}}}**@{},
<-10mm,-1.5mm>*{};<-12mm,-7mm>*{}**@{-},
  <-10mm,-1.5mm>*{};<-11mm,-7mm>*{}**@{-},
  <-10mm,-1.5mm>*{};<-9.5mm,-6mm>*{}**@{-},
  <-10mm,-1.5mm>*{};<-8mm,-7mm>*{}**@{-},
 <-10mm,-1.5mm>*{};<-9.5mm,-6.6mm>*{.\hspace{-0.4mm}.\hspace{-0.4mm}.}**@{},
 <0mm,0mm>*{};<-6.5mm,-3.6mm>*{.\hspace{-0.1mm}.\hspace{-0.1mm}.}**@{},
  <0mm,-1.5mm>*{};<-5mm,-7mm>*{}**@{-},
  <0mm,-1.5mm>*{};<-4mm,-7mm>*{}**@{-},
  <0mm,-1.5mm>*{};<-2.5mm,-6mm>*{}**@{-},
  <0mm,-1.5mm>*{};<-1.3mm,-7mm>*{}**@{-},
 <0mm,-1.5mm>*{};<-2.5mm,-6.6mm>*{.\hspace{-0.4mm}.\hspace{-0.4mm}.}**@{},
  <0mm,-1.5mm>*{};<1.3mm,-7mm>*{}**@{-},
  <0mm,-1.5mm>*{};<2.5mm,-6mm>*{}**@{-},
  <0mm,-1.5mm>*{};<4mm,-7mm>*{}**@{-},
 <0mm,-1.5mm>*{};<2.5mm,-6.6mm>*{.\hspace{-0.4mm}.\hspace{-0.4mm}.}**@{},
 <0mm,0mm>*{};<6mm,-3.6mm>*{.\hspace{-0.1mm}.\hspace{-0.1mm}.}**@{},
<10mm,-1.5mm>*{};<8mm,-7mm>*{}**@{-},
  <10mm,-1.5mm>*{};<9mm,-7mm>*{}**@{-},
  <10mm,-1.5mm>*{};<10.5mm,-6mm>*{}**@{-},
  <10mm,-1.5mm>*{};<12mm,-7mm>*{}**@{-},
 <10mm,-1.5mm>*{};<10.5mm,-6.6mm>*{.\hspace{-0.4mm}.\hspace{-0.4mm}.}**@{},
%
<0mm,0mm>*{};<-9.5mm,-9.2mm>*{^{J_{ 1}}}**@{},
<0mm,0mm>*{};<-3mm,-9.2mm>*{^{J_{ j}}}**@{},
<0mm,0mm>*{};<2.2mm,-9.2mm>*{^{\sqcup J_{ j+1}}}**@{},
<0mm,0mm>*{};<10mm,-9.2mm>*{^{J_{ n}}}**@{},
\endxy.
$$
\no{\bf 3.5.1 Proposition.}
{\em $H(\cQ)$ is a free prop generated by an $\bS$-bimodule $\{\sgn_m\ot \sgn_n[m+n-2]\}_{m,n\geq 1}$.}
\sip

{\no\bf Proof}. Use Theorem 3.2.4 and the K\"unneth theorem. \hfill $\Box$

\sip

The complex of representations, $(\mbox{Rep}(\cQ), d)$, 
is isomorphic as a { graded vector
space}
to
$
C^{\bullet,\bullet}(\bar{\f}_V):=\bigoplus_{m,n\geq 1} \Hom(\bar{\f}_V^{\ot m}, \bar{\f}_V^{\ot n})[m+n-2].
$
 The latter has a
well-known {\em Gerstenhaber-Schack differential}\, \cite{GS},
$$
d_{GS}:  \Hom(\bar{\f}_V^{\ot n}, \bar{\f}_V^{\ot m})\stackrel{d_{GS}^1 \oplus d_{GS}^2}{\lon}
  \Hom(\bar{\f}_V^{\ot n+1},
\bar{\f}_V^{\ot m}) \oplus
\Hom(\bar{\f}_V^{\ot n}, \bar{\f}_V^{\ot m+1}),
$$
with $d_{GS}^1$ given on an arbitrary $\Phi\in \Hom(V^{\ot n}, V^{\ot m})$  by
$$
(d_{GS}^1\Phi)(f_0,..., f_n):=-\Delta^{m-1}(f_0)\cdot \Phi(f_1,..., f_n) +
\sum_{i=0}^{n-1} (-1)^if(f_0,..., f_if_{i+1}, ..., f_n)
\vspace{-1.5mm}
$$
$$
\hspace{20mm} + (-1)^{n}\Phi(f_1, f_2, ..., f_{n-1})\cdot\Delta^{m-1}(f_n),  \ \ \
\forall\ f_0, f_1,\ldots, f_n\in \bar{\f}_V,
$$
where the multiplication in $\bar{\f}_V$ is denoted by juxtaposition, the induced multiplication in the algebra
$\bar{\f}_V^{\ot m}$ by
$\cdot$, the comultiplication in $\bar{\f}_V$ by $\Delta$, and
$$
\Delta^{m-1}: (\Delta\ot \Id^{\ot m-2 })\circ (\Delta\ot \Id^{\ot m-3})\circ \ldots \circ \Delta:
\bar{\f}_V \rar \bar{\f}_V^{\ot m},
$$
for $m\geq 2$ while $\Delta^0:=\Id$.
The expression for $d^2_{GS}$ is an obvious dual analogue of the one for $d^1_{GS}$.

It is evident, however, that $(\mbox{Rep}(\cQ), d) \neq (C^{\bullet,\bullet}(\bar{\f}_V), d_{GS})$.
What is then the meaning of the {\em naturally}\, constructed
 complex $(\mbox{Rep}(\cQ), d)$?
\sip

\no{\bf 3.5.2\, Definition-proposition.} {\em Let $C^{\bullet,\bullet}_\diff(\bar{\f}_V)$ be a subspace of
 $C^{\bullet,\bullet}(\bar{\f}_V)$
spanned by {\em polydifferential}\, operators of the form,
$$
\Ba{rccc}
\Phi: & \bar{\f}_V^{\ot m} &  \lon &  \bar{\f}_V^{\ot n}\\
& f_1\ot \ldots \ot f_m & \lon & \Gamma(f_1, \ldots, f_m),
\Ea
$$
with
$
\Phi(f_1, \ldots, f_m)= 
x^{J_1} \ot \ldots \ot x^{J_n}\cdot \Delta^{n-1}\left( \frac{\p^{|I_1|}f_1}{\p x^{I_1}}\right)\cdot\ldots \cdot
 \Delta^{n-1}\left( \frac{\p^{|I_1|}f_m}{\p x^{I_m}}\right)$ for some families of nonempty multi-indexes $I_\bullet$ and $J_\bullet$.
Then $C^{\bullet,\bullet}_\diff(\bar{\f}_V)$ is a subcomplex of the Gerstenhaber-Schack complex
 $(C^{\bullet,\bullet}(\bar{\f}_V), d_{GS})$.
}

\sip

\no{\bf Proof.}
Proving that $C^{\bullet,\bullet}_\diff(\bar{\f}_V)$ is a subcomplex of
$(C^{\bullet,\bullet}(\bar{\f}_V), d^1_{GS})$ is very similar to the Hochschild complex case.
So we omit these details and concentrate instead on showing that $C^{\bullet,\bullet}_\diff(\bar{\f}_V)$ is a subcomplex of
$(C^{\bullet,\bullet}(\bar{\f}_V), d^2_{GS})$. If, for arbitrary $f\in \bar{\f}_V$, we use Sweedler's  notation,
$\Delta f=  \sum f'\ot f''$, for the coproduct in $\bar{\f}_V$, then, for an operator $\Phi$
as above, one has,
\Beqrn
d_{GS}^2\Phi(f_1, ...,f_m)&=&- \sum f_1'\cdots f_m'\ot \Phi(f_1'', ..., f_m'') - \sum_{i=1}^n (-1)^i
\Delta_i \Phi(f_1,...,f_m)\\
&& + (-1)^n\, \sum \Phi(f_1',...,f_m')\ot f_1''\cdots f_m''
\Eeqrn
where $\Delta_i$ means $\Delta$ applied to the $i$-tensor factor in the space of values, $\bar{\f}_V^{\ot n}$, of $\Phi$.
Taking into account the particular structure of  $\Phi$, one can see that $d_{GS}^2\Phi$ is a linear combination of
polydifferential operators  if and only if an equality holds,
$$
\Delta^2 \frac{\p^{|I|} f}{\p x^I}= \sum f'\ot \Delta \frac{\p^{|I|}f''}{\p x^I},
$$
for arbitrary $f\in \bar{\f}_V$ and arbitrary non-empty multi-index $I$. As product and coproduct in
$ \bar{\f}_V$ are consistent, it is enough to check this equality under the assumption that $\dim V=1$
in which case it is straightforward. \hfill $\Box$

\sip

\no{\bf 3.5.3 \, Proposition.} {\em (i) The complexes $(\mbox{Rep}(\cQ), d)$ and
$(C^{\bullet,\bullet}_\diff(\bar{\f}_V), d_{GS})$
are canonically isomorphic. (ii)  $HC^{\bullet,\bullet}_\diff(\bar{\f}_V)=
\wedge^{\bullet \geq 1} V\ot \odot^{\bullet \geq 1} V^*$.
}
\sip

\no{\bf Proof.} (i) Any representation $\rho\in \mbox{Rep}(\cQ)$ is uniquely determined by its values on the
generators,\vspace{-2mm}
$$
\rho\left(
\xy
 <0mm,0mm>*{\mbox{$\xy *=<14mm,3mm>\txt{}*\frm{-}\endxy$}};<0mm,0mm>*{}**@{},
  <-7mm,1.5mm>*{};<-9mm,7mm>*{}**@{-},
  <-7mm,1.5mm>*{};<-8mm,7mm>*{}**@{-},
  <-7mm,1.5mm>*{};<-6.5mm,6mm>*{}**@{-},
  <-7mm,1.5mm>*{};<-5mm,7mm>*{}**@{-},
 <-7mm,1.5mm>*{};<-6.5mm,6.6mm>*{.\hspace{-0.4mm}.\hspace{-0.4mm}.}**@{},
 <0mm,0mm>*{};<-3.5mm,3.6mm>*{.\hspace{-0.1mm}.\hspace{-0.1mm}.}**@{},
  <0mm,1.5mm>*{};<-2mm,7mm>*{}**@{-},
  <0mm,1.5mm>*{};<-1mm,7mm>*{}**@{-},
  <0mm,1.5mm>*{};<0.5mm,6mm>*{}**@{-},
  <0mm,1.5mm>*{};<2mm,7mm>*{}**@{-},
 <0mm,1.5mm>*{};<0.5mm,6.6mm>*{.\hspace{-0.4mm}.\hspace{-0.4mm}.}**@{},
 <0mm,0mm>*{};<3.5mm,3.6mm>*{.\hspace{-0.1mm}.\hspace{-0.1mm}.}**@{},
<7mm,1.5mm>*{};<5mm,7mm>*{}**@{-},
  <7mm,1.5mm>*{};<6mm,7mm>*{}**@{-},
  <7mm,1.5mm>*{};<7.5mm,6mm>*{}**@{-},
  <7mm,1.5mm>*{};<9mm,7mm>*{}**@{-},
 <7mm,1.5mm>*{};<7.5mm,6.6mm>*{.\hspace{-0.4mm}.\hspace{-0.4mm}.}**@{},
<-9.8mm,8.2mm>*{^{{ 1}}}**@{},
<-7.8mm,8.2mm>*{^{{ 2}}}**@{},
<9.3mm,8.2mm>*{^{{ m}}}**@{},
<0mm,8.2mm>*{\cdots}**@{},
<-9.5mm,-8.3mm>*{_{{ 1}}}**@{},
<9.8mm,-8.3mm>*{_{{ n}}}**@{},
<-7.8mm,-8.2mm>*{_{{ 2}}}**@{},
<0mm,-8.8mm>*{\cdots}**@{},
  %
 <-7mm,-1.5mm>*{};<-9mm,-7mm>*{}**@{-},
  <-7mm,-1.5mm>*{};<-8mm,-7mm>*{}**@{-},
  <-7mm,-1.5mm>*{};<-6.5mm,-6mm>*{}**@{-},
  <-7mm,-1.5mm>*{};<-5mm,-7mm>*{}**@{-},
 <-7mm,-1.5mm>*{};<-6.5mm,-6.6mm>*{.\hspace{-0.4mm}.\hspace{-0.4mm}.}**@{},
 <0mm,0mm>*{};<-3.5mm,-3.6mm>*{.\hspace{-0.1mm}.\hspace{-0.1mm}.}**@{},
  <0mm,-1.5mm>*{};<-2mm,-7mm>*{}**@{-},
  <0mm,-1.5mm>*{};<-1mm,-7mm>*{}**@{-},
  <0mm,-1.5mm>*{};<0.5mm,-6mm>*{}**@{-},
  <0mm,-1.5mm>*{};<2mm,-7mm>*{}**@{-},
 <0mm,-1.5mm>*{};<0.5mm,-6.6mm>*{.\hspace{-0.4mm}.\hspace{-0.4mm}.}**@{},
 <0mm,0mm>*{};<3.5mm,-3.6mm>*{.\hspace{-0.1mm}.\hspace{-0.1mm}.}**@{},
<7mm,-1.5mm>*{};<5mm,-7mm>*{}**@{-},
  <7mm,-1.5mm>*{};<6mm,-7mm>*{}**@{-},
  <7mm,-1.5mm>*{};<7.5mm,-6mm>*{}**@{-},
  <7mm,-1.5mm>*{};<9mm,-7mm>*{}**@{-},
 <7mm,-1.5mm>*{};<7.5mm,-6.6mm>*{.\hspace{-0.4mm}.\hspace{-0.4mm}.}**@{},
  %
\endxy
\right)
=\sum \Gamma_{J_1, \ldots, J_l}^{I_1,\ldots, I_k}
x^{J_1}\ot \ldots \ot x^{J_l}\ot e_{I_1}\ot\ldots\ot e_{I_k}.
\vspace{-2mm}
$$
for some $\Gamma_{J_1, \ldots, J_l}^{I_1,\ldots, I_k}\in \K$.
It is a straightforward calculation to check that the map $i: \mbox{Rep}(\cQ)\rar
C^{\bullet,\bullet}_\diff(\bar{\f}_V)$ given by
$$
i(\rho):= \sum  \frac{\Gamma_{J_1, \ldots, J_l}^{I_1,\ldots, I_k}}{|I_1|!\cdots |J_l|!}\,
x^{J_1}\ot \ldots \ot x^{J_l}\cdot
\Delta^{l-1}\left( \frac{\p^{|I_1|}}{\p x^{I_1}}\right)\cdot ... \cdot
 \Delta^{l-1}\left( \frac{\p^{|I_l|}}{\p x^{I_l}}\right)
$$
satisfy the condition $\rho\circ d= d_{GS}\circ \rho$. Now 3.5.3(ii) follows immediately
from isomorphism (\ref{main}) and Proposition 3.5.1. \hfill $\Box$

\sip

\no{\bf 3.6 Gerstenhaber-Schack complex.} It is not hard to guess which dg prop,
$(\hat{\cQ}, \hat{d})$, has the property
that its associated dg space of representations, $(\mbox{Rep}(\hat{\cQ}), \hat{d})$, is exactly the
Gerstenhaber-Schack complex  $(C^{\bullet,\bullet}(\bar{\f}_V), d_{GS})$. As a prop, $\hat{\cQ}$ is, by definition,
 the same as $\cQ$
above, but the differential differs from $d$ by the following four groups of terms,
$$
(\hat{d}-d)\hspace{-1mm}
\xy
 <0mm,0mm>*{\mbox{$\xy *=<8mm,3mm>\txt{}*\frm{-}\endxy$}};<0mm,0mm>*{}**@{},
  <-4mm,1.5mm>*{};<-6mm,7mm>*{}**@{-},
  <-4mm,1.5mm>*{};<-5mm,7mm>*{}**@{-},
  <-4mm,1.5mm>*{};<-3.5mm,6mm>*{}**@{-},
  <-4mm,1.5mm>*{};<-2mm,7mm>*{}**@{-},
 <-3.5mm,6.6mm>*{.\hspace{-0.4mm}.\hspace{-0.4mm}.}**@{},
<0mm,3.6mm>*{.\hspace{-0.1mm}.\hspace{-0.1mm}.}**@{},
<4mm,1.5mm>*{};<2mm,7mm>*{}**@{-},
  <4mm,1.5mm>*{};<3mm,7mm>*{}**@{-},
  <4mm,1.5mm>*{};<4.5mm,6mm>*{}**@{-},
  <4mm,1.5mm>*{};<6mm,7mm>*{}**@{-},
<4.5mm,6.6mm>*{.\hspace{-0.4mm}.\hspace{-0.4mm}.}**@{},
<-4.5mm,8.2mm>*{^{I_{ 1}}}**@{},
<4.5mm,8.2mm>*{^{I_{ k}}}**@{},
 <-4mm,-1.5mm>*{};<-6mm,-7mm>*{}**@{-},
  <-4mm,-1.5mm>*{};<-5mm,-7mm>*{}**@{-},
  <-4mm,-1.5mm>*{};<-3.5mm,-6mm>*{}**@{-},
  <-4mm,-1.5mm>*{};<-2mm,-7mm>*{}**@{-},
 <-3.5mm,-6.6mm>*{.\hspace{-0.4mm}.\hspace{-0.4mm}.}**@{},
<0mm,-3.6mm>*{.\hspace{-0.1mm}.\hspace{-0.1mm}.}**@{},
<4mm,-1.5mm>*{};<2mm,-7mm>*{}**@{-},
  <4mm,-1.5mm>*{};<3mm,-7mm>*{}**@{-},
  <4mm,-1.5mm>*{};<4.5mm,-6mm>*{}**@{-},
  <4mm,-1.5mm>*{};<6mm,-7mm>*{}**@{-},
 <4.5mm,-6.6mm>*{.\hspace{-0.4mm}.\hspace{-0.4mm}.}**@{},
<-4.5mm,-8.7mm>*{_{J_{ 1}}}**@{},
<4.5mm,-8.7mm>*{_{J_{ l}}}**@{},
\endxy
=
-\hspace{-4mm} \sum_{J_\bullet=J_\bullet'\sqcup J_\bullet''\atop \sum|J_\bullet'|=|I_1|}
\xy
<-9mm,-5mm>*{};<-9mm,5mm>*{}**@{-},
<-10mm,-5mm>*{};<-10mm,5mm>*{}**@{-},
<-11mm,-5mm>*{};<-11mm,5mm>*{}**@{-},
<-13mm,0mm>*{.\hspace{-0.1mm}.\hspace{-0.1mm}.}**@{},
<-15mm,-5mm>*{};<-15mm,5mm>*{}**@{-},
<-16mm,-5mm>*{};<-16mm,5mm>*{}**@{-},
<-17mm,-5mm>*{};<-17mm,5mm>*{}**@{-},
 <-13mm,6.5mm>*{\overbrace{\ \ \ \ \ \ \ }};
<-13mm,8.5mm>*{^{I_1}};
<-13mm,-6.5mm>*{{\underbrace{\ \ \ \ \ \ \ }}};
<-13mm,-9.1mm>*{_{J_1'\sqcup \ldots \sqcup J_l'}};
 <0mm,0mm>*{\mbox{$\xy *=<8mm,3mm>\txt{}*\frm{-}\endxy$}};<0mm,0mm>*{}**@{},
  <-4mm,1.5mm>*{};<-6mm,7mm>*{}**@{-},
  <-4mm,1.5mm>*{};<-5mm,7mm>*{}**@{-},
  <-4mm,1.5mm>*{};<-3.5mm,6mm>*{}**@{-},
  <-4mm,1.5mm>*{};<-2mm,7mm>*{}**@{-},
 <-3.5mm,6.6mm>*{.\hspace{-0.4mm}.\hspace{-0.4mm}.}**@{},
<0mm,3.6mm>*{.\hspace{-0.1mm}.\hspace{-0.1mm}.}**@{},
<4mm,1.5mm>*{};<2mm,7mm>*{}**@{-},
  <4mm,1.5mm>*{};<3mm,7mm>*{}**@{-},
  <4mm,1.5mm>*{};<4.5mm,6mm>*{}**@{-},
  <4mm,1.5mm>*{};<6mm,7mm>*{}**@{-},
<4.5mm,6.6mm>*{.\hspace{-0.4mm}.\hspace{-0.4mm}.}**@{},
<-4.5mm,8.2mm>*{^{I_{ 2}}}**@{},
<4.5mm,8.2mm>*{^{I_{ k}}}**@{},
 <-4mm,-1.5mm>*{};<-6mm,-7mm>*{}**@{-},
  <-4mm,-1.5mm>*{};<-5mm,-7mm>*{}**@{-},
  <-4mm,-1.5mm>*{};<-3.5mm,-6mm>*{}**@{-},
  <-4mm,-1.5mm>*{};<-2mm,-7mm>*{}**@{-},
 <-3.5mm,-6.6mm>*{.\hspace{-0.4mm}.\hspace{-0.4mm}.}**@{},
<0mm,-3.6mm>*{.\hspace{-0.1mm}.\hspace{-0.1mm}.}**@{},
<4mm,-1.5mm>*{};<2mm,-7mm>*{}**@{-},
  <4mm,-1.5mm>*{};<3mm,-7mm>*{}**@{-},
  <4mm,-1.5mm>*{};<4.5mm,-6mm>*{}**@{-},
  <4mm,-1.5mm>*{};<6mm,-7mm>*{}**@{-},
 <4.5mm,-6.6mm>*{.\hspace{-0.4mm}.\hspace{-0.4mm}.}**@{},
<-4.5mm,-8.9mm>*{_{J''_{ 1}}}**@{},
<4.5mm,-8.9mm>*{_{J''_{ l}}}**@{},
\endxy
-\hspace{-2mm}  \sum_{J_\bullet=J_\bullet'\sqcup J_\bullet''\atop \sum|J_\bullet''|=|I_k|}
\hspace{-1.5mm}
(-1)^k
\xy
<9mm,-5mm>*{};<9mm,5mm>*{}**@{-},
<10mm,-5mm>*{};<10mm,5mm>*{}**@{-},
<11mm,-5mm>*{};<11mm,5mm>*{}**@{-},
<13mm,0mm>*{.\hspace{-0.1mm}.\hspace{-0.1mm}.}**@{},
<15mm,-5mm>*{};<15mm,5mm>*{}**@{-},
<16mm,-5mm>*{};<16mm,5mm>*{}**@{-},
<17mm,-5mm>*{};<17mm,5mm>*{}**@{-},
 <13mm,6.5mm>*{\overbrace{\ \ \ \ \ \ \ }};
<13mm,8.5mm>*{^{I_k}};
<13mm,-6.5mm>*{{\underbrace{\ \ \ \ \ \ \ }}};
<13mm,-9.1mm>*{_{J_1''\sqcup \ldots \sqcup J_l''}};
 <0mm,0mm>*{\mbox{$\xy *=<8mm,3mm>\txt{}*\frm{-}\endxy$}};<0mm,0mm>*{}**@{},
  <-4mm,1.5mm>*{};<-6mm,7mm>*{}**@{-},
  <-4mm,1.5mm>*{};<-5mm,7mm>*{}**@{-},
  <-4mm,1.5mm>*{};<-3.5mm,6mm>*{}**@{-},
  <-4mm,1.5mm>*{};<-2mm,7mm>*{}**@{-},
 <-3.5mm,6.6mm>*{.\hspace{-0.4mm}.\hspace{-0.4mm}.}**@{},
<0mm,3.6mm>*{.\hspace{-0.1mm}.\hspace{-0.1mm}.}**@{},
<4mm,1.5mm>*{};<2mm,7mm>*{}**@{-},
  <4mm,1.5mm>*{};<3mm,7mm>*{}**@{-},
  <4mm,1.5mm>*{};<4.5mm,6mm>*{}**@{-},
  <4mm,1.5mm>*{};<6mm,7mm>*{}**@{-},
<4.5mm,6.6mm>*{.\hspace{-0.4mm}.\hspace{-0.4mm}.}**@{},
<-4.5mm,8.2mm>*{^{I_{ 1}}}**@{},
<4.5mm,8.2mm>*{^{I_{ k-1}}}**@{},
 <-4mm,-1.5mm>*{};<-6mm,-7mm>*{}**@{-},
  <-4mm,-1.5mm>*{};<-5mm,-7mm>*{}**@{-},
  <-4mm,-1.5mm>*{};<-3.5mm,-6mm>*{}**@{-},
  <-4mm,-1.5mm>*{};<-2mm,-7mm>*{}**@{-},
 <-3.5mm,-6.6mm>*{.\hspace{-0.4mm}.\hspace{-0.4mm}.}**@{},
<0mm,-3.6mm>*{.\hspace{-0.1mm}.\hspace{-0.1mm}.}**@{},
<4mm,-1.5mm>*{};<2mm,-7mm>*{}**@{-},
  <4mm,-1.5mm>*{};<3mm,-7mm>*{}**@{-},
  <4mm,-1.5mm>*{};<4.5mm,-6mm>*{}**@{-},
  <4mm,-1.5mm>*{};<6mm,-7mm>*{}**@{-},
 <4.5mm,-6.6mm>*{.\hspace{-0.4mm}.\hspace{-0.4mm}.}**@{},
<-4.1mm,-8.7mm>*{_{J'_{ 1}}}**@{},
<4.5mm,-8.7mm>*{_{J'_{ l}}}**@{},
\endxy\vspace{-1mm}
$$
$$
\ \ \ \ \ \ \
- \sum_{I_\bullet=I_\bullet'\sqcup I_\bullet''\atop \sum|I_\bullet'|=|J_1|}
\xy
<-9mm,-5mm>*{};<-9mm,5mm>*{}**@{-},
<-10mm,-5mm>*{};<-10mm,5mm>*{}**@{-},
<-11mm,-5mm>*{};<-11mm,5mm>*{}**@{-},
<-13mm,0mm>*{.\hspace{-0.1mm}.\hspace{-0.1mm}.}**@{},
<-15mm,-5mm>*{};<-15mm,5mm>*{}**@{-},
<-16mm,-5mm>*{};<-16mm,5mm>*{}**@{-},
<-17mm,-5mm>*{};<-17mm,5mm>*{}**@{-},
 <-13mm,6.5mm>*{\overbrace{\ \ \ \ \ \ \ }};
<-13mm,-8.9mm>*{_{J_1}};
<-13mm,-6.5mm>*{{\underbrace{\ \ \ \ \ \ \ }}};
<-13mm,8.7mm>*{_{I_1'\sqcup \ldots \sqcup I_k'}};
 <0mm,0mm>*{\mbox{$\xy *=<8mm,3mm>\txt{}*\frm{-}\endxy$}};<0mm,0mm>*{}**@{},
  <-4mm,1.5mm>*{};<-6mm,7mm>*{}**@{-},
  <-4mm,1.5mm>*{};<-5mm,7mm>*{}**@{-},
  <-4mm,1.5mm>*{};<-3.5mm,6mm>*{}**@{-},
  <-4mm,1.5mm>*{};<-2mm,7mm>*{}**@{-},
 <-3.5mm,6.6mm>*{.\hspace{-0.4mm}.\hspace{-0.4mm}.}**@{},
<0mm,3.6mm>*{.\hspace{-0.1mm}.\hspace{-0.1mm}.}**@{},
<4mm,1.5mm>*{};<2mm,7mm>*{}**@{-},
  <4mm,1.5mm>*{};<3mm,7mm>*{}**@{-},
  <4mm,1.5mm>*{};<4.5mm,6mm>*{}**@{-},
  <4mm,1.5mm>*{};<6mm,7mm>*{}**@{-},
<4.5mm,6.6mm>*{.\hspace{-0.4mm}.\hspace{-0.4mm}.}**@{},
<-4.5mm,8.2mm>*{^{I_{ 2}''}}**@{},
<4.5mm,8.2mm>*{^{I_{ k}''}}**@{},
 <-4mm,-1.5mm>*{};<-6mm,-7mm>*{}**@{-},
  <-4mm,-1.5mm>*{};<-5mm,-7mm>*{}**@{-},
  <-4mm,-1.5mm>*{};<-3.5mm,-6mm>*{}**@{-},
  <-4mm,-1.5mm>*{};<-2mm,-7mm>*{}**@{-},
 <-3.5mm,-6.6mm>*{.\hspace{-0.4mm}.\hspace{-0.4mm}.}**@{},
<0mm,-3.6mm>*{.\hspace{-0.1mm}.\hspace{-0.1mm}.}**@{},
<4mm,-1.5mm>*{};<2mm,-7mm>*{}**@{-},
  <4mm,-1.5mm>*{};<3mm,-7mm>*{}**@{-},
  <4mm,-1.5mm>*{};<4.5mm,-6mm>*{}**@{-},
  <4mm,-1.5mm>*{};<6mm,-7mm>*{}**@{-},
 <4.5mm,-6.6mm>*{.\hspace{-0.4mm}.\hspace{-0.4mm}.}**@{},
<-4.1mm,-8.9mm>*{_{J_{ 2}}}**@{},
<4.5mm,-8.9mm>*{_{J_{ l}}}**@{},
\endxy
-\hspace{-1mm} \sum_{I_\bullet=I_\bullet'\sqcup I_\bullet''\atop \sum |I_\bullet''|=|J_l|}\hspace{-1.5mm}
(-1)^l
\xy
<9mm,-5mm>*{};<9mm,5mm>*{}**@{-},
<10mm,-5mm>*{};<10mm,5mm>*{}**@{-},
<11mm,-5mm>*{};<11mm,5mm>*{}**@{-},
<13mm,0mm>*{.\hspace{-0.1mm}.\hspace{-0.1mm}.}**@{},
<15mm,-5mm>*{};<15mm,5mm>*{}**@{-},
<16mm,-5mm>*{};<16mm,5mm>*{}**@{-},
<17mm,-5mm>*{};<17mm,5mm>*{}**@{-},
 <13mm,6.5mm>*{\overbrace{\ \ \ \ \ \ \ }};
<13mm,-8.7mm>*{_{J_l}};
<13mm,-6.5mm>*{{\underbrace{\ \ \ \ \ \ \ }}};
<13mm,9.1mm>*{_{I_1''\sqcup \ldots \sqcup I_k''}};
 <0mm,0mm>*{\mbox{$\xy *=<8mm,3mm>\txt{}*\frm{-}\endxy$}};<0mm,0mm>*{}**@{},
  <-4mm,1.5mm>*{};<-6mm,7mm>*{}**@{-},
  <-4mm,1.5mm>*{};<-5mm,7mm>*{}**@{-},
  <-4mm,1.5mm>*{};<-3.5mm,6mm>*{}**@{-},
  <-4mm,1.5mm>*{};<-2mm,7mm>*{}**@{-},
 <-3.5mm,6.6mm>*{.\hspace{-0.4mm}.\hspace{-0.4mm}.}**@{},
<0mm,3.6mm>*{.\hspace{-0.1mm}.\hspace{-0.1mm}.}**@{},
<4mm,1.5mm>*{};<2mm,7mm>*{}**@{-},
  <4mm,1.5mm>*{};<3mm,7mm>*{}**@{-},
  <4mm,1.5mm>*{};<4.5mm,6mm>*{}**@{-},
  <4mm,1.5mm>*{};<6mm,7mm>*{}**@{-},
<4.5mm,6.6mm>*{.\hspace{-0.4mm}.\hspace{-0.4mm}.}**@{},
<-4.5mm,8.4mm>*{^{I_{ 1}'}}**@{},
<4.1mm,8.4mm>*{^{I_{ k}'}}**@{},
 <-4mm,-1.5mm>*{};<-6mm,-7mm>*{}**@{-},
  <-4mm,-1.5mm>*{};<-5mm,-7mm>*{}**@{-},
  <-4mm,-1.5mm>*{};<-3.5mm,-6mm>*{}**@{-},
  <-4mm,-1.5mm>*{};<-2mm,-7mm>*{}**@{-},
 <-3.5mm,-6.6mm>*{.\hspace{-0.4mm}.\hspace{-0.4mm}.}**@{},
<0mm,-3.6mm>*{.\hspace{-0.1mm}.\hspace{-0.1mm}.}**@{},
<4mm,-1.5mm>*{};<2mm,-7mm>*{}**@{-},
  <4mm,-1.5mm>*{};<3mm,-7mm>*{}**@{-},
  <4mm,-1.5mm>*{};<4.5mm,-6mm>*{}**@{-},
  <4mm,-1.5mm>*{};<6mm,-7mm>*{}**@{-},
 <4.5mm,-6.6mm>*{.\hspace{-0.4mm}.\hspace{-0.4mm}.}**@{},
<-4.5mm,-8.7mm>*{_{J_{ 1}}}**@{},
<4.3mm,-8.7mm>*{_{J_{ l-1}}}**@{},
\endxy
$$
Using a spectral sequence argument very similar to the one used in the proof of Proposition 3.4.1, one
easily obtains the following

\sip

\no{\bf 3.6.1 Proposition.} $H(\hat{\cQ})= H(\cQ)$.

\sip

\no{\bf 3.6.2 Corollary.} $H(C^{\bullet,\bullet}(\bar{\f}_V))= H(C^{\bullet,\bullet}_\diff(\bar{\f}_V))=
\wedge^{\bullet \geq 1} V\ot \odot^{\bullet \geq 1} V^*$.

\sip

\no The latter result together with the standard results from the theory of simplicial modules \cite{Lo} imply
formula (\ref{introduction-KTH-GS}).

\sip

\no{\bf 3.7 On the Etingof-Kazhdan quantization.} Note that the Gerstenhaber-Schack complex
$C^{\bullet,\bullet}(\bar{\f}_V)$
has a structure of prop, the endomorphism prop of $\bar{\f}_V$. Moreover, it is easy to see that
$C^{\bullet,\bullet}_\diff(\bar{\f}_V)$ is also closed under prop compositions so that the natural inclusion,
$j: C^{\bullet,\bullet}_\diff(\bar{\f}_V)\rar C^{\bullet,\bullet}(\bar{\f}_V)$, is a morphism of props
\cite{Me4}.
A choice of a minimal resolution, $\cA ss\cB_\infty$, of the prop, $\cA ss\cB$, of associative bialgebras,
induces \cite{MV} on $C^{\bullet,\bullet}(\bar{\f}_V)$ (resp.\ on $C^{\bullet,\bullet}_\diff(\bar{\f}_V)$)
the structure of a filtered $L_\infty$-algebra whose Maurer-Cartan elements describe deformations
of the standard bialgebra structure on $\bar{\f}_V$ in the class of (resp.\ polydifferential) strongly homotopy
bialgebra structures. Moreover \cite{MV}, the initial term of this induced  $L_\infty$-structure is precisely
the Gerstenhaber-Schack differential. The inclusion map
$j: C^{\bullet,\bullet}_\diff(\bar{\f}_V)\rar C^{\bullet,\bullet}(\bar{\f}_V)$ extends to a morphism of $L_\infty$-algebras
which, by isomorphisms (\ref{introduction-KTH-GS-poly}) and (\ref{introduction-KTH-GS}), is a quasi-isomorphism.
The Etingof-Kazhdan universal quantization \cite{EK} of (possibly, infinite-dimensional) Lie bialgebra structures
on $V$ associates to such a structure, say $\nu$, a Maurer-Cartan element, $\ga_\nu$, in the $L_\infty$-algebra
$C^{\bullet,\bullet}(\bar{\f}_V)$. As $L_\infty$ quasi-isomorphisms are invertible \cite{Ko}, there is always an associated Maurer-Cartan
element $j^{-1}(\ga^{EK}_\nu)$ which, for degree reasons, describes an associted to $\nu$  {\em polydifferential}\, bialgebra
structure on $\bar{\f}_V$. Thus we proved that  {\em for any Lie bialgebra structure, $\nu$,
 on a vector space $V$ there exists its bialgebra quantization, $j^{-1}(\ga^{EK}_\nu)$, within the class
of polydifferential operators from $\cC^{\bullet,\bullet}_{poly}(\f_V)$.}

\section{Dg prop of unital $A_\infty$-structures}
\no{\bf 4.1 Differential in a free prop.}
A differential in a free prop $\cF^\uparrow\langle E\rangle$ can be decomposed into a sum,
$\delta=\sum_{p\geq 1}\delta_{(p)}$, where
$
\delta^{(p)}: E \stackrel{\delta}{\lon}
\cF^\uparrow\langle E\rangle \stackrel{pr_p}{\lon} \cF^\uparrow_{(p)}\langle E\rangle
$
 is the composition
of $\delta$ with the projection to the subspace,  $\cF^\uparrow_{(p)}\langle E\rangle\subset
\cF^\uparrow\langle E\rangle$, spanned by decorated graphs with precisely $p$ vertices. We studied in \S 3 free props
 equipped with
differentials of the form $\delta=\delta_{(1)}$ which preserve the number of vertices  of decorated graphs, and heavily used
the fact that $\delta$ makes the associated space of
representations, $\mbox{Rep}_V( \cF^\uparrow\langle E\rangle)\simeq \Hom(E, \End_V)$, into a complex whose cohomology
one can easily read from the cohomology of $(\cF^\uparrow\langle E\rangle, \delta)$.
Remarkably \cite{MV}, a generic  differential $\delta$ in  $\cF^\uparrow\langle E\rangle$
makes the vector space $\mbox{Rep}_V( \cF^\uparrow\langle E\rangle)[1]$ into a
$L_\infty$-{\em algebra}\, whose $p$th homotopy Lie bracket is completely determined by $p$-th summand, $\delta_{(p)}$, of the
differential $\delta$. In particular, if $\delta$ has $\delta_{(p)}=0$ for all $p\geq 3$, then
$\mbox{Rep}_V( \cF^\uparrow\langle E\rangle)[1]$ is canonically a dg Lie algebra with the differential determined
by $\delta_{(1)}$ and Lie brackets determined by $\delta_{(2)}$. Thus, if we want to extend isomorphisms
(\ref{introduction-KTH}) and (\ref{introduction-KTH-poly}) into isomorphisms of Lie algebras, we have to look for more
complicated (than the ones studied in \S 3) dg props  canonically associated with the (polydifferential)
Hochschild complex for $\f_V$.

\sip

\no{\bf 4.2 Dg prop of polyvector fields.} Let $\PV$ be a dg free prop  generated by
the $\bS$-module, $X[-1]=\{X(m,n)[-1]\}_{m\geq 1,n\geq 0}$,  \vspace{-2mm}
$$
X(m,n)[-1]= \sgn_m\ot \id_n[m-2]=\mbox{span} \langle
\begin{xy}
 <0mm,0mm>*{\bullet};<0mm,0mm>*{}**@{},
 <0mm,0mm>*{};<-8mm,5mm>*{}**@{-},
 <0mm,0mm>*{};<-4.5mm,5mm>*{}**@{-},
 <0mm,0mm>*{};<-1mm,5mm>*{\ldots}**@{},
 <0mm,0mm>*{};<4.5mm,5mm>*{}**@{-},
 <0mm,0mm>*{};<8mm,5mm>*{}**@{-},
   <0mm,0mm>*{};<-8.5mm,5.5mm>*{^1}**@{},
   <0mm,0mm>*{};<-5mm,5.5mm>*{^2}**@{},
   <0mm,0mm>*{};<9.0mm,5.5mm>*{^m}**@{},
 <0mm,0mm>*{};<-8mm,-5mm>*{}**@{-},
 <0mm,0mm>*{};<-4.5mm,-5mm>*{}**@{-},
 <0mm,0mm>*{};<-1mm,-5mm>*{\ldots}**@{},
 <0mm,0mm>*{};<4.5mm,-5mm>*{}**@{-},
 <0mm,0mm>*{};<8mm,-5mm>*{}**@{-},
   <0mm,0mm>*{};<-8.5mm,-6.9mm>*{^1}**@{},
   <0mm,0mm>*{};<-5mm,-6.9mm>*{^2}**@{},
   <0mm,0mm>*{};<9.0mm,-6.9mm>*{^n}**@{},
 \end{xy}
\rangle
\vspace{-1mm}
$$
which is obtained from the $\bS$-module $X$ of Proposition~3.3.1 by a degree shift.
The differential in $\PV$ is defined as follows (cf.\ \cite{Me0}),  \vspace{-6mm}
$$
\p \begin{xy}
 <0mm,0mm>*{\bullet};<0mm,0mm>*{}**@{},
 <0mm,0mm>*{};<-8mm,5mm>*{}**@{-},
 <0mm,0mm>*{};<-4.5mm,5mm>*{}**@{-},
 <0mm,0mm>*{};<-1mm,5mm>*{\ldots}**@{},
 <0mm,0mm>*{};<4.5mm,5mm>*{}**@{-},
 <0mm,0mm>*{};<8mm,5mm>*{}**@{-},
   <0mm,0mm>*{};<-8.5mm,5.5mm>*{^1}**@{},
   <0mm,0mm>*{};<-5mm,5.5mm>*{^2}**@{},
   <0mm,0mm>*{};<9.0mm,5.5mm>*{^m}**@{},
 <0mm,0mm>*{};<-8mm,-5mm>*{}**@{-},
 <0mm,0mm>*{};<-4.5mm,-5mm>*{}**@{-},
 <0mm,0mm>*{};<-1mm,-5mm>*{\ldots}**@{},
 <0mm,0mm>*{};<4.5mm,-5mm>*{}**@{-},
 <0mm,0mm>*{};<8mm,-5mm>*{}**@{-},
   <0mm,0mm>*{};<-8.5mm,-6.9mm>*{^1}**@{},
   <0mm,0mm>*{};<-5mm,-6.9mm>*{^2}**@{},
   <0mm,0mm>*{};<9.0mm,-6.9mm>*{^n}**@{},
 \end{xy} =
 \sum_{[m]=I_1\sqcup I_2\atop {[n]=J_1\sqcup J_2\atop
 {|I_1|\geq 0, |I_2|\geq 1 \atop
 |J_1|\geq 1, |J_2|\geq 0}}
}\hspace{0mm} (-1)^{\sigma(I_1\sqcup I_2) + |I_1|(|I_2|+1)}
 \begin{xy}
 <0mm,0mm>*{\bullet};<0mm,0mm>*{}**@{},
 <0mm,0mm>*{};<-8mm,5mm>*{}**@{-},
 <0mm,0mm>*{};<-4.5mm,5mm>*{}**@{-},
 <0mm,0mm>*{};<0mm,5mm>*{\ldots}**@{},
 <0mm,0mm>*{};<4.5mm,5mm>*{}**@{-},
 <0mm,0mm>*{};<13mm,5mm>*{}**@{-},
     <0mm,0mm>*{};<-2mm,7mm>*{\overbrace{\ \ \ \ \ \ \ \ \ \ \ \ }}**@{},
     <0mm,0mm>*{};<-2mm,9mm>*{^{I_1}}**@{},
 <0mm,0mm>*{};<-8mm,-5mm>*{}**@{-},
 <0mm,0mm>*{};<-4.5mm,-5mm>*{}**@{-},
 <0mm,0mm>*{};<-1mm,-5mm>*{\ldots}**@{},
 <0mm,0mm>*{};<4.5mm,-5mm>*{}**@{-},
 <0mm,0mm>*{};<8mm,-5mm>*{}**@{-},
      <0mm,0mm>*{};<0mm,-7mm>*{\underbrace{\ \ \ \ \ \ \ \ \ \ \ \ \ \ \
      }}**@{},
      <0mm,0mm>*{};<0mm,-10.6mm>*{_{J_1}}**@{},
 <13mm,5mm>*{};<13mm,5mm>*{\bullet}**@{},
 <13mm,5mm>*{};<5mm,10mm>*{}**@{-},
 <13mm,5mm>*{};<8.5mm,10mm>*{}**@{-},
 <13mm,5mm>*{};<13mm,10mm>*{\ldots}**@{},
 <13mm,5mm>*{};<16.5mm,10mm>*{}**@{-},
 <13mm,5mm>*{};<20mm,10mm>*{}**@{-},
      <13mm,5mm>*{};<13mm,12mm>*{\overbrace{\ \ \ \ \ \ \ \ \ \ \ \ \ \ }}**@{},
      <13mm,5mm>*{};<13mm,14mm>*{^{I_2}}**@{},
 <13mm,5mm>*{};<8mm,0mm>*{}**@{-},
 <13mm,5mm>*{};<12mm,0mm>*{\ldots}**@{},
 <13mm,5mm>*{};<16.5mm,0mm>*{}**@{-},
 <13mm,5mm>*{};<20mm,0mm>*{}**@{-},
     <13mm,5mm>*{};<14.3mm,-2mm>*{\underbrace{\ \ \ \ \ \ \ \ \ \ \ }}**@{},
     <13mm,5mm>*{};<14.3mm,-4.5mm>*{_{J_2}}**@{},
 \end{xy}
$$
where $\sigma(I_1\sqcup I_2)$ is the sign of the permutation $[n]\rar I_1\sqcup I_2$. This differential is quadratic,
$\p=\p_{(2)}$, so that, according to the general theory (see Theorem 60 in \cite{MV}), the space
 $\mbox{Rep}(\PV)_V[1]=
\wedge^{\bullet \geq 1}V \ot \odot^\bullet V \simeq \wedge^{\bullet\geq 1} \cT_V$  comes equipped with a Lie algebra structure
which, as it is not hard to check (cf.\ \cite{Me0}), is precisely the Schouten bracket.

\sip

\no{\bf 4.3 Unital $A_\infty$-structures on $\f_V$.} It is well-known that the vector space
 $\bar{C}^\bullet(\f_V):=\oplus_{k\geq 1} \Hom(\f_V^{\ot k}, \f_V)[1-k]$ has a natural graded Lie algebra structure
with respect to the Gerstenhaber brackets, $[\ ,\ ]_G$. By definition, an $A_\infty$-{\em algebra structure}\,
on the space $\f_V$ is a Maurer-Cartan element in this Lie algebra, that is, a total degree 1 element
$\Ga\in \bar{C}^\bullet(\f_V)$ such that $[\Ga, \Ga]_G=0$. Such an element, $\Ga$, is equivalent to a sequence
of homogeneous linear maps, $\{\Ga_k: \f_V^{\ot k}\rar \f_V[2-k]\}$ satisfying a sequence of quadratic equations
(cf.\ \cite{St}).
An $A_\infty$-algebra structure is called {\em unital}\, if, for every $k\geq 3$, the map $\Ga_k$ factors
 through the composition
$\f_V^{\ot k}\rar \bar{\f}_V^{\ot k}\rar \f_V[k-2]$ and $\Ga_2(1,f)=\Ga_2(f,1)=f$.
The following lemma is obvious.

\sip

\no{\bf 4.3.1 Lemma.} {\em There is a one-to-one correspondence between unital $A_\infty$-structures on $\f_V$
and Maurer-Cartan elements, \vspace{-1mm}
$$
\{\Ga\in \bar{C}^\bullet(\bar{\f}_V):\ |\Ga|=1\ \mbox{and}\ d_H\Ga+ \frac{1}{2}[\Ga,\Ga]_G=0\},
\vspace{-1mm}
$$
in the Hochschild dg Lie algebra, $\bar{C}^\bullet(\bar{\f}_V)$, for the ring $\bar{\f}_V\subset \f_V$.
}
\sip

\no{\bf 4.4\,  Dg prop of unital $A_\infty$-structures.} Consider a dg free prop, $(\DefQ, d)$,
generated by corollas (\ref{corollas-D}) (to which we assign now degree $2-k$)
and equipped with the differential
given by\footnote{We have to assume
that  $\DefQ$ is completed with respect to the genus filtration.}(cf.\ \cite{Me1}) \vspace{-0.5mm}
$$
d
\xy
 <0mm,0mm>*{\mbox{$\xy *=<20mm,3mm>\txt{}*\frm{-}\endxy$}};<0mm,0mm>*{}**@{},
  <-10mm,1.5mm>*{};<-12mm,7mm>*{}**@{-},
  <-10mm,1.5mm>*{};<-11mm,7mm>*{}**@{-},
  <-10mm,1.5mm>*{};<-9.5mm,6mm>*{}**@{-},
  <-10mm,1.5mm>*{};<-8mm,7mm>*{}**@{-},
 <-10mm,1.5mm>*{};<-9.5mm,6.6mm>*{.\hspace{-0.4mm}.\hspace{-0.4mm}.}**@{},
 <0mm,0mm>*{};<-6.4mm,3.6mm>*{.\hspace{-0.1mm}.\hspace{-0.1mm}.}**@{},
  <-3mm,1.5mm>*{};<-5mm,7mm>*{}**@{-},
  <-3mm,1.5mm>*{};<-4mm,7mm>*{}**@{-},
  <-3mm,1.5mm>*{};<-2.5mm,6mm>*{}**@{-},
  <-3mm,1.5mm>*{};<-1mm,7mm>*{}**@{-},
 <-3mm,1.5mm>*{};<-2.5mm,6.6mm>*{.\hspace{-0.4mm}.\hspace{-0.4mm}.}**@{},
  <2mm,1.5mm>*{};<0mm,7mm>*{}**@{-},
  <2mm,1.5mm>*{};<1mm,7mm>*{}**@{-},
  <2mm,1.5mm>*{};<2.5mm,6mm>*{}**@{-},
  <2mm,1.5mm>*{};<4mm,7mm>*{}**@{-},
 <2mm,1.5mm>*{};<2.5mm,6.6mm>*{.\hspace{-0.4mm}.\hspace{-0.4mm}.}**@{},
 <0mm,0mm>*{};<6mm,3.6mm>*{.\hspace{-0.1mm}.\hspace{-0.1mm}.}**@{},
<10mm,1.5mm>*{};<8mm,7mm>*{}**@{-},
  <10mm,1.5mm>*{};<9mm,7mm>*{}**@{-},
  <10mm,1.5mm>*{};<10.5mm,6mm>*{}**@{-},
  <10mm,1.5mm>*{};<12mm,7mm>*{}**@{-},
 <10mm,1.5mm>*{};<10.5mm,6.6mm>*{.\hspace{-0.4mm}.\hspace{-0.4mm}.}**@{},
 <-10mm,-1.5mm>*{};<-12mm,-6mm>*{}**@{-},
 <-7mm,-1.5mm>*{};<-8mm,-6mm>*{}**@{-},
 <-4mm,-1.5mm>*{};<-4.5mm,-6mm>*{}**@{-},
 <0mm,0mm>*{};<0mm,-4.6mm>*{.\hspace{0.1mm}.\hspace{0.1mm}.}**@{},
<10mm,-1.5mm>*{};<12mm,-6mm>*{}**@{-},
 <7mm,-1.5mm>*{};<8mm,-6mm>*{}**@{-},
  <4mm,-1.5mm>*{};<4.5mm,-6mm>*{}**@{-},
<0mm,0mm>*{};<-9.5mm,8.2mm>*{^{I_{ 1}}}**@{},
<0mm,0mm>*{};<-3mm,8.2mm>*{^{I_{ i}}}**@{},
<0mm,0mm>*{};<2mm,8.2mm>*{^{I_{ i+1}}}**@{},
<0mm,0mm>*{};<10mm,8.2mm>*{^{I_{ k}}}**@{},
<0mm,0mm>*{};<-12mm,-7.4mm>*{_1}**@{},
<0mm,0mm>*{};<-8mm,-7.4mm>*{_2}**@{},
<0mm,0mm>*{};<-4mm,-7.4mm>*{_3}**@{},
<0mm,0mm>*{};<6mm,-7.4mm>*{\ldots}**@{},
<0mm,0mm>*{};<12.5mm,-7.4mm>*{_n}**@{},
\endxy
 =
 \sum_{i=1}^k(-1)^{i+1}
\xy
 <0mm,0mm>*{\mbox{$\xy *=<20mm,3mm>\txt{}*\frm{-}\endxy$}};<0mm,0mm>*{}**@{},
  <-10mm,1.5mm>*{};<-12mm,7mm>*{}**@{-},
  <-10mm,1.5mm>*{};<-11mm,7mm>*{}**@{-},
  <-10mm,1.5mm>*{};<-9.5mm,6mm>*{}**@{-},
  <-10mm,1.5mm>*{};<-8mm,7mm>*{}**@{-},
 <-10mm,1.5mm>*{};<-9.5mm,6.6mm>*{.\hspace{-0.4mm}.\hspace{-0.4mm}.}**@{},
 <0mm,0mm>*{};<-5.5mm,3.6mm>*{.\hspace{-0.1mm}.\hspace{-0.1mm}.}**@{},
%
  <0mm,1.5mm>*{};<-4mm,7mm>*{}**@{-},
  <0mm,1.5mm>*{};<-2.0mm,6mm>*{}**@{-},
  <0mm,1.5mm>*{};<-1mm,7mm>*{}**@{-},
 <0mm,1.5mm>*{};<-2.3mm,6.6mm>*{.\hspace{-0.4mm}.\hspace{-0.4mm}.}**@{},
%
  <0mm,1.5mm>*{};<1mm,7mm>*{}**@{-},
  <0mm,1.5mm>*{};<2.0mm,6mm>*{}**@{-},
  <0mm,1.5mm>*{};<4mm,7mm>*{}**@{-},
 <0mm,1.5mm>*{};<2.3mm,6.6mm>*{.\hspace{-0.4mm}.\hspace{-0.4mm}.}**@{},
 <0mm,0mm>*{};<6mm,3.6mm>*{.\hspace{-0.1mm}.\hspace{-0.1mm}.}**@{},
<10mm,1.5mm>*{};<8mm,7mm>*{}**@{-},
  <10mm,1.5mm>*{};<9mm,7mm>*{}**@{-},
  <10mm,1.5mm>*{};<10.5mm,6mm>*{}**@{-},
  <10mm,1.5mm>*{};<12mm,7mm>*{}**@{-},
 <10mm,1.5mm>*{};<10.5mm,6.6mm>*{.\hspace{-0.4mm}.\hspace{-0.4mm}.}**@{},
 <-10mm,-1.5mm>*{};<-12mm,-6mm>*{}**@{-},
 <-7mm,-1.5mm>*{};<-8mm,-6mm>*{}**@{-},
 <-4mm,-1.5mm>*{};<-4.5mm,-6mm>*{}**@{-},
 <0mm,0mm>*{};<0mm,-4.6mm>*{.\hspace{0.1mm}.\hspace{0.1mm}.}**@{},
<10mm,-1.5mm>*{};<12mm,-6mm>*{}**@{-},
 <7mm,-1.5mm>*{};<8mm,-6mm>*{}**@{-},
  <4mm,-1.5mm>*{};<4.5mm,-6mm>*{}**@{-},
<0mm,0mm>*{};<-9.5mm,8.2mm>*{^{I_{ 1}}}**@{},
<0mm,0mm>*{};<-2.5mm,8.2mm>*{{^{I_{ i}\sqcup}}}**@{},
<0mm,0mm>*{};<2.7mm,8.2mm>*{^{I_{ i+1}}}**@{},
<0mm,0mm>*{};<10mm,8.2mm>*{^{I_{ k}}}**@{},
<0mm,0mm>*{};<-12mm,-7.4mm>*{_1}**@{},
<0mm,0mm>*{};<-8mm,-7.4mm>*{_2}**@{},
<0mm,0mm>*{};<-4mm,-7.4mm>*{_3}**@{},
<0mm,0mm>*{};<6mm,-7.4mm>*{\ldots}**@{},
<0mm,0mm>*{};<12.5mm,-7.4mm>*{_n}**@{},
\endxy
 + \sum_{p+q=k+1\atop p\geq 1,q\geq 0}\sum_{i=0}^{p-1}
\sum_{  {I_{i+1}=I_{i+1}'\sqcup I''_{i+1}\atop .......................}
\atop
I_{i+q}=I_{i+q}'\sqcup I''_{i+q}}
$$
\vspace{-1mm}
\Beq\label{differential-DefQ}
\sum_{[n]=J_1\sqcup J_2 \atop {s\geq 0}}\hspace{-1mm}
(-1)^{(p+1)q + i(q-1)}
\frac{1}{s!}
\xy
 <19mm,0mm>*{\mbox{$\xy *=<58mm,3mm>\txt{}*\frm{-}\endxy$}};<0mm,0mm>*{}**@{},
  <-10mm,1.5mm>*{};<-12mm,7mm>*{}**@{-},
  <-10mm,1.5mm>*{};<-11mm,7mm>*{}**@{-},
  <-10mm,1.5mm>*{};<-9.5mm,6mm>*{}**@{-},
  <-10mm,1.5mm>*{};<-8mm,7mm>*{}**@{-},
 <-10mm,1.5mm>*{};<-9.5mm,6.6mm>*{.\hspace{-0.4mm}.\hspace{-0.4mm}.}**@{},
 <0mm,0mm>*{};<-6.5mm,3.6mm>*{.\hspace{-0.1mm}.\hspace{-0.1mm}.}**@{},
  <-3mm,1.5mm>*{};<-5mm,7mm>*{}**@{-},
  <-3mm,1.5mm>*{};<-4mm,7mm>*{}**@{-},
  <-3mm,1.5mm>*{};<-2.5mm,6mm>*{}**@{-},
  <-3mm,1.5mm>*{};<-1mm,7mm>*{}**@{-},
 <-3mm,1.5mm>*{};<-2.5mm,6.6mm>*{.\hspace{-0.4mm}.\hspace{-0.4mm}.}**@{},
  <10mm,1.5mm>*{};<0mm,7mm>*{}**@{-},
  <10mm,1.5mm>*{};<4mm,7mm>*{}**@{-},
<10mm,1.5mm>*{};<7.3mm,5.9mm>*{.\hspace{-0.0mm}.\hspace{-0.0mm}.}**@{},
   <10mm,1.5mm>*{};<3.8mm,6.0mm>*{}**@{-},
 <10mm,1.5mm>*{};<2.5mm,6.6mm>*{.\hspace{-0.4mm}.\hspace{-0.4mm}.}**@{},
%
%
<10mm,1.5mm>*{};<9mm,7mm>*{}**@{-},
  <10mm,1.5mm>*{};<10.5mm,6mm>*{}**@{-},
  <10mm,1.5mm>*{};<12mm,7mm>*{}**@{-},
 <10mm,1.5mm>*{};<10.5mm,6.6mm>*{.\hspace{-0.4mm}.\hspace{-0.4mm}.}**@{},
 <-10mm,-1.5mm>*{};<-12mm,-6mm>*{}**@{-},
 <-7mm,-1.5mm>*{};<-8mm,-6mm>*{}**@{-},
 <-4mm,-1.5mm>*{};<-4.5mm,-6mm>*{}**@{-},
 <-1mm,-1.5mm>*{};<-1.1mm,-6mm>*{}**@{-},
 <2mm,-1.5mm>*{};<2.0mm,-6mm>*{}**@{-},
 <0mm,0mm>*{};<20mm,-4.6mm>*{.\hspace{2mm}.\hspace{2mm}.\hspace{2mm}
 .\hspace{2mm}.\hspace{2mm}.\hspace{2mm}.\hspace{2mm}
 .\hspace{2mm}.\hspace{2mm}}**@{},
<48mm,-1.5mm>*{};<50mm,-6mm>*{}**@{-},
 <45mm,-1.5mm>*{};<46mm,-6mm>*{}**@{-},
  <42mm,-1.5mm>*{};<42.5mm,-6mm>*{}**@{-},
  <39mm,-1.5mm>*{};<39.2mm,-6mm>*{}**@{-},
   <36mm,-1.5mm>*{};<36mm,-6mm>*{}**@{-},
 <20mm,-1.5mm>*{};<20.0mm,-8mm>*{\underbrace{\hspace{66mm}}}**@{},
 <20mm,-1.5mm>*{};<20.0mm,-11mm>*{_{J_1}}**@{},
<0mm,0mm>*{};<-9.5mm,8.4mm>*{^{I_{ 1}}}**@{},
<0mm,0mm>*{};<-3mm,8.4mm>*{^{I_{ i}}}**@{},
<0mm,0mm>*{};<2mm,8.6mm>*{^{I_{ i+1}'}}**@{},
<0mm,0mm>*{};<10.5mm,8.6mm>*{^{I_{ i+q}'}}**@{},
<10mm,1.5mm>*{};<18mm,12mm>*{}**@{-},
<10mm,1.5mm>*{};<20.0mm,12mm>*{}**@{-},
<10mm,1.5mm>*{};<25mm,12mm>*{}**@{-},
<10mm,1.5mm>*{};<18.7mm,8.6mm>*{.\hspace{-0.4mm}.\hspace{-0.4mm}.}**@{},
<10mm,1.5mm>*{};<19.6mm,9.0mm>*{^s}**@{},
<25mm,13.75mm>*{\mbox{$\xy *=<14mm,3mm>\txt{}*\frm{-}\endxy$}};
<0mm,0mm>*{}**@{},
 <18mm,15mm>*{};<16mm,20.5mm>*{}**@{-},
 <18mm,15mm>*{};<17mm,20.5mm>*{}**@{-},
 <18mm,15mm>*{};<18.5mm,19.6mm>*{}**@{-},
 <18mm,15mm>*{};<20mm,20.5mm>*{}**@{-},
 <18mm,15mm>*{};<18.6mm,20.3mm>*{.\hspace{-0.4mm}.\hspace{-0.4mm}.}**@{},
<22mm,15mm>*{};<25.5mm,17.7mm>*{\cdots}**@{},
 <32mm,15.2mm>*{};<30mm,20.5mm>*{}**@{-},
 <32mm,15.2mm>*{};<31mm,20.5mm>*{}**@{-},
 <32mm,15.2mm>*{};<32.5mm,19.6mm>*{}**@{-},
 <32mm,15mm>*{};<34mm,20.5mm>*{}**@{-},
 <32mm,15mm>*{};<32.3mm,20.3mm>*{.\hspace{-0.4mm}.\hspace{-0.4mm}.}**@{},
%
<0mm,0mm>*{};<18mm,22.6mm>*{^{I_{ i+1}''}}**@{},
<0mm,0mm>*{};<32.5mm,22.6mm>*{^{I_{ i+q}''}}**@{},
 <26mm,12mm>*{};<25mm,9mm>*{}**@{-},
 <27mm,12mm>*{};<26.8mm,9mm>*{}**@{-},
 <29mm,12mm>*{};<29.2mm,10mm>*{.\hspace{-0.1mm}.\hspace{-0.1mm}.}**@{},
 <32mm,12mm>*{};<33mm,9mm>*{}**@{-},
  <31mm,12mm>*{};<31.5mm,9mm>*{}**@{-},
 <29mm,12mm>*{};<29mm,8mm>*{\underbrace{\ \ \ \ \ \ \ \  }}**@{},
 <29mm,12mm>*{};<29mm,5.3mm>*{_{J_2}}**@{},
<38mm,1.5mm>*{};<36mm,7mm>*{}**@{-},
<38mm,1.5mm>*{};<37mm,7mm>*{}**@{-},
<38mm,1.5mm>*{};<38.5mm,6mm>*{}**@{-},
<38mm,1.5mm>*{};<40mm,7mm>*{}**@{-},
<38mm,1.5mm>*{};<38.5mm,6.6mm>*{.\hspace{-0.4mm}.\hspace{-0.4mm}.}**@{},
<38mm,1.5mm>*{};<43mm,4mm>*{.\hspace{-0.0mm}.\hspace{-0.0mm}.}**@{},
<48mm,1.5mm>*{};<46mm,7mm>*{}**@{-},
<48mm,1.5mm>*{};<47mm,7mm>*{}**@{-},
<48mm,1.5mm>*{};<48.5mm,6mm>*{}**@{-},
<48mm,1.5mm>*{};<50mm,7mm>*{}**@{-},
<48mm,1.5mm>*{};<48.5mm,6.6mm>*{.\hspace{-0.4mm}.\hspace{-0.4mm}.}**@{},
<0mm,0mm>*{};<40.3mm,8.6mm>*{^{I_{ i+q+1}}}**@{},
<0mm,0mm>*{};<48.5mm,8.6mm>*{^{I_{ k}}}**@{},
\endxy
\Eeq
It was shown in \cite{Me1} that there is a one-to-one correspondence between degree $0$ representations
of the dg prop $(\DefQ, d)$ in a dg vector space $V$ and Maurer-Cartan elements in the Hochschild dg Lie algebra
$(\bar{C}^\bullet(\bar{\f}_V), [\ ,\ ]_G, d_H)$, i.e.\ with unital $A_\infty$-structures on $\f_V$.
Put another way, the dg Lie algebra induced on $\mbox{Rep}_V(\DefQ)[1]$ from the above differential $d$ is precisely
the Hochschild dg Lie algebra.

 Consider now a filtration, $F_{-p}:=\{ \mbox{span}\langle G\rangle: \mbox{number of vertices in}\ G\geq p\}$,
of the complex $(\DefQ, d)$.
It is clear that $0$-th term, $(E_0,\delta)$, of the associated spectral sequence, $\{E_r,d_r\}_{r\geq 1}$,
 isomorphic
(modulo an inessential shift of degree)
to the prop $(\caD, \delta)$ introduced in \S 3.3 so that, by Proposition~3.3.1, we conclude that
 $E_1=H(E_0)$ is isomorphic as a free
prop to $\PV$ whose shifted representation space, $\mbox{Rep}_V(\PV)[1]$,
 is $H(\bar{C}^\bullet(\bar{\f}_V))=\wedge^{\bullet\geq 1}\cT_V$. The  Lie algebra
structure on $H(\bar{C}^\bullet(\bar{\f}_V))$ induced from the Gerstenhaber brackets on $\bar{C}^\bullet(\bar{\f}_V)$
is then given by the differential, $d_1$, induced on the next term of the spectral sequence, $E_1=\PV$,
 from the differential $d$ in $\DefQ$. A direct inspection of formula (\ref{differential-DefQ})
implies that $d_1$ is precisely $\p$ which in turn implies by \S 4.2 that the induced Lie algebra structure
on  $H(\bar{C}^\bullet(\bar{\f}_V)$ is indeed given by Schouten brackets. It is worth noting in conclusion
that  $L_\infty$-morphisms (in the sense of Kontsevich \cite{Ko}) between dg Lie algebras $\bar{C}^\bullet(\bar{\f}_V)$
and $\wedge^{\bullet\geq 1}\cT_V$ can be equivalently understood as morphisms of dg props, $\DefQ\rar \PV^\circlearrowright$,
where $\PV^\circlearrowright$ is the wheeled completion of the prop of polyvector fields (by definition,
 $\PV^\circlearrowright$
is the smallest wheeled prop containing $\PV$ as a subspace). This point of view on quantizations was discussed
in more detail in \cite{Me1,Me3}.

%
%
%
\input{referenc}



\printindex
\end{document}

%% file: referenc.tex
%
%
%

%
%

%% file: sm.bbl
\begin{thebibliography}{[KLR73]}
%
%
%







\bibitem[Ba]{Ba} Baranovsky, V.:
\newblock { A universal enveloping for $L_\infty$-algebras.}. Preprint
\newblock{arXiv:0706.1396} (2006).

\bibitem[BM]{BM} Borisov, D.,  Manin, Yu.I.:
{Internal cohomomorphisms for operads}, arXiv:math.CT/0609748 (2006).


\bibitem[CFL]{CFL}
Cattaneo, A.S., Fiorenza, D.,  Longoni, R.:
{On the Hochschild-Kostant-Rosenberg map for graded manifolds},
arXiv:math/0503380 (2006).

\bibitem[Ch]{Ch} Chapoton, F.:
\newblock {Op\'{e}rades diff\'{e}rentielles gradu\'{e}es sur les simplexes et les permuto\`{e}dres }. Preprint
\newblock{arXiv:math/0102172} (2001).






\bibitem[GS1]{GS} Gerstenhaber, M., Schack, S.D.:
\newblock {Bialgebra cohomology, deformations, and quantum groups}.
\newblock{Proc.\ Nat.\ Acad.\ Sci.\ USA, } {\bf 87}, 478-481 (1990).

\bibitem[GS2]{GS2} Gerstenhaber, M., Schack, S.D.:
\newblock {Algebras, bialgebras, quantum groups, and algebraic
deformations}. In:
\newblock {Contemporary Mathematics} {\bf 134}, 51-92 (1992).


\bibitem[GK]{GeKa}
Getzler E., Kapranov, M.M.:
\newblock {Modular operads}.
\newblock { Compositio Math.}, {\bf 110}, 65-126 (1998).



\bibitem[EK]{EK} Etingof, P., Kazhdan, D.:
\newblock { Quantization of Lie bialgebras, I}.
\newblock{ Selecta Math. (N.S.)} {\bf 2} , 1-41 (1996)


\bibitem[Ko]{Ko} Kontsevich, M.:
\newblock {Deformation quantization
 of Poisson manifolds I},
 math/9709040.

\bibitem[LM]{LM} Lazarev, A.Yu., Movshev, M.V.:
\newblock {Deformations of Hopf algebras},
\newblock{Uspekhi Mat.\ Nauk}, {\bf 46}, 211-212 (1991)

\bibitem[Lo]{Lo} Loday, J.L.: { Cyclic homology}, Springer (1998)



\bibitem[MSS]{MSS}
Markl, M., Shnider,S., Stasheff, J.D.:
\newblock {Operads in Algebra, Topology and Physics}.
\newblock American Mathematical Society, Providence, Rhode Island (2002)




\bibitem[Me1]{Me0} Merkulov, S.A.:
\newblock { Prop profile of Poisson geometry},
\newblock math.DG/0401034, { Commun. Math. Phys.} {\bf 262}, 117-135 (2006)

\bibitem[Me2]{Me1} Merkulov, S.A.:
\newblock {Graph complexes with
  loops and wheels}. To appear in: Algebra, Arithmetic and Geometry - Manin Festschrift (2007)



\bibitem[Me3]{Me3} Merkulov, S.A.: { Lectures on props, Poisson geometry and deformation quantization},
to appear in {Poisson Geometry in Mathematics and Physics}, Contemporary Mathematics
(eds.\ G.\ Dito, J.H.\ Lu, Y.\ Maeda and A.\ Weinstein), AMS (2007).

\bibitem[Me4]{Me4} Merkulov, S.A.: {Quantization of strongly homotopy Lie bialgebras},
arXiv:math/0612431  (2006)



\bibitem[MV]{MV} Merkulov, S.A., Vallette, B.:
\newblock{Deformation theory of representations of prop(erad)s},
\newblock  preprint arXiv:0707.0889 (2007).


\bibitem[St]{St} Stasheff, J.D.: {On the homotopy
  associativity of $H$-spaces, I \@ II}, Trans.\ Amer.\
  Math.\ Soc.\ {\bf 108}, 272-292 \& 293-312 (1963)

\bibitem[SU]{SU}
Saneblidze, S.,  Umble, R.: {Diagonals on the permutahedra, multiplihedra and
associahedra}. Homology Homotopy Appl. {\bf 6}, 363-411 (2004).


\bibitem[Va]{Va} Vallette, B.:  {A Koszul duality for
props}.  Trans.\ AMS., {\bf 359}, 4865-4943 (2007)








\end{thebibliography}
